%% file: main.tex
\pdfoutput=1
\PassOptionsToPackage{hyphens}{url}
\documentclass[11pt]{article}

\usepackage[margin=1in]{geometry}
\usepackage[T1]{fontenc}
\usepackage[utf8]{inputenc}
\usepackage{lmodern}
\usepackage{microtype}
\usepackage{amsmath,amssymb,amsthm}
\usepackage{mathtools}
\usepackage{booktabs}
\usepackage{colortbl}
\usepackage{array}
\usepackage{graphicx}
\usepackage{xcolor}
\usepackage[numbers,sort&compress]{natbib}
\usepackage[hidelinks]{hyperref}
\usepackage[capitalise,noabbrev]{cleveref}

\newcommand{\edit}[1]{\textsf{#1}}

\usepackage{longtable}
\usepackage{algorithm}
\usepackage[noend]{algpseudocode}
\usepackage{fancyvrb}
\newenvironment{kind}[1]{%
  \def\kindshape{#1}%
  \par\medskip\noindent
  \begin{minipage}[t]{0.6\linewidth}\vspace{0pt}}{%
  \end{minipage}\hfill
  \begin{minipage}[t]{0.36\linewidth}\vspace{0pt}\centering
    \includegraphics[width=\linewidth]{figures/diagrams/pdf/\kindshape}%
  \end{minipage}\par\medskip}
\usepackage{flafter}

\AddToHook{cmd/subsection/before}{\FloatBarrier}
\AddToHook{cmd/subsubsection/before}{\FloatBarrier}
\usepackage[section]{placeins}

\title{Reclassified to Conform}
\author{Chris Martin\\
  \normalsize Kepler AI\\
  \normalsize\texttt{chris.martin@kepler.ai}}
\date{September 22, 2026}

\begin{document}
\maketitle

\begin{abstract}
A company's financial statements over many years must be assembled from
many filings. Each filing may rename, combine, split or move its lines,
and restates only the last two or three years to match. Filings provide
no consistent, machine-readable mapping between successive statements,
and neither a line's label nor its XBRL tag is fixed. We treat the reconstruction as an
optimization problem. A tag that keeps its numbers from one filing to
the next is taken to be the same line, and such runs are joined into
strands before any edit is weighed. At
each boundary between filings, the candidate edits that pass the check
their kind requires are listed, each with the strand endpoints it would
explain and a price for how much it asks the reader to
believe beyond what the filings show. An integer program picks the
cheapest set of the listed edits that accounts for every line exactly
once, and the chosen edits trace each line through the filings. The
result is a statement whose rows follow one line across every year,
comparable within each filing's restated window and marked where a
number was restated or derived, in any
filing's layout, with each number cited to the filing that stated it
and each declared total re-checked. On 200 companies drawn at random,
99.9\% of the totals that can be checked add up. Among the 124 with a
usable vendor comparison, the lines the method follows agree with the
vendor's series on 96\% of the cells both fill, reach 10\% more cells
than a tag alone and agree with the vendor on 90\% of those. The method also yields a
record of what each company changed in its reporting.
\end{abstract}

\vspace{1.5em}
\begin{flushright}
  \begin{minipage}{0.46\linewidth}
    \raggedleft\small
    \itshape When I use a word, it means just what I choose it to mean,
    neither more nor less.\\[0.4em]
    \upshape --- Humpty Dumpty, in \emph{Through the Looking-Glass}
  \end{minipage}
\end{flushright}
\vspace{1.5em}

\input{sections/01-problem}
\input{sections/02-guarantees}
\input{sections/03-formulation}
\input{sections/04-optimization}
\input{sections/05-results}
\input{sections/06-next}
\input{sections/07-coda}

\clearpage
\appendix
\begin{center}
  {\LARGE\bfseries Appendices}
\end{center}
\vspace{1em}
\input{sections/a-model}
\input{sections/b-vocabulary}
\input{sections/c-examples}
\input{sections/d-hard-cases}
\input{sections/e-measurements}

\bibliographystyle{plainnat}
\bibliography{references}

\end{document}

%% file: sections/01-problem.tex
\section{The problem}\label{sec:problem}

Company filings are the case this paper takes up. Anyone analyzing a
company wants to see how it has fared over time. A snapshot is not
enough; a forecast or a model needs the course of the business across
many years, and so many years of reporting.

Reported data is plentiful: public companies file their financial
statements with the SEC every quarter and every year, every number
tagged and machine-readable. The measures are standardized by GAAP and
audited, so each number is reliable for the period it was filed. But
which lines a company reports, and how finely it breaks them down, is
largely its own decision~\cite{reg-s-x,asc205,edgar-xbrl-guide}, and
those decisions change.

\begin{figure}[htbp]
  \small
  \newcommand{\cutgap}{\\[3pt]{\color{gray!60}\dotfill}\\[3pt]}
  \begingroup\offinterlineskip\centering
  \textsf{\textbf{(a)} 10-K for fiscal 2021, cash flow statement}\\[4pt]
  \includegraphics[width=\linewidth]{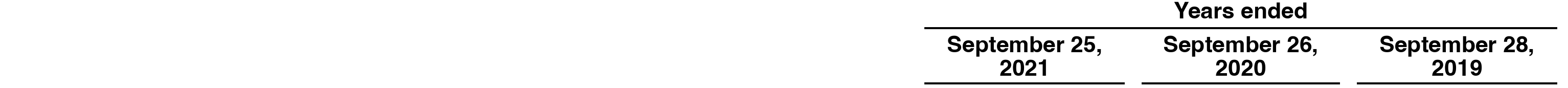}\cutgap
  \includegraphics[width=\linewidth]{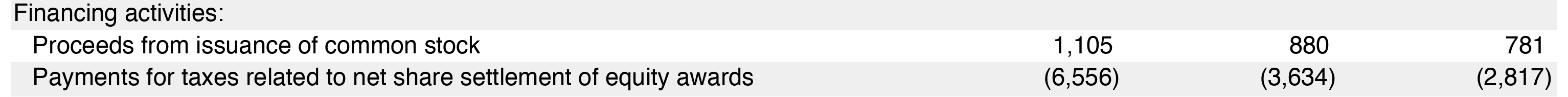}\cutgap
  \includegraphics[width=\linewidth]{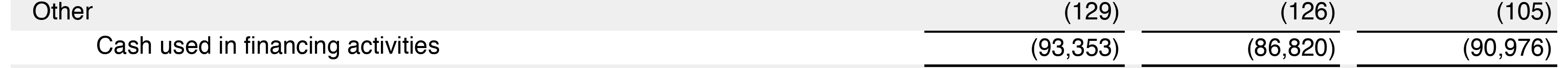}\\[10pt]
  \textsf{\textbf{(b)} 10-K for fiscal 2022, cash flow statement}\\[4pt]
  \includegraphics[width=\linewidth]{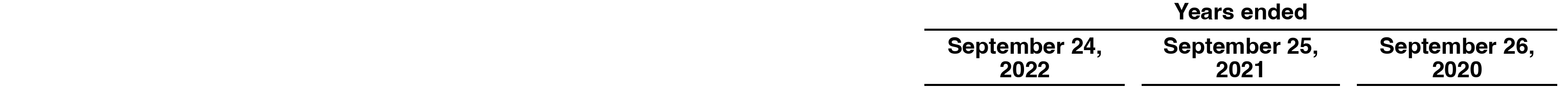}\cutgap
  \includegraphics[width=\linewidth]{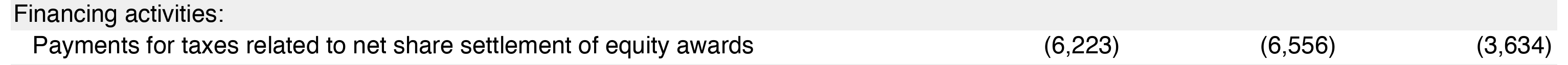}\cutgap
  \includegraphics[width=\linewidth]{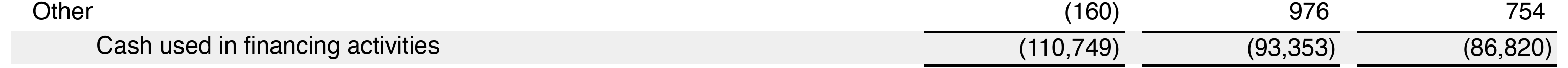}\par
  \endgroup
  \caption{The financing section of Apple's cash flow statement in two
    consecutive annual reports, cut from the filings; the dotted lines
    mark omitted rows. Each report shows three years; the two share 2021
    and 2020, and 2019 appears only in the
    first~\cite{aapl-10k-2021,aapl-10k-2022}.}
  \label{fig:window}
\end{figure}

Apple's cash flow statement is the small case that runs through the
paper (\cref{fig:window}). Its report for fiscal 2021 shows
\emph{proceeds from issuance of common stock} of 880 for 2020. The
report for fiscal 2022 has no such line, and \emph{other} for 2020 is
754, the old 880 plus the old
(126)~\cite{aapl-10k-2021,aapl-10k-2022}. Both numbers were
true when filed. This is the easiest case: one line moved, and the
arithmetic says where. Usually several lines move at once and more than
one explanation fits, and sometimes the company gets it wrong
(\cref{sec:optimization}, Appendix~\ref{sec:solutions}). This paper
states the question precisely and answers it for every line of every
statement.

The figure also shows what makes the question answerable. Two filings
a year apart differ because the presentation moved and because the
business moved, and comparing their newest years cannot tell one from
the other. But each filing shows the prior years again: the fiscal
2022 report presents 2020 under its own lines, and the fiscal 2021
report presented the same 2020 under its. With one period under two
presentations, the change in presentation can be read apart from the
change in the business. That is the evidence the paper works from. It
constrains the answer but does not always fix it: several potential
sets of edits can fit the same numbers, and \cref{sec:example} shows
a boundary that admits three.

The paper offers a way to discover what changed and a method for
rebuilding comparable history where it exists. Its title is the
sentence companies print whenever they change their
reporting~\cite{farmhouse-10k-2024}:
\begin{quote}
  \itshape Certain amounts in prior periods have been reclassified to
  conform to the current period presentation.
\end{quote}

\subsection{The mask}\label{sec:mask}

A company has financial properties that outsiders never see directly. A
filing shows what the company chooses to show: each line is a total
over things the reader cannot see. Call that choice the \emph{mask}.
The mask is what changes when a company changes its reporting; what we
want is the properties underneath, over time. Names, labels and
position on the page say something about how two masks differ, but
not much. What separates them is the re-presented years: this year's
filing applies this year's mask to last year's state, and last year's
filing applied last year's mask to the same state.

The mask can move in many ways. A few recur, and
each is drawn in Appendix~\ref{app:vocabulary}. Three change the
detail. A line continues under a new name or tag, as when Netflix's
deferred revenue moved to the tag of a new accounting
standard~\cite{nflx-10k-2017,nflx-10k-2018}. One line becomes several,
as when JPMorgan separated asset-management fees from
commissions~\cite{jpm-10k-2022,jpm-10k-2023}. Several lines become one,
as when Apple's stock proceeds disappeared into
\emph{other}~\cite{aapl-10k-2021,aapl-10k-2022}.

Two change the numbers. Money moves between lines: a line keeps its
name, part of what it held now sits elsewhere, and the prior years are
restated to match, as when Tesla took leasing out of automotive
revenue~\cite{tsla-10k-2015,tsla-10k-2016}. Or the company gets it
wrong: Palantir attached three years of numbers to the wrong
lines~\cite{pltr-10k-2020,pltr-10k-2021,pltr-10k-2023}. And a number can
be revised with nothing about the layout changing.

Which lines exist is not enough to tell these apart. A line gone from
the newer filing may have ended, gone into a neighbour, been combined
with others into a new line, or been rolled into its total, and all
four leave the newer filing looking the same. Only the numbers on the
year both filings report tell them apart. \Cref{sec:formulation}
turns this into a precise list of the kinds of edit, and draws the four
(\cref{fig:disappear-named}).

\subsection{Analogues in other fields}\label{sec:related}

The problem is not specific to accounting. Wherever something is
observed only through a series of reports, each under a classification
the reporter is free to change, anyone who wants the thing itself over
time has to reconcile the reports first. Every organization that
maintains a classification builds the same three things: a list of the
kinds of change, a record of which old codes correspond to which new
ones, and a warning that combining cannot be undone. The Census Bureau
publishes relationship files whenever tracts are redrawn, typed and
weighted~\cite{census-relationship-files}, and the mapping from ICD-9
to ICD-10 medical codes is published in both directions and documented
as lossy~\cite{cms-gems}. Economists keep every published version of a
statistical series, so that what was known on a given date can be
recovered~\cite{alfred}; that answers what was known and when, and
reconciles nothing, since a vintage archive keeps every label as it
was.\footnote{Also: git guesses that a file was renamed from how
similar its contents are~\cite{git-diff-rename}, and Mercurial records
the rename~\cite{osullivan2009mercurial}; country codes are never
reused once retired, after \texttt{CS} was assigned
twice~\cite{iso3166-1}; trade economists group product codes into
families that stay comparable across
revisions~\cite{pierce2012concordance}; and finance has a permanent
identifier for a company~\cite{crsp-permno} and nothing like it for a
line of its statements.} The difference from company filings is who
controls the classification. The Census Bureau and the code-set
maintainers publish their mappings because they made the changes. In
company filings there is no such authority. Thousands of companies
each change their own list of lines on their own schedule, and what
they share is a vocabulary of tags, not a list of lines. A company
knows how its lines changed, since it restates the prior years to
conform, but its filings state that only as a sentence in a note, if
at all, and no one outside the company has it. So the mapping has to
be recovered from the filings themselves.

%% file: sections/02-guarantees.tex
\section{The filings as data}\label{sec:provide}

\subsection{The structure of a filing}\label{sec:xbrl}

A company files its financial statements in two forms: the document
people read, the 10-K or 10-Q, and XBRL, the eXtensible Business
Reporting Language, in which every number is a \emph{fact} a machine
can read~\cite{xbrl21}. A fact is a number plus what it measures (a
tag, from the shared \texttt{us-gaap} dictionary or the company's own),
the period it covers, its unit, and how it was rounded. A breakdown by
product or segment is a fact on the same tag, marked with the product
or segment, not a tag of its own. \Cref{fig:xbrl} shows two facts as
filed. A line of a statement is one tag's facts, one for each period
the statement shows.

\begin{figure}[htbp]
  \small
  \textsf{\textbf{(a)} Apple, 10-K for fiscal 2021: the number printed as ``880'',
  the context it points to, and its unit}
\begin{Verbatim}
<ix:nonFraction unitRef="usd"
                contextRef="i725ce765160e48a2a3b36989bda2ca94_D20190929-20200926"
                decimals="-6"
                name="us-gaap:ProceedsFromIssuanceOfCommonStock"
                scale="6"
                id="id3VybDovL2…3642a79d-e216-4edd-aaef-8b22f84f7992">880</ix:nonFraction>

<xbrli:context id="i725ce765160e48a2a3b36989bda2ca94_D20190929-20200926">
  <xbrli:entity>
    <xbrli:identifier scheme="http://www.sec.gov/CIK">0000320193</xbrli:identifier>
  </xbrli:entity>
  <xbrli:period>
    <xbrli:startDate>2019-09-29</xbrli:startDate>
    <xbrli:endDate>2020-09-26</xbrli:endDate>
  </xbrli:period>
</xbrli:context>

<xbrli:unit id="usd">
  <xbrli:measure>iso4217:USD</xbrli:measure>
</xbrli:unit>
\end{Verbatim}
  \medskip
  \textsf{\textbf{(b)} Tesla, 10-K for fiscal 2021: automotive leasing revenue, a
  revenue tag whose context names a product}
\begin{Verbatim}
<ix:nonFraction id="F_abccb4fb-2668-4799-8607-ec8ce640f676"
                contextRef="C_bae6b89d-012c-4639-b3ad-49711e152c22"
                name="us-gaap:RevenueFromContractWithCustomerExcludingAssessedTax"
                unitRef="U_USD"
                scale="6"
                decimals="-6"
                format="ixt:numdotdecimal">1,052</ix:nonFraction>

<xbrli:context id="C_bae6b89d-012c-4639-b3ad-49711e152c22">
  <xbrli:entity>
    <xbrli:identifier scheme="http://www.sec.gov/CIK">0001318605</xbrli:identifier>
    <xbrli:segment>
      <xbrldi:explicitMember dimension="srt:ProductOrServiceAxis">
        tsla:AutomotiveLeasingMember
      </xbrldi:explicitMember>
    </xbrli:segment>
  </xbrli:entity>
  <xbrli:period>
    <xbrli:startDate>2020-01-01</xbrli:startDate>
    <xbrli:endDate>2020-12-31</xbrli:endDate>
  </xbrli:period>
</xbrli:context>
\end{Verbatim}
  \caption{XBRL facts as filed. A fact is a number with a tag (\texttt{name}),
    a context (whose number, and for which period), a unit, and its rounding
    (\texttt{decimals}; \texttt{scale="6"} means the page shows millions). A
    breakdown is the same kind of fact, with the product or segment named
    in its context~\cite{aapl-10k-2021,tsla-10k-2021}.}
  \label{fig:xbrl}
\end{figure}

\subsection{The parts of a filing and their reliability}\label{sec:parts}

Four parts of a filing matter here: the facts, the layout, the
arithmetic, and the prior periods the filing shows again. The facts
give the values. The other three tell us different things and can be
relied on to different degrees.

The layout gives the order of the page and nothing more. A heading is
an abstract element, which the XBRL specification bars from carrying a
fact~\cite{xbrl21}, and none does; indentation says where a
line sits, not what it adds up to.

The arithmetic says which lines add up to which totals, and with what
sign, stated once per tag in a separate file of the filing. A filing
must declare it wherever a
statement shows lines and their total~\cite{edgar-xbrl-guide}, but a
total that does not add up has no effect on whether EDGAR accepts the
filing~\cite{edgar-xbrl-guide,efm-vol2}, so how often the equations
hold is a measurement. Over the 200 companies drawn at random for
\cref{sec:results} and their 2,826 annual reports, 87.7\% of the
104,408 declared totals add up exactly as declared, 89.8\% once
breakdown rows are folded into their line and the filer's rounding is
allowed, and 98.7\% among totals whose every part carries a number;
most of the rest have a declared part with no number in that period
(Appendix~\ref{app:footing}). Breakdowns are the exception: a
declared sum relates one tag to another in the same context, so a
line's breakdown is never declared to add up to it~\cite{xbrl21}, and a
breakdown can confirm a number but cannot prove one wrong.

The prior periods tell us the most. When a filing shows last year
again, it says that these are last year's numbers under this year's
lines: the two masks over one state that \cref{sec:mask} rests on. The window
is short: two balance sheets and three years of income and cash
flows~\cite{reg-s-x-periods}. Older years are never redone, so they
have no line in the newest filing. The history we want is ten or
fifteen years, so no single filing can supply it; each boundary between
consecutive filings has to be read on its own, from the two masks that
pair shows, and a line's history is the chain of those readings.
\Cref{fig:window} in \cref{sec:problem} shows the window as printed:
two consecutive Apple 10-Ks, each with three years, sharing two.

That is all that carries from one filing to the next. No rule requires
a filing to map its lines onto the one before. The accounting standard
requires only a note saying that prior periods were reclassified to
conform, with the nature and size of the change, in
words~\cite{asc205}. Nothing in a filing states, in a consistent,
machine-readable form, which of its lines continue an earlier filing's
or what became of a line that is gone; the note, where there is one,
says so in words. The equations and the re-presented years are the only
commitments a company makes across filings, and the rest of the paper
builds on those two alone.

%% file: sections/03-formulation.tex
\section{Formulating the problem}\label{sec:formulation}

For every measure a line of a company's statements represents, we want
its history: which line it was in each earlier filing, which of its
numbers were restated, and where the comparison breaks. The filings
do not state it, so it has to be inferred, and usually more than one
explanation of a boundary fits the numbers. That is the shape of an
optimization problem: an edit involves only a few lines, so the
candidates can be listed; a solution is a set of edits that accounts
for each line of both filings exactly once, an explanation of the
boundary; and among the explanations we want the one that asks the
least of the reader beyond what the filings show.

\subsection{Input and output}\label{sec:have}

The input is one company's filings of one statement, in the order
published. Each presents the statement as a list of lines; a line has a
tag, a label, a place in the layout, and a number for the current year
and one or two prior years. Each filing also declares which lines add
up to which totals. So every number has two dates, the period it
describes and the filing that stated it, and most years are stated by
two or three filings that do not always agree. Notes may describe how
lines changed, but provide no consistent, machine-readable mapping
between filings.

The output is the statement's history: between each pair of
consecutive filings, a set of edits that turns the older list of lines
into the newer one, with every line of both accounted for exactly once.
A line continues, under its old name or a new one; or it takes part in
an edit with other lines; or it ended or started. Nothing is left over
and nothing is used twice. The rest follows: rows that run across
years, the filing behind each number, where a comparison breaks, and
the table as first or as last reported, which differ only in which
filing each year's number is taken from.

\subsection{The kinds of edit}\label{sec:kinds}

\begin{figure}[htbp]
  \centering
  \includegraphics[width=\linewidth,height=0.8\textheight,keepaspectratio]{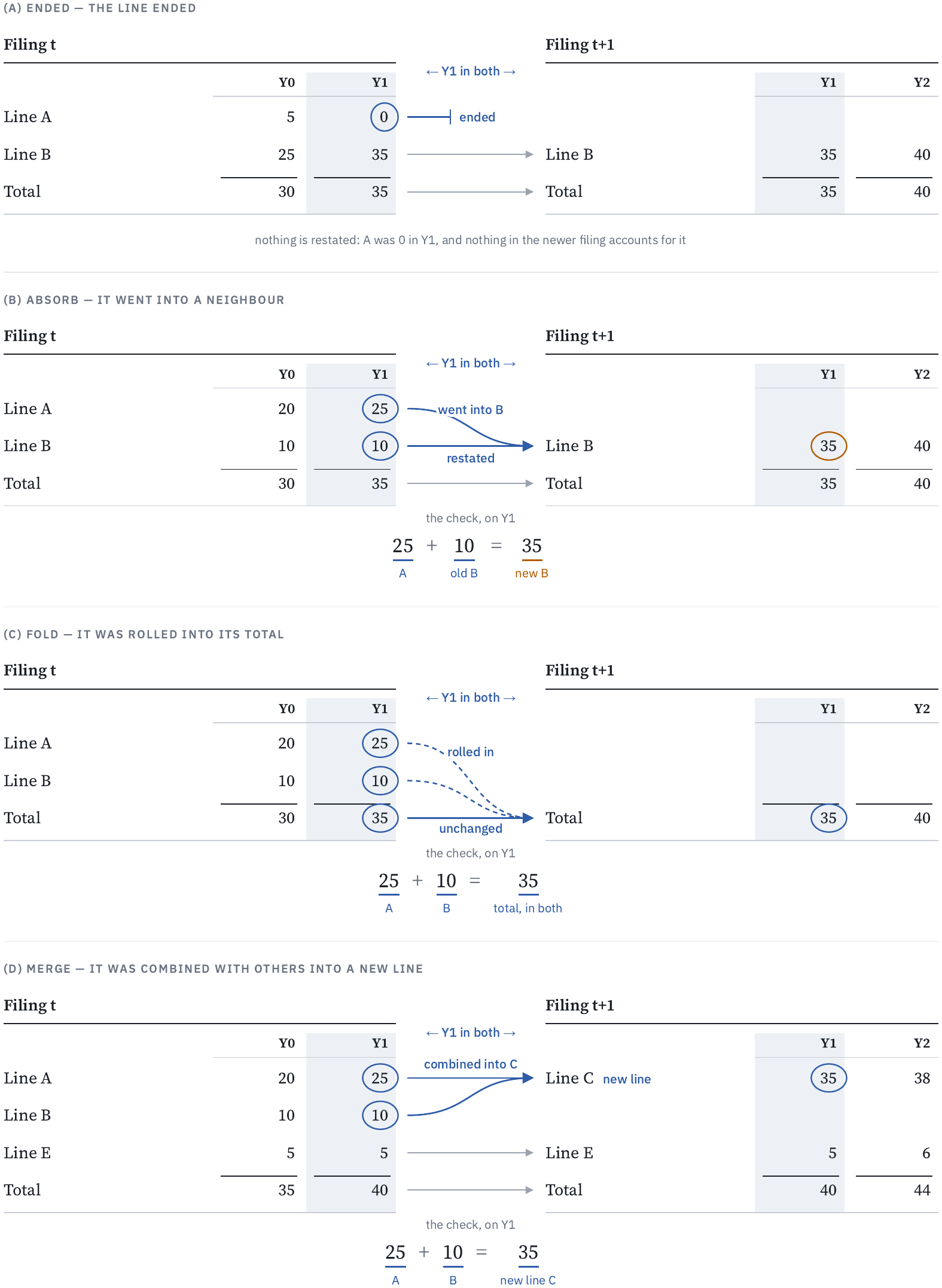}
  \caption{Four ways a line disappears, each named as a kind of edit
    from \cref{tab:edits}. In every case line A is gone from the newer
    filing; only the shared year, Y1, says which happened: nothing
    restated (\edit{ended}); a neighbour restated up by exactly the
    line's amount (\edit{absorb}); the total unchanged and the parts gone
    (\edit{fold}); a new line equal to the old ones' sum (\edit{merge}).
    Each side is one filing's statement, older on the left; the shaded
    column is the year both report; the edit is drawn in blue (dashed
    where parts are shown or hidden and no number moves); a restated
    number is circled in orange.}
  \label{fig:disappear-named}
\end{figure}

An edit says which lines of the older filing correspond to which of
the newer. Each kind is defined by what the numbers must show on the
years both filings report, and each moves the mask of \cref{sec:mask}
in its own way (\cref{tab:edits}). Three kinds add detail: a
\edit{split}, a \edit{carve-out} and a \edit{break out}. Three remove
it, each undoing one of the first three: a \edit{merge} undoes a split,
an \edit{absorb} a carve-out, a \edit{fold} a break out. A carve-out
and an absorb move money between lines, so a continuing line is
restated; the other four leave every continuing number as it was. A
\edit{rename} moves nothing. \edit{Equal copies} and a \edit{swap} are
rarer, and a line nothing explains has \edit{ended} or \edit{started}.
\Cref{fig:disappear-named} draws the four ways a line can be gone from
the newer filing that \cref{sec:problem} listed, each named: a line that
ended, an absorb, a fold and a merge.

\begin{table}[htbp]
  \centering
  \includegraphics[width=\linewidth]{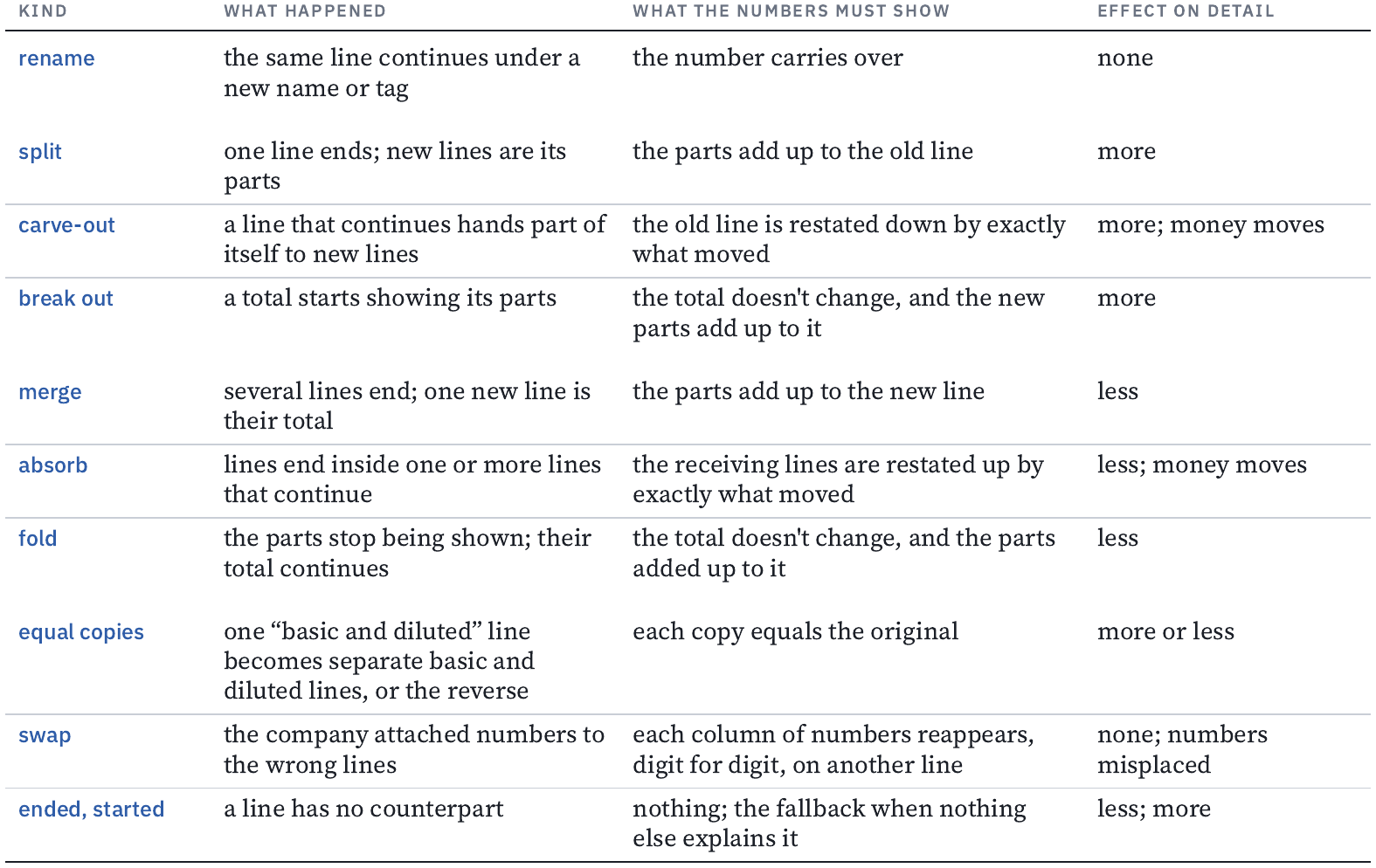}
  \caption{The kinds of edit. Two more touch numbers only: a
    \edit{restatement}, a prior year's number revised, and a
    \edit{transfer}, part of a split line moving to another section while
    that section's total is restated by the same amount. Every kind is
    drawn, with a made-up and a real example, in
    Appendix~\ref{app:vocabulary}.}
  \label{tab:edits}
\end{table}

\subsection{Multiple explanations}\label{sec:hard}

Several explanations can fit the same numbers; Tesla's 2016 income
statement admits at least three (\cref{sec:example}). Whether a line
is accounted for depends on every other decision at the boundary, so
the explanations cannot be settled line by line, and
\cref{sec:optimization} develops the price used to choose among them.

%% file: sections/04-optimization.tex
\section{Selecting the explanation}\label{sec:optimization}

\subsection{Choosing edits together}\label{sec:rules}

Candidate edits compete for the same lines, so selecting them is a
weighted set-packing problem. Although this problem is hard in general,
the instances here decompose into small connected components and are
inexpensive to solve exactly. We therefore use an integer program,
avoiding a problem-specific ordering heuristic. Across the random sample,
all component solves together took a quarter of a second
(\cref{sec:choice}).
Think of a market where edits compete to explain the lines. Each bids a
price, no line can be sold twice, and the market clears at the lowest total
cost, including the cost of lines left unexplained.

\subsection{A simple example}\label{sec:example}

\begin{figure}[htbp]
  \centering
  \includegraphics[width=\linewidth]{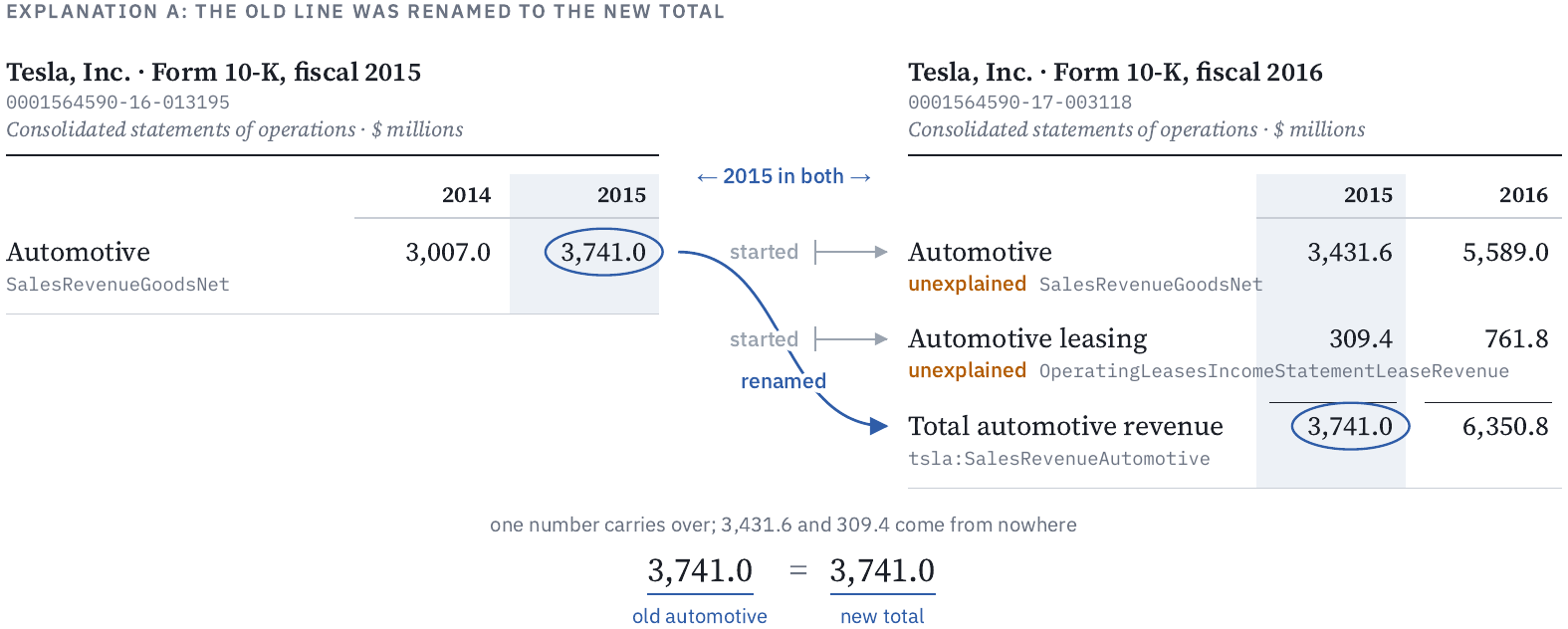}
  \caption{Explanation A, a rename alone.}
  \label{fig:tesla-reading-a}
\end{figure}

\begin{figure}[htbp]
  \centering
  \includegraphics[width=\linewidth]{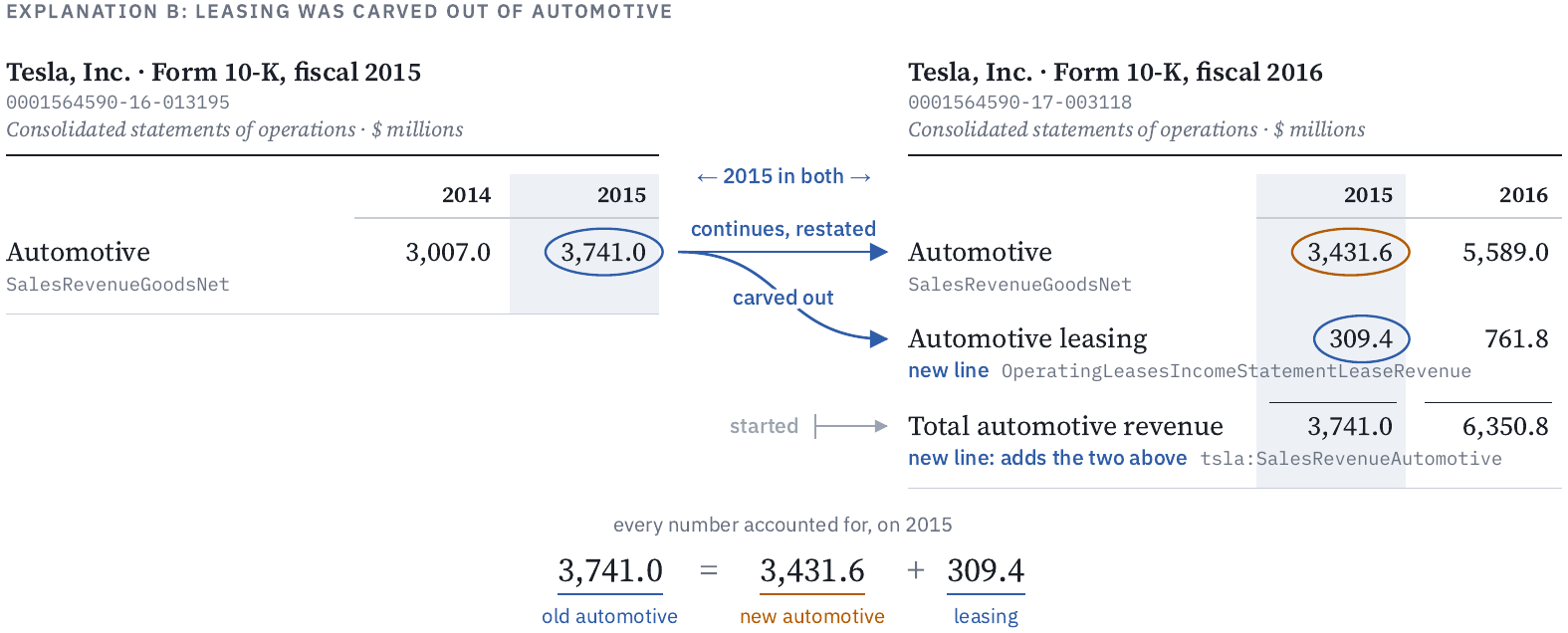}
  \caption{Explanation B, a carve-out.}
  \label{fig:tesla-reading-b}
\end{figure}

\begin{figure}[htbp]
  \centering
  \includegraphics[width=\linewidth]{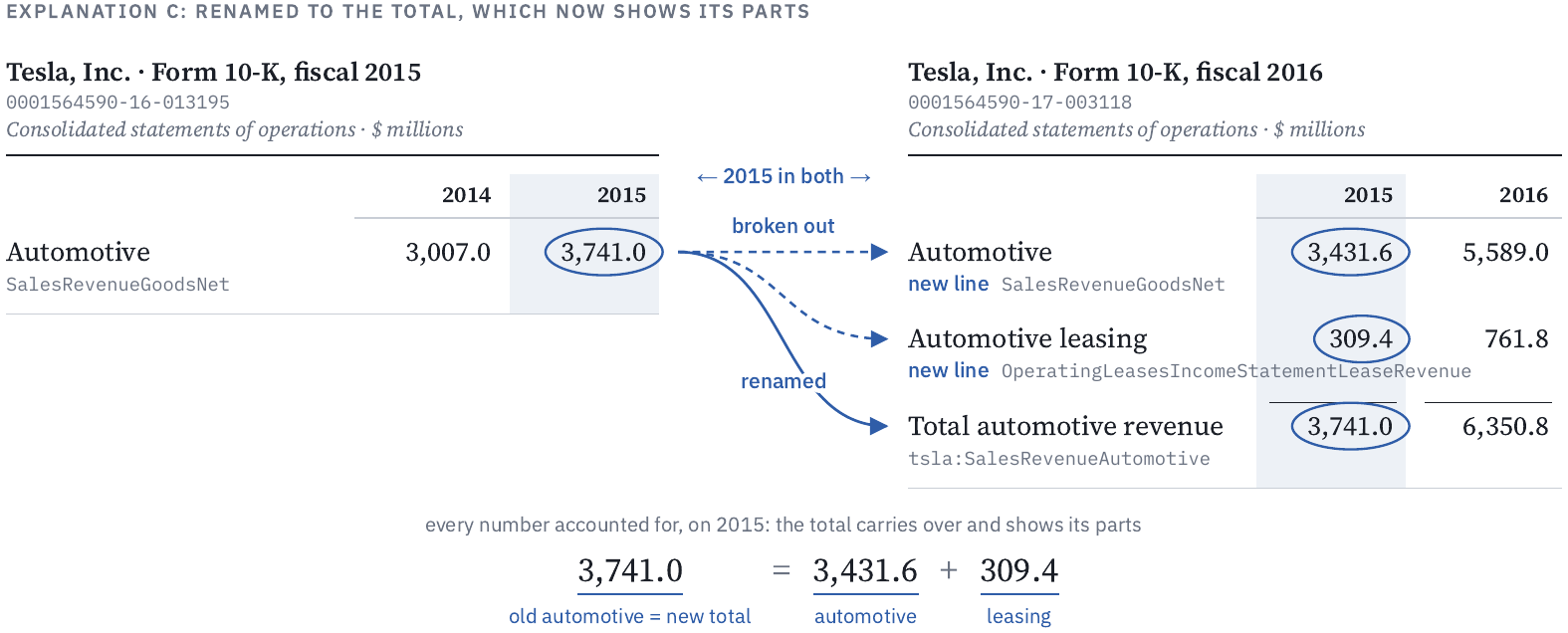}
  \caption{Explanation C, a rename plus a break out.}
  \label{fig:tesla-reading-c}
\end{figure}

A small case, to show the kinds of edit that can explain one boundary
and how they compete for the same lines. Tesla's 2015 and 2016 annual
reports are the case this section and Appendix~\ref{sec:solutions}
return to. The 2015 report shows Automotive revenue of 3,741.0. The 2016 report shows,
for the same year, Automotive 3,431.6 under the same tag, a new
Automotive leasing line of 309.4, and a new Total automotive revenue of
3,741.0 under a new tag~\cite{tsla-10k-2015,tsla-10k-2016}.
Three explanations fit, each a set of edits accounting for every line
of both reports
(\cref{fig:tesla-reading-a,fig:tesla-reading-b,fig:tesla-reading-c}).

\textbf{A, a rename alone.} The old line became ``Total automotive
revenue''; its number carries over exactly. The 2016 Automotive line
then keeps the old tag but has to be a new line, and neither its
3,431.6 nor the 309.4 is explained by anything.

\textbf{B, a carve-out.} Leasing was taken out of Automotive, which
continues restated ($3{,}741.0 = 3{,}431.6 + 309.4$), and the total is a
new line that adds the two up.

\textbf{C, a rename plus a break out.} The old line became the total,
with its number carried over, and the total now shows its parts: two new
lines that add up to it.

B and C both account for every number and differ in what they ask the
reader to believe: B, that Automotive was restated and a total
appeared; C, that the line changed tag and label and two detail lines
appeared. They leave different histories. Under C the total runs
1,921.9, 3,007.0, 3,741.0, 6,350.8 like for
like~\cite{tsla-10k-2015,tsla-10k-2016}, because Tesla's
older filings counted leasing inside automotive revenue; under B the
Automotive row mixes older years that include leasing with restated
years that do not.

\subsection{The procedure}\label{sec:how}

The procedure has four steps (Appendix~\ref{app:algorithm} states
them as pseudocode).
\begin{enumerate}
  \item \textbf{Strands.} Join, before anything is weighed, the same tag
    in consecutive filings with the same numbers on the shared years. Each
    run is a strand; where a strand ends or starts is an \emph{endpoint},
    and endpoints are what is left to explain (\cref{sec:strands}).
  \item \textbf{Candidates.} At each boundary, list the candidate edits
    of every kind: each could explain some of the endpoints and passes
    the check its kind requires, and is listed with the endpoints it
    would cover and a price (\cref{sec:candidates}).
  \item \textbf{Choose the edits.} Find the cheapest set of candidates that
    covers no endpoint twice, counting an endpoint no candidate covers at
    the price of leaving it unexplained (\cref{sec:choice}). Candidates
    that share no endpoint never interact, so the choice splits into
    independent problems, one per connected group of candidates, each
    solved exactly.
  \item \textbf{Rows.} Read the history off the choice: strands joined by
    the chosen renames are rows, and every year's numbers come from one
    filing (\cref{sec:output}).
\end{enumerate}
Only the third step involves a decision.
Appendix~\ref{sec:worked} walks through the steps on a small example.

\subsection{The formulation}\label{sec:details}

Write $F_1, \dots, F_T$ for one company's filings of one statement
(\cref{sec:have}), in the order published, with boundary $t$ between
$F_{t-1}$ and $F_t$. Within one filing, $v(k)$ is a tag's number for
period $k$ and $d(k)$ its \texttt{decimals}, with a subscript to say
which tag when it matters; each tag also has a label, the total it is
declared to add into, and a weight of $+1$ or $-1$.

\subsubsection{Strands}\label{sec:strands}

A strand is a run of one tag across consecutive filings joined on the
strongest evidence the filings offer, the same numbers on the shared
years, before any edit is weighed. It is an assumption: a
line's definition can change without changing the numbers both filings
show, and such a change is invisible to the method. Two numbers
\emph{agree} when they
differ by no more than half the coarser of their two rounding units, or
by no more than $\varepsilon$ if that is larger:
\[
  |v - v'| \le \max\bigl(\tfrac12\, 10^{-\min(d, d')},\ \varepsilon\bigr),
  \qquad \varepsilon = 0.5 \text{ (0.005 per share)}.
\]
Here $d$ is the fact's \texttt{decimals}, the number of significant
digits after the decimal point; Apple's 880 million in \cref{fig:xbrl}
carries $d = -6$, rounded to the nearest million, so any number within
half a million of it agrees. A number is
\emph{material} when $|v| > \varepsilon$.

A tag that appears in $F_i$ and $F_{i+1}$ joins across the boundary when
all three hold:
\begin{enumerate}
  \item the two unit classes do not conflict;
  \item no period both filings report disagrees;
  \item some material period agrees: a number carried over.
\end{enumerate}
Anything else breaks the run: a restated number, a missing filing, a
changed unit, or nothing to compare, when the two filings share no
period with a material number: a line that is always zero, a line
stated for one date only, a change of fiscal year, or a filing whose
comparatives our copy lacks. A break does not mean the two are not the
same line; it leaves a pair of endpoints to explain, and a cheap edit
may explain them. A heading carries no number and continues with its
tag: it is layout, not a line. The rule is conservative because a join
is irrevocable, while a missed join is made again whenever the
corresponding edit is in the solution.

Each strand $q$ runs from filing $\mathit{first}(q)$ to
$\mathit{last}(q)$. The endpoints to explain are
\[
  N = E \cup S, \qquad
  E = \{\, \mathit{end}(q) : \mathit{last}(q) < T \,\}, \qquad
  S = \{\, \mathit{start}(q) : \mathit{first}(q) > 1 \,\}.
\]

\subsubsection{Candidate edits}\label{sec:candidates}

A candidate is one possible edit, written as the set of endpoints it
would explain, with a price $w_c$ and a check it must pass to exist at
all (\cref{tab:edits}). A candidate is all or nothing: it covers every
endpoint in its set or none. It is a set rather than a pair because a
merge, a split, a carve-out or an absorb pairs one line with several;
matched one to one, the extra lines come out as ended or started and
the totals they take part in go unexplained. Each kind has
its own check, rather than one shared measure of what fails to add up,
because an explanation that simply ends lines is never tested against
any arithmetic and under a shared measure would win by default. Every
candidate but a rename sits at one boundary $t$: its ends are strands whose last filing is
$F_{t-1}$ and its starts are strands whose first filing is $F_t$. Each
kind is drawn beside its description: a filled dot is where a line's run
ends in filing $t$, an open dot where a run starts in filing $t{+}1$, and
a diamond is one edit, covering every endpoint it touches
(\cref{fig:edit-shapes} in Appendix~\ref{app:vocabulary} shows them all together).
Every check below is a sum or an equality on periods both reports
show, and a sum is evidence only through the terms that could have made
it fail: a term inside the sum's tolerance on the checked period adds
nothing, so a part counts only where it carries a material number
there. Parts $u_1, \dots, u_m$ \emph{add up} to $v$ when
\[
  \Bigl|\sum_j u_j - v\Bigr| \le \max\bigl(\tfrac{m}{2}\, 10^{-d},\ \varepsilon\bigr),
\]
with $d$ the coarsest precision among the parts and the total, so that
each rounded part may be off by half a unit; every check below that
asks for a sum uses this. Each kind is described here by its shape and
its check; Appendix~\ref{app:candidates} states every rule and bound
in full.

\begin{kind}{shape-rename}
\textbf{Rename.} One end and one start, $\{e, s\}$. Any end may pair
with any later start, at the next boundary or across a gap, and a tag
re-tied to itself after a restatement is a rename too. A pair is
refused when the units, the period types (a balance never continues as
a flow) or the per-share status conflict, when two different tags carry
opposite words (gross and net, current and noncurrent), or, across a
gap, when either tag or label appears in a filing the rename skips.
\end{kind}

\begin{kind}{shape-merge-split}
\textbf{Merge and split.} Several ends and one start, or one end and
several starts, at one boundary: two or three material parts that add
up to a material whole on every period all of them report, at least
one of those periods material. No part may be left out with the rest
still adding up, and parts the filer itself declares under the whole
are a break out or a fold, not a merge or split. \emph{Equal copies}
have the same shape, with every part equal to the whole. The search is
bounded (Appendix~\ref{app:candidates}).
\end{kind}

\begin{kind}{shape-carve-absorb}
\textbf{Carve-out and absorb.} A rename $e \to s$ whose numbers
disagree on a shared period, plus one to three lines that start
(carve-out) or end (absorb) whose weighted numbers explain the
restatement exactly on every period each child reports:
\[
  \begin{aligned}
    \text{carve-out: } v_{\mathrm{old}}(k) &= v_{\mathrm{new}}(k) + \textstyle\sum_j r_j\, u_j(k),\\
    \text{absorb: } v_{\mathrm{new}}(k) &= v_{\mathrm{old}}(k) + \textstyle\sum_j r_j\, u_j(k),
  \end{aligned}
\]
where $r_j$ is the ratio of the child's declared weight to the
parent's when both add into one total, and otherwise $-1$ when their
balance attributes (debit or credit) are opposite and $1$ when not. A
declared total above the parent or a child whose restatement equals the
movement the edit implies into it is explained by the edit and joins
the candidate; a declared total is never itself the parent. One line
can also be absorbed into two.
\end{kind}

\begin{kind}{shape-split-transfer}
\textbf{Split with transfer.} A split whose parts do not add up to the old
line, where the difference, period by period, is exactly the restatement of
the section total they add into. The candidate adds that total's rename,
$\{e_P\} \cup \{s_j\} \cup \{e_T, s_T\}$, with the total's restatement
credited: the same rule a carve-out or absorb applies to every total it
moves lines into or out of.
\end{kind}

\begin{kind}{shape-swap}
\textbf{Swap.} Two or more lines whose numbers moved onto each other's
tags, as a cycle, $\{e_1, s_1, \dots, e_m, s_m\}$: every leg's numbers
continue verbatim on the other tag, and every leg's old place is
accounted for.
\end{kind}

\begin{kind}{shape-fold-breakout}
\textbf{Fold and break out.} The parts' ends (a fold) or starts (a break
out) under a total that runs on across the boundary: at least two
parts, each material on a shared period, all declared under the total
and adding up to it, with the total's numbers unchanged. If the total's
own tag changes at the same boundary the candidate holds that rename
too, as explanation C of Tesla does.
\end{kind}

\begin{kind}{shape-ended-started}
\textbf{Ended and started} are not candidates. They are what an endpoint
costs when nothing covers it.
\end{kind}

\subsubsection{The optimization}\label{sec:choice}

Given the endpoints $N$ and the candidates $C$, each $c \in C$ a set of
endpoints with a price $w_c$, choose candidates that use no endpoint
twice, at the lowest total price, where an endpoint nothing covers pays
its own price $p_n$ for having ended or started:
\[
  \min_{x \in \{0,1\}^{C}} \;
  \underbrace{\sum_{c \in C} w_c\, x_c}_{\text{price of the chosen edits}}
  \;+\;
  \underbrace{\sum_{n \in N} p_n \Bigl(1 - \sum_{c \ni n} x_c\Bigr)}_{\text{price of the strand ends left uncovered}}
\]
subject to
\[
  \underbrace{\sum_{c \ni n} x_c}_{\substack{0:\ n\ \text{uncovered}\\ 1:\ n\ \text{covered}}} \le 1
  \quad \text{for every } n \in N .
\]
$x_c = 1$ chooses candidate $c$. The constraint says no endpoint is
covered twice, so the inner sum is $1$ if $n$ is covered and $0$ if
not.

Regrouping the objective by candidate instead of by endpoint
(Appendix~\ref{app:savings}) shows the same problem as weighted set
packing: each candidate saves $s_c = \sum_{n \in c} p_n - w_c$, the
price of leaving its endpoints unexplained minus its own price, and
the choice is the disjoint set of candidates with the largest total
saving. Candidates that save nothing are dropped
before solving, since removing one from any answer never makes it worse.
The guarantee is narrow: the chosen set is the cheapest among the
candidates \cref{sec:candidates} generates, whose listing is bounded
and pruned, under the prices of \cref{sec:cost}. An explanation the
listing did not produce cannot be chosen, and a different set of
prices can choose differently.

The problem splits into independent groups. Two candidates interact
only through a shared endpoint, so the connected components of the
candidates, two adjacent when they share an endpoint, are independent
problems, and each is solved on its own, exactly, as a 0--1 integer
program with one constraint per shared endpoint; a group with one
candidate takes it without solving anything, and choosing nothing is
always allowed, so every group has an answer. The groups are small. Over the 200 companies
of \cref{sec:results}, the 596 statement histories with a boundary, the
choice splits into
27,969 groups; 86\% hold one candidate and 95\% five or fewer, the
largest holds 3,351 candidates over 158 endpoints and solves in 41~ms,
and all the solves together take a quarter of a second. The answer need not be unique: two sets of candidates can save exactly
the same, and the implementation then picks one by a fixed order, which
only makes the choice repeatable.

Appendix~\ref{app:tesla-choice} lays out Tesla's 2015 to 2016 boundary
as this problem, with its candidates as sets and what each explanation
leaves uncovered.

\subsubsection{The prices}\label{sec:cost}

A price is roughly how much one would have to write down to describe
the edit beyond what the filings already say. There are many ways to
write such a price, and choosing one is a matter of judgment. The
prices decide the answer, since candidates compete for the same
endpoints: raising one kind's price makes it lose more often, lowering
it makes it win. Both the form of the prices and their constants took
trial and error against cases we knew the answer to, and many others
would work. What follows is the structure we settled on.

Prices are counted in digits. $\sigma(v)$ is the number of significant
figures of $v$ at its stated precision, trailing zeros removed (1,105
million has four; 880 million has two). Writing a number down costs
$c(v) = \alpha\,\sigma(v) + \beta$, and every price below is counted in
that unit, so different kinds of edit compete on equal terms. A number
that carries over is free and a changed or new number costs its
digits, so an edit cannot win by hiding numbers; evidence that singles
out one pairing earns a credit and evidence the filer guaranteed earns
none; one event is charged once, however many numbers it touched; and
no kind is rewarded for size (the six properties are listed in
Appendix~\ref{app:prices}).

Every price has the same shape:
\[
  w_c = \underbrace{A_c}_{\text{the act: how rare this kind should be}}
    + \underbrace{D_c}_{\text{what differs: names and numbers}}
    - \underbrace{K_c}_{\text{what matches: the evidence}} .
\]
The act is a flat price per kind of edit: a rename has none, a merge,
a carve-out or a swap has one, and raising it means that kind needs
more evidence before it wins. What differs is charged at its digits,
every word of a label and part of a tag that changed and every number
that is new or changed. What matches is credited, the same tag or
label, a number that agrees, a sum that checks, and a match many lines
share is worth little, one that singles out this pairing a lot, and
one the filer guaranteed nothing. Charging only for what
differs, with no credit, would leave the price saying nothing about
evidence: a tag whose numbers moved verbatim to another tag (Palantir,
Appendix~\ref{sec:palantir}) would price the same as that tag continuing with a
wholesale restatement, and the tie would go to keeping the tag. The
credit for numbers that agree is what lets the numbers override the
tag. \Cref{tab:prices} says what
the three terms are for each kind of edit; Appendix~\ref{app:prices} writes
them out in full, and Appendix~\ref{app:constants} gives the constants, which
set how much each term weighs. An uncovered endpoint is priced the same
way with no act to buy: a start pays for every number it shows, an end
only for the numbers that turn up on some other line in the next
filing.

\begin{table}[htbp]
  \centering
  \includegraphics[width=\linewidth]{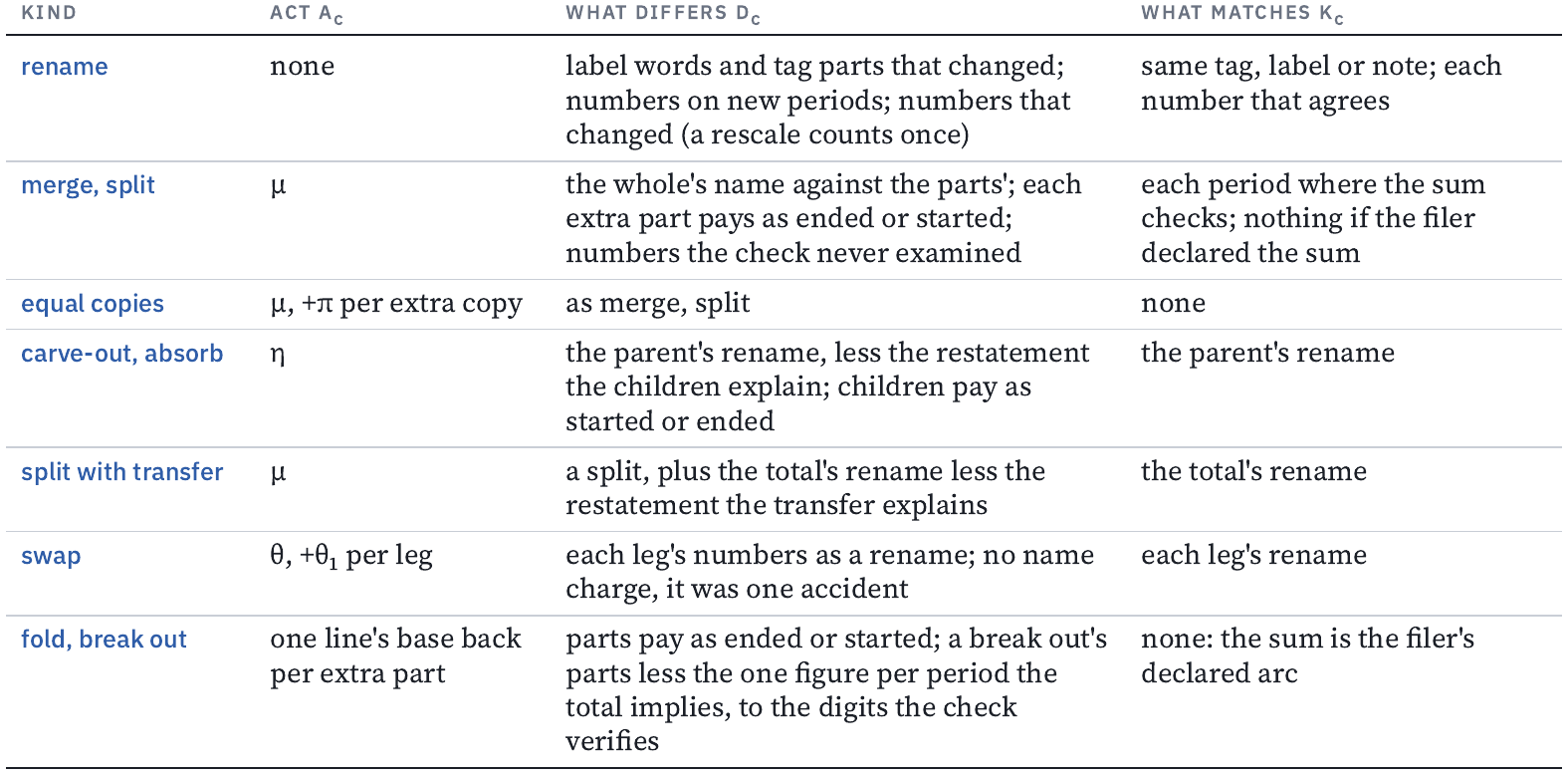}
  \caption{The three terms for each kind of edit. ``The parent's rename''
    means that line's rename price, computed the same way.}
  \label{tab:prices}
\end{table}

\subsection{The rows}\label{sec:output}

The solution is the cheapest set of edits between each pair of
successive filings. The renames among them, including those inside
carve-outs and absorbs, join strands end to start, and the joined
strands are the rows of the statement. A merge or a split supplies a
row with a derived number where no filing states one, the sum of the
parts, marked as derived and cited to the lines it adds up
(Appendix~\ref{sec:worked} and \ref{sec:jpm}). That number rests on
the merge's check having held on the periods both reports show, two
years for an income statement or cash flow statement and one instant
for a balance sheet; beyond them it is the check extended, which is
why it is marked.

\begin{figure}[htbp]
  \centering
  \includegraphics[width=\linewidth]{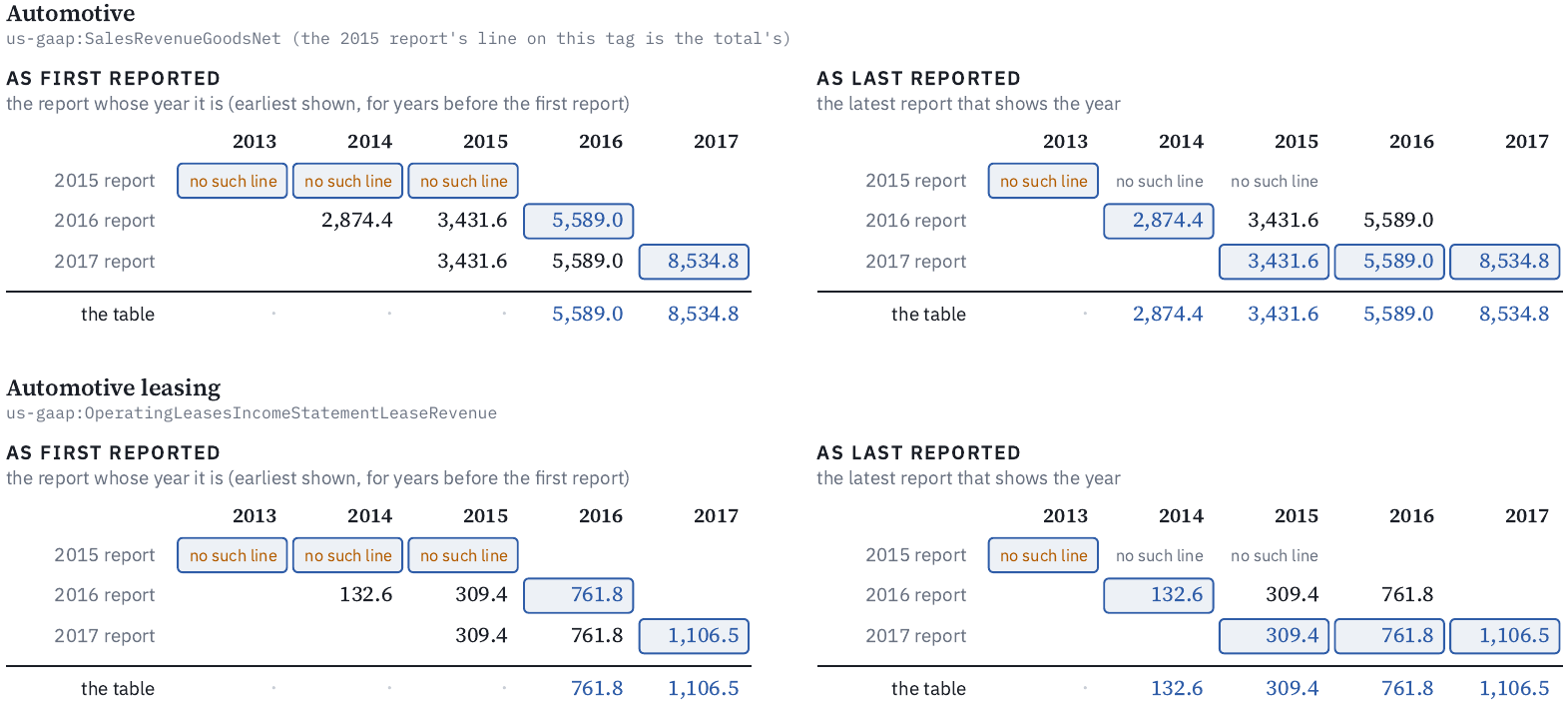}
  \caption{Which filing rules each year's cell, for two of Tesla's rows.
    Each row of a stack is one annual filing and each column a year; a
    filing shows its own year and the two before. The boxed cell is the
    ruling filing's, and the table row beneath is what it gives. As first
    reported takes the filing whose year it is (for 2013 and 2014, before
    the first filing shown, the earliest that shows them); as last
    reported takes the latest filing that shows the year. Where the ruling
    filing has no such line, the cell is blank rather than filled from a
    later filing.}
  \label{fig:ruling-filing}
\end{figure}

A row spans several filings, which may give different numbers for the
same year. So each year of the table is assigned a \emph{ruling}
filing, and its column contains only that filing's numbers. Each row's
number is looked up by the tag the row had in that filing. Every
number in a column comes from one filing, so the column should add up
as that filing did; \cref{sec:results-measured} checks that it does.
The table as first reported and the table as last reported differ only
in which filing rules each year: the
filing for that year, or the latest filing that still shows it
(\cref{fig:ruling-filing}). If the ruling filing has no such line, the
cell is blank. A number is never taken from a later filing.

A row is a lineage, and its numbers need not be comparable
throughout. The output tells three things apart. The
\emph{lineage} is the strands the chosen edits joined: one line
followed through the filings, whatever it was called. \emph{Comparable}
is narrower. A carve-out or an absorb changes what a continuing line
contains, and the filing that made the change restates only the years
it shows again, so a row is like for like within each ruling filing's
window and may change definition at a boundary where such an edit
sits. Apple's ``Other'' is one lineage with two definitions: 2019,
ruled by the fiscal 2021 report, excludes the proceeds from common
stock, and 2020 onwards includes them. The table marks every cell a
restatement changed, and the edits chosen at each boundary
(\cref{sec:results-measured}) say which kind of edit sits at each join
of a row, so the boundary where a definition changed can be read off;
the table does not restate the earlier years itself, since no filing
states those numbers. \emph{Reconstructed} is a number made from
others: the derived cell of a merge or split, the sum of its parts,
marked and cited to the parts.

Appendix~\ref{sec:solutions} shows all of this on real filings: the
candidates, their prices, the choice and the rows for each. Tesla's
rows, as first and as last reported, are in Appendix~\ref{sec:tesla}.

%% file: sections/05-results.tex
\section{Results and applications}\label{sec:results}

\subsection{Results}\label{sec:results-measured}

\paragraph{The two sets.} We use two sets of companies. The
\emph{specimen set}: thirty-two companies, three statements each (income
statement, balance sheet and cash flow statement), up to fifteen annual
reports apiece, 96 statement histories and 800 boundaries between one
annual report and the next. The method was developed against
them; the specimens of Appendix~\ref{sec:solutions} are among them, and
the constants of \cref{sec:cost} were settled by looking at how they came
out. The \emph{random sample}: two hundred companies drawn at random from
every company our data pipeline has processed with at least eight annual
reports: 597 statement histories, one of them a single report with no
boundary, and 6,103 boundaries. The sample was held out until the method
had been built and first evaluated; after that, its cases informed
three revisions to the method, listed in Appendix~\ref{app:revisions}.
The results below are the final method on both sets. The draw includes companies
that have since been acquired, gone private or deregistered; a sample of
survivors would have been easier to compare with a vendor's series, but
it would not represent the filings. Every number below comes from one
run of the implementation, with nothing set by hand for any company.

\paragraph{Changes of shape.} \Cref{tab:changes} counts, by statement,
the boundaries at which the optimization chose at least one change of
shape: a merge, split, carve-out, absorb, fold, break-out, swap or
transfer. These are chosen edits, not candidates. A change of shape is
one the tag cannot follow, and 30 of the 32 specimen companies and 194
of the 200 in the random sample have one somewhere in their history;
the cash flow statement carries them most.

\begin{table}[htbp]
  \centering
  \includegraphics[width=\linewidth]{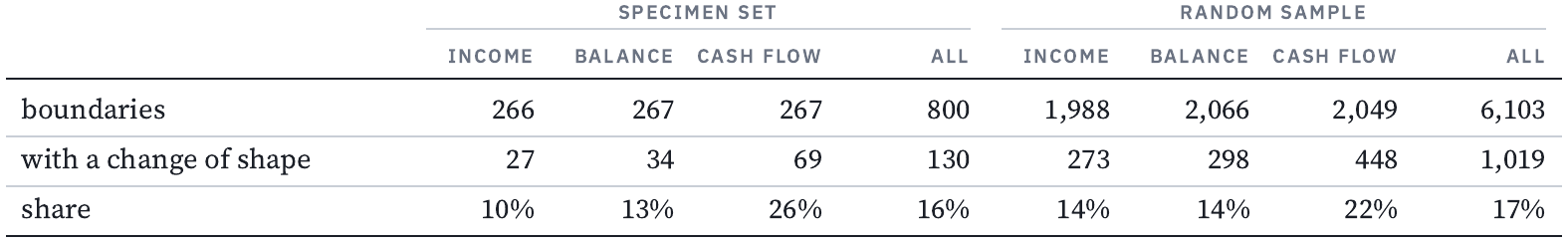}
  \caption{Boundaries with a change of shape, by statement. A boundary is
    one pair of consecutive annual reports of one statement; it counts
    when the optimization chose at least one edit other than a
    continuation there.}
  \label{tab:changes}
\end{table}

\paragraph{The edits chosen.} \Cref{tab:counts} counts every edit the
method chose, from the same runs. Across a break the
optimization makes 3,431 continuations on the specimen set and 39,882 on
the random sample, a third and a quarter of them changing the tag and
the rest keeping it across a restatement or a gap in the numbers; 3,048
and 16,638 numbers are restated. Continuations and restatements
dominate; every other kind is rare, and each occurs on both sets.
Written out for one company, Apple's cash flow statement, the boundaries
read: \emph{2011 annual report}: renamed the tax paid on net share
settlement to a company tag; absorbed \emph{gain or loss on property,
plant and equipment} and \emph{other current assets} into \emph{other
operating assets}; restated 2 numbers. \emph{2022 annual report}:
absorbed \emph{proceeds from issuance of common stock} (880 for 2020)
into \emph{other financing} (754); absorbed \emph{purchases} and
\emph{proceeds from sale and maturity of other investments} into
\emph{other investing}. The 2022 boundary answers the
question \cref{sec:problem} opened with. Appendix~\ref{app:hard} gives
one case per kind, chosen where the names give no help.

\begin{table}[htbp]
  \centering
  \includegraphics[width=\linewidth]{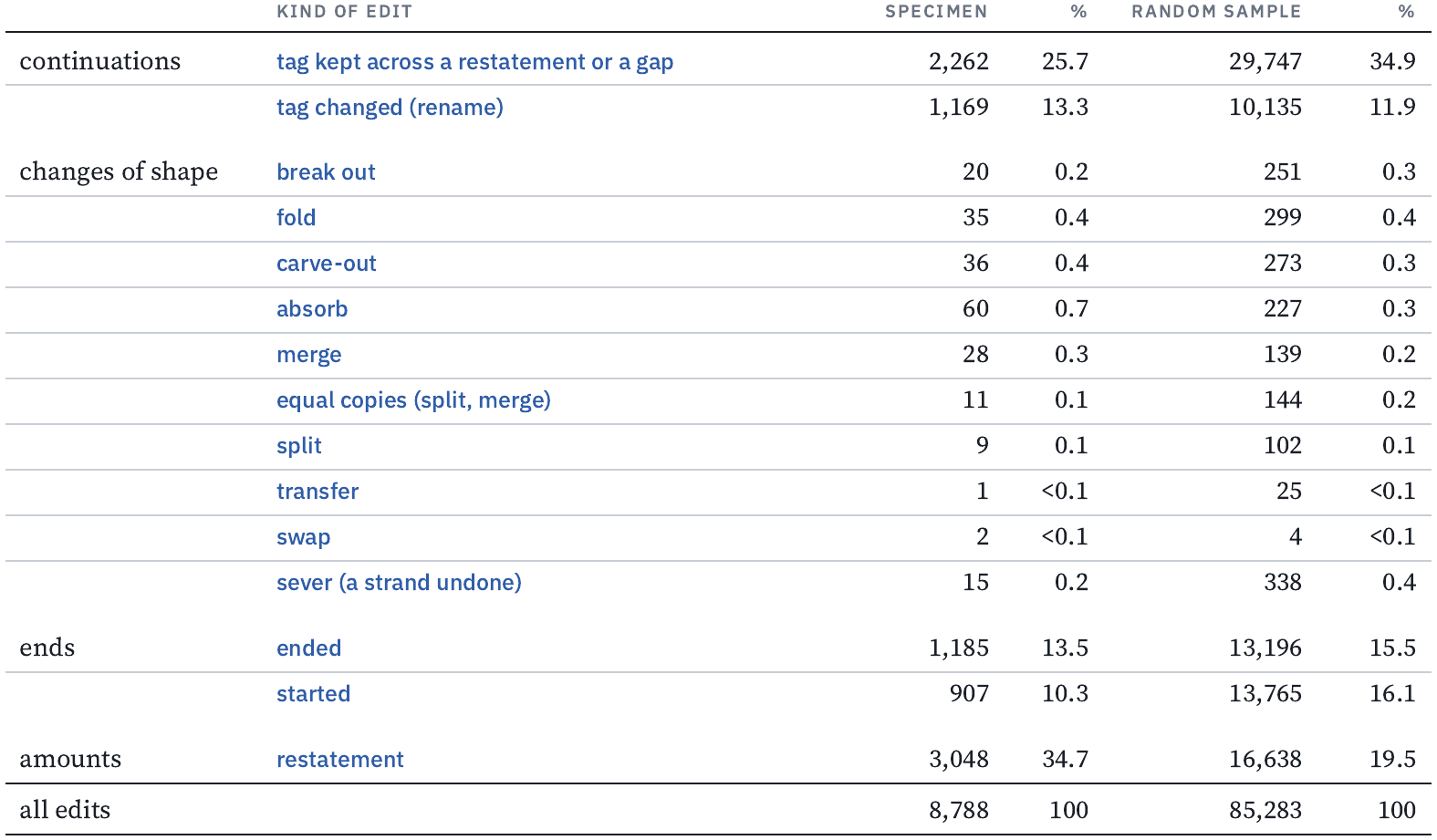}
  \caption{Every edit the method chose, all three statements: the
    specimen set (32 companies, 800 boundaries) and the random sample (200
    companies, 6,103 boundaries), each kind with its share of all the
    edits chosen on that set. A continuation is one line carried across a
    boundary where the strands did not already carry it.}
  \label{tab:counts}
\end{table}

\paragraph{Two more cases.} The changes of shape are few, so the best
way to see what the optimization finds is to follow some of them. Apple's
absorbs (\cref{sec:problem}) and Tesla's boundary (\cref{sec:example})
have been followed already; Appendix~\ref{sec:solutions} gives all
four in full, as facts, strands, candidates, choice and rows. Palantir
filed three years of cash flow numbers on the wrong tags in its fiscal
2021 report and moved them back in 2023: a swap of three legs beats the
four same-tag continuations, and the rows are the true series, each
named by the tag it ends on (Appendix~\ref{sec:palantir}). JPMorgan's
2023 report replaces one line of non-interest revenue with two and
restates both comparative years into the parts: a split beats a rename
with a carve-out, and the old line continues as the sum of its parts,
a derived row cited to both (Appendix~\ref{sec:jpm}).

\paragraph{Arithmetic checks.} Two checks, both against what the filings
declare. First, the rule of \cref{sec:have} that every line is
accounted for exactly once: every line of every filing ends up in
exactly one row, note or flag, and nothing is used twice; 0 violations
in 96 statements and in the 596 with a boundary. Second, the declared totals:
the finished table must still add up, year by year, once the lines are
joined across filings. A check is one declared total in one year; it
passes, passes within the filer's rounding, fails, or cannot be run
because a line it needs has no number in that column. \Cref{tab:arms}
counts them.

\begin{table}[htbp]
  \centering
  \includegraphics[width=\linewidth]{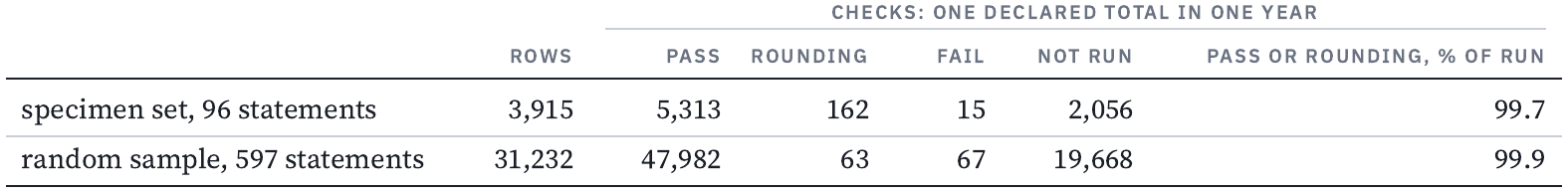}
  \caption{Declared totals re-checked in the finished tables. The method
    also marks 18 and 115 checks not comparable, where the total and a
    part are stated by different filings.}
  \label{tab:arms}
\end{table}

Of the checks run, 99.7\% pass or pass within rounding on the specimen
set and 99.9\% on the random sample, and no statement type, taken over
all companies, is below 99.5\%; single company statements go lower, to
93.7\% for GE's balance sheet. The specimen set's 15 failures are all in GE and
Berkshire, whose statements carry more than one entity in their columns
(GE Industrial and GE Capital; Berkshire's segments): the model treats
a breakdown row as a line of its own, and a total whose parts sit in
another column does not add up on our page. The random sample's 67 are
in 33 companies, 33 of them within a tenth of a percent of the
total. The checks that cannot be run are mostly years for which no
filing declares arithmetic for the line, and otherwise years where a
part has no number in that column (Appendix~\ref{app:layout}); those
years are left unchecked rather than checked against another filing's
equations.

\paragraph{Layout invariance.} Each history can be shown in any filing's
layout. Rendered in the newest filing's layout and in the oldest
presenter's, all 96 specimen statements are identical in every cell and
every check over a five-year window, and 94 of 96 over fifteen; on the
random sample, 564 of 597 over fifteen. Of the 35 that differ, 21
differ only in which checks run; in the other 14, every differing row
was read, and the same numbers sit under a differently named row, or
a merged line that one layout derives as the sum of its parts appears
in the other as the parts alone (Appendix~\ref{app:layout}). No cell
holds two different numbers for the same line. An earlier run of this
check found three wrong renames of a balance to a flow, and the rule
of \cref{sec:candidates} that refuses the pairing came from them
(Appendix~\ref{app:revisions}). The rows come from the strands, so
the layout decides how the numbers are ordered, labelled and grouped,
and which derived rows exist, not what the numbers are.

\paragraph{Agreement with a vendor.} A commercial provider (FMP)
publishes standardized series keyed by stock ticker, for companies
listed today. Our sets span fifteen years, and 4 of the specimen
companies and 75 of the random sample have since deregistered or been
delisted, or file without a listing; one more in each set has no series
our anchoring could use. That leaves 27 and 124. We anchor
each of its fields to one of our lines in the newest year both cover,
walk our history backwards under the method and under a floor, the
tag alone, the anchored tag followed through every filing whatever
its numbers do, and ask, year by year, whether the line we followed
carries the number the vendor shows. It is a test of the lineage, not of the rendered
cell. Any filing's version of the number counts, since the vendor holds
as-filed numbers and we show the ruling filing's, and two numbers
agree within one part in ten thousand or half a million dollars,
whichever is larger (one cent for a per-share figure). A derived cell,
the sum a merge or split supplies, is in no filing and is not
confirmed by this test. A vendor value of zero is skipped, and a vendor
number no filing states is the vendor's own arithmetic and is left out
(1,999 and 7,947 of them), a test made against the vendor's
as-reported feed~\cite{fmp-as-reported}. \Cref{tab:fmp} is the
result, with the tag alone as the floor.

\begin{table}[htbp]
  \centering
  \includegraphics[width=\linewidth]{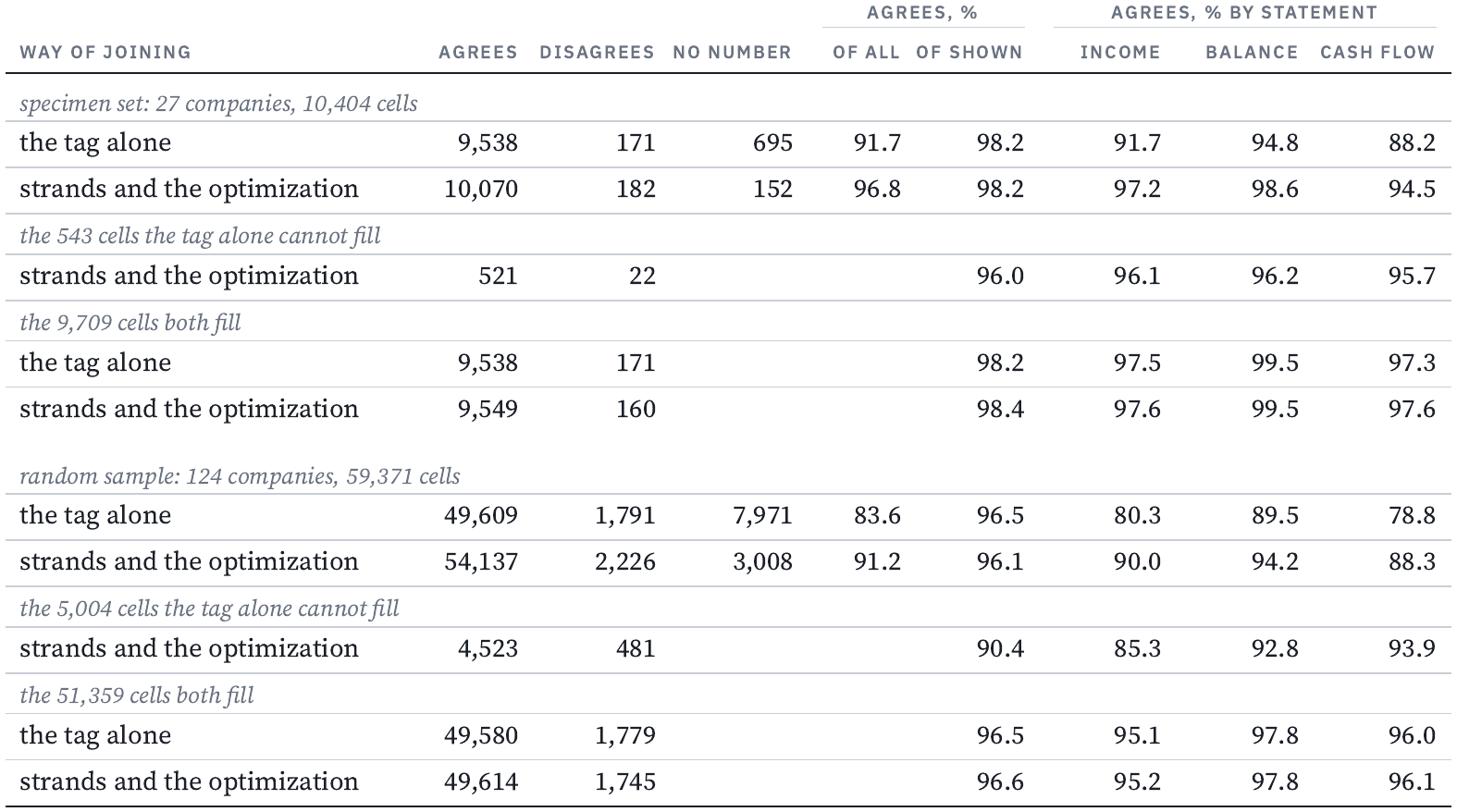}
  \caption{Agreement with a data vendor's series along the
    reconstructed histories. A cell agrees if the line we followed
    carries, in some filing, the number the vendor shows for that year;
    ``no number'' means our history has no value there. ``Of shown''
    counts only cells where we show a number. The rows for the cells
    the tag alone cannot fill score the optimization where it adds
    coverage; the rows for the cells both fill score the two ways of
    joining on the same cells.}
  \label{tab:fmp}
\end{table}

The tag alone is a strong floor. On the cells it reaches it agrees
with the vendor 98.2\% and 96.5\% of the time, and the method keeps
that (98.4 and 96.6\% on the same cells): where a tag runs on with
its numbers unchanged the strand rule joins it, so there the two are
the same chain. What the optimization adds is the cells the tag
cannot reach: 543 on the specimen set, 5.6\% more than the tag alone
fills, and 5,004 on the random sample, 9.7\% more (41 go the other
way, so 4,963 net), which takes agreement over all cells from 91.7 to
96.8\% and from 83.6 to 91.2\%. Every one of those cells lies past a
boundary where the tag changed or the line was reshaped, and the edits
of \cref{tab:counts} are what reach them: 3,431 continuations and 202
changes of shape on the specimen set, 39,882 and 1,464 on the random
sample. On the added cells the method agrees with the vendor 96.0\%
and 90.4\% of the time, where the tag alone shows nothing; these are
the hard cells, and the rate over each method's own cells (98.2\% for
both, 96.1 against 96.5\%) is lower for the method only because it
includes them. The cash flow statement is the hard one for both, on
both sets; it is where filers reshape most: 26\% of its boundaries on
the specimen set and 22\% on the random sample carry a change of shape,
against 10 to 14\% for the other two statements
(\cref{tab:changes}).

The vendor's series are not ground truth, and the vendor has a
different goal. It standardizes so that companies can be compared with
one another~\cite{fmp-income-statement}, the purpose the SEC gave for
standard tags~\cite{sec-33-9002}; we standardize one company across
time and keep its own lines. Most of the misses are the vendor's own
arithmetic: it nets some lines, gathers others into residual buckets,
adjusts per-share figures for later splits, and applies its own
definitions of net income and of payables, each of which puts a number
in its series that no filing states for the line we anchored.
Appendix~\ref{app:misses} sorts every miss; the identity errors the
comparison finds or cannot rule out come to 61 cells on the specimen
set, 0.6\%, and 1,844 on the random sample, 3.1\%, most of them one
line that the filer kept under two tags for a while and the
optimization left in two rows. These are agreement rates, not
accuracy: two series can agree on a wrong number, and a cell that
agrees by coincidence is not found by examining disagreements. What
the comparison does establish is the coverage gain, and that every
cell where the vendor's standardization departs from the filings
counts against us. The method and the tag alone are scored against the
same series with the same anchors, so the rows of \cref{tab:fmp} can
be compared; the absolute level is set by the vendor's
standardization, not by us.

\subsection{Applications}\label{sec:uses}

The output is the strands with the edits between them, and several
products follow. Each needs both the numbers and the links, with the
filing behind each: the filings have every number but no link between
one year's filing and the next (\cref{sec:parts}), and a vendor has
links, made its own way, but keeps one version of each number,
standardized onto its own fields, without citation. What the output is
not: it is not standardized across companies, since a row is one
company's line, not ``revenue'' in general, and putting two companies
side by side is the vendor's job; it does not say whether a change
mattered to the business; and its checks reach only as far as the
arithmetic the filers declare.

\paragraph{Cited tables.} A company's statements over fifteen years in
any filing's layout, where every cell names the filing and the fact it
came from, every break in a line is marked, and a rebuilt total is a
checked formula over the lines it adds up (\cref{sec:output}). What
this adds over a vendor's series is the citation on each number and a
mark where a line ends.

\paragraph{Two axes of time.} Every number a company reports is
reported again, for two or three years, by later filings, and
sometimes changed. The strands keep every version, each dated by its
filing: a lattice with two axes of time, the period a number is for
and the date it was said, of which any table is a slice: the table as
it stood on any past date, the history as first reported against as
last reported (\cref{sec:output}), or the restatement record of one
line. Economists keep such vintages for macroeconomic series
(\cref{sec:related}); nothing published keeps them for the lines of a
company's statements.

\paragraph{The record of changes.} The method's by-product is a record
of what each company changed in how it reports, boundary by boundary
and by kind: on the specimen set, 130 of 800 boundaries carry a change
of shape, and \cref{tab:counts} says which. Companies work this mapping
out every time they restate, and the filings carry no consistent,
machine-readable form of it (\cref{sec:related}); this is the nearest
thing, built from the outside. It says which lines are stable, when a
company reorganized its reporting, and how one company's habits
compare with another's. It does not say why; the reason is in the
filing's own words, and the record says which boundary and which
lines to read it for.

\paragraph{An interface for programs.} A program asking about a
company's history needs a stable identity for a line across years, a
source for every number, and a way to check a table it did not build.
The strands are the identity, the cells carry the source, and the
declared totals are the check, so the natural interface is a query
over strands, cells and edits, with the tables one view of it:
\emph{revenue for the last ten years as first reported, with the
filing for each}; \emph{does this table an assistant produced agree
with the filings, cell by cell}.

%% file: sections/06-next.tex
\section{Extensions}\label{sec:extensions}

The method needs no mapping written for a company and no dictionary of
related tags. It does not know that \texttt{OtherOperatingActivities}
is an operating item or that commissions are a kind of non-interest
revenue. The evidence it works from is in the filings: the numbers on
each line, the prior years each filing shows again, the units, whether
a number is stated at a date or over a period, the arithmetic each
filing declares, and the debit or credit mark on each line. Names are compared as strings, by the words two labels share and
the characters two tags share, with one exception: a short list of word
pairs (gross and net, current and noncurrent) that a rename may not
cross. A tag the company invented is handled the same way as a standard
one. The method should carry over to any series of reports in which
each report overlaps the one before.

The method's input is XBRL facts: a tag, a period, a number and the
arithmetic the filing declares. Filings that carry none, annual reports
before 2009, most forms other than the 10-K and 10-Q, and reports filed
outside the SEC, could be brought in by extracting their numbers and
assigning each a tag first. That is the task recent work has set for
language models: FinTagging benchmarks the extraction of numbers from
the text and tables of a report and their linking to the US GAAP
taxonomy~\cite{wang2025fintagging}, and XBRLTagRec assigns tags to
numbers by retrieval and re-ranking~\cite{hu2026xbrltagrec}. Facts
produced that way would feed the method unchanged, with two losses: a
tag assigned by a model is weaker evidence of identity than one the
filer chose, so the prices of \cref{sec:cost} would need to weigh it
less, and a report with no declared arithmetic gives the checks
nothing to foot against, so folds, break-outs and the re-checked
totals would need the arithmetic inferred as well.

The record of changes (\cref{sec:uses}) also supports questions this
paper does not ask. Do reorganizations cluster around new versions of
the tag dictionary? Are companies reporting less detail over time?
Which lines are stable across a whole industry, and which never are?

%% file: sections/07-coda.tex
\section{Conclusion}\label{sec:coda}

A financial statement that covers many years has to be built from many
filings, and the filings report the same measures differently as the
company changes what it shows and how it labels it. Most lines keep
their label and tag from one filing to the next and can be followed by
name. But neither is fixed: both change while the line continues, and
both are reused for lines that are new. A name is only evidence.
The filings carry no consistent, machine-readable record of the
changes, at most a sentence in a note, so the changes have to be
inferred from the filings, and the only certain evidence is the
arithmetic each filing declares and the prior years each shows again
under its new lines.

We set that inference up as an optimization problem. At each boundary
between consecutive filings, the candidate edits of each known kind
that pass its check are listed, each with the lines it would explain
and a price for how much it asks the reader to believe beyond what the
filings show. An integer program
picks the cheapest set of the listed edits that accounts for every line
exactly once, and the chosen edits trace each line through the filings.
The
result is a statement whose rows run across every year the filings
cover, in any filing's layout, with each number cited to the filing
that stated it and each declared total re-checked.

We ran it on 232 companies. Where the line the method follows and a
commercial vendor's series both give a number, they agree 96 to 98\%
of the time. The table reaches 6 to 10\% more cells than the tag
alone, every one past a change of tag or of shape, and agrees with the
vendor on 96\% and 90\% of those; 99.9\% of the declared totals that
can be checked add up. A company works out
this mapping every time it restates, and its filings carry no
machine-readable form of it. The method recovers it from what the
filings do state: their arithmetic and the years they show again.

%% file: sections/a-model.tex
\section{The method in full}\label{app:model}

\subsection{The procedure as pseudocode}\label{app:algorithm}

\Cref{alg:method} states the four steps of \cref{sec:how} as one
procedure, from the filings to the edits and the rows.
\begin{algorithm}[H]
  \caption{Rebuilding a statement's history from its filings.}
  \label{alg:method}
  \begin{algorithmic}[1]
    \Require filings $F_1, \dots, F_T$ of one company and statement, in the order published
    \Ensure the edits between each pair of filings, and the rows that follow from them
    \Statex
    \State $S \gets$ the strands of $F_1, \dots, F_T$: runs of one tag joined on the same numbers in consecutive filings
    \State $N \gets$ every point where a strand ends before $F_T$ or starts after $F_1$
    \State $C \gets \varnothing$
    \For{each filing boundary $t$}
      \For{each kind of edit $k$: rename, merge, split, carve-out, absorb, fold, \dots}
        \State $C \gets C \cup \Call{Candidates}{k, S, t}$, each with its endpoints, its check and its price
      \EndFor
    \EndFor
    \State $C \gets \{\, c \in C : c \text{ is priced below leaving its endpoints unexplained} \,\}$
    \State $\mathcal{G} \gets$ the connected components of $C$, two candidates adjacent when they share an endpoint
    \State $X \gets \varnothing$
    \For{each component $G \in \mathcal{G}$}
      \State $X \gets X \cup \Call{Choose}{G}$, the cheapest set in $G$ using no endpoint twice, solved exactly
    \EndFor
    \State every endpoint not covered by $X$ has ended or started
    \State rows $\gets$ strands linked end to start by the renames in $X$
    \State \Return $X$, the rows, and every candidate that lost, with its price
  \end{algorithmic}
\end{algorithm}

\subsection{The savings form}\label{app:savings}

The optimization of \cref{sec:choice} minimizes the price of the chosen
edits plus the price of the endpoints left uncovered,
\[
  \sum_{c \in C} w_c\, x_c + \sum_{n \in N} p_n \Bigl(1 - \sum_{c \ni n} x_c\Bigr),
\]
subject to no endpoint being covered twice. Multiplying out the second
term gives $\sum_n p_n - \sum_n p_n \sum_{c \ni n} x_c$. The double sum
has one term $p_n x_c$ for every endpoint $n$ and every candidate $c$
that holds it, so it can be added up by candidate instead:
\[
  \sum_{n \in N} p_n \sum_{c \ni n} x_c
  \;=\; \sum_{\substack{(n,\,c)\\ n \in c}} p_n\, x_c
  \;=\; \sum_{c \in C} x_c \sum_{n \in c} p_n .
\]
The objective is then
\[
  \sum_{n \in N} p_n + \sum_{c \in C} x_c \Bigl(w_c - \sum_{n \in c} p_n\Bigr)
  = \sum_{n \in N} p_n - \sum_{c \in C} s_c\, x_c ,
  \qquad s_c = \sum_{n \in c} p_n - w_c ,
\]
where $s_c$ is what candidate $c$ saves. The first term does not depend
on $x$, so minimizing the price is maximizing the total saving of a
disjoint set of candidates: weighted set packing.

\subsection{The candidates in full}\label{app:candidates}

The rules of \cref{sec:candidates}, in full: for each kind of edit, the
shape of the candidate, the check it must pass to exist, and the bounds
on its search.
\paragraph{Rename.} One end and one start: $\{e, s\}$.
\begin{itemize}
  \item Any end may pair with any later start, at the next boundary or across
    a gap. A tag re-tied to itself, after a restatement or after a break
    with nothing to compare, is a rename too.
  \item Refused when the unit classes conflict, when one line is per share
    and the other is not, when one line is stated at instants and the
    other over periods (a balance never continues as a flow: the two
    share no period, so no number could ever test the pairing), or when
    two different tags carry opposite words (gross and net, current and
    noncurrent; a tag carrying both words of a pair takes no side).
  \item A rename across a gap is refused if, in any filing it skips, either
    tag or either label appears: the line was not missing there, it was
    somewhere else.
\end{itemize}

\paragraph{Merge and split.} Several ends and one start, or one end and
several starts, at one boundary.
\begin{itemize}
  \item Two or three parts, every one of them material, and the whole
    material.
  \item The parts add up to the whole on every period where the whole and
    every part report a number, and at least one of those periods is material
    (two, for three parts).
  \item No part can be left out with the rest still adding up, and the parts are not
    simply the lines the filer itself declares under the whole: that is a
    break out or a fold, not a merge or split.
  \item \emph{Equal copies} have the same shape, with every part equal to the
    whole instead of adding up to it.
  \item The search is bounded: pairs and triples only where a boundary has few
    enough ends (or starts), and only the cheapest few merges and splits per
    boundary are kept (Appendix~\ref{app:constants}).
\end{itemize}

\paragraph{Carve-out and absorb.} A rename $e \to s$ across boundary $t$,
plus new lines (carve-out) or ended lines (absorb): $\{e, s\} \cup$ the
children's starts or ends.
\begin{itemize}
  \item The rename's numbers must disagree on some shared period, and not all
    by one ratio (that is a rescale, priced in the rename).
  \item One, two or three children whose weighted numbers explain the
    restatement exactly. On every period each child reports,
    \[
      \begin{aligned}
        \text{carve-out: } v_{\mathrm{old}}(k) &= v_{\mathrm{new}}(k) + \textstyle\sum_j r_j\, u_j(k),\\
        \text{absorb: } v_{\mathrm{new}}(k) &= v_{\mathrm{old}}(k) + \textstyle\sum_j r_j\, u_j(k),
      \end{aligned}
    \]
    as sums that agree, where $r_j = \omega_j / \omega_{\mathrm{parent}}$ if the
    child and the parent add into the same total; otherwise $r_j = -1$
    when the child's balance attribute (debit or credit) is the opposite
    of the parent's, and $1$ when it is the same or undeclared. So
    interest expense ($-1$ into pre-tax income) absorbed by other income
    ($+1$) enters with $-1$, and a liability carved out of a receivable
    it was netted inside enters with $-1$.
  \item Children are tried one at a time first; a pair only from lines no
    single child explained, and a triple only from pairs that failed.
  \item A declared total above the parent or a child, at any height,
    whose own strand broke at the boundary and whose restatement equals
    the movement the edit implies into it along the declared arcs on the
    footing periods, is explained by the edit: its rename joins the
    candidate and that restatement is credited, as in a split with
    transfer. Every other continuing line of that total must be unmoved
    on those periods, and at least one must exist.
  \item A declared total is never the parent. Its restatement is
    explained by the movement of its parts, and a total carving a line
    out of itself would count that movement twice.
  \item \emph{Absorbed into two}: one ended line $c$ and two restated renames
    $e_1 \to s_1$, $e_2 \to s_2$ whose changes add up to $c$ on every period it
    reported: $\{e_c, e_1, s_1, e_2, s_2\}$.
\end{itemize}

\paragraph{Split with transfer.} A split whose parts do not add up to the old
line, where the difference, period by period, is exactly the restatement of
the section total they add into. The candidate adds that total's rename,
$\{e_P\} \cup \{s_j\} \cup \{e_T, s_T\}$, with the total's restatement
credited: the same rule a carve-out or absorb applies to every total it
moves lines into or out of.

\paragraph{Swap.} Two or more lines whose numbers moved onto each other's
tags, as a cycle. The candidate is every leg's end and start,
$\{e_1, s_1, \dots, e_m, s_m\}$.
\begin{itemize}
  \item Each leg's numbers continue verbatim on the other tag.
  \item Each leg's old place is accounted for: its tag is another leg's
    destination, or it is absent from the new filing, or its run there breaks.
\end{itemize}

\paragraph{Fold and break out.} The parts' ends (a fold) or starts (a break
out), under a total that runs on across the boundary.
\begin{itemize}
  \item At least two parts, each material on a period both reports show,
    all declared under the total, and the declared sum adds up. Two new
    lines that are zero in every shared year add up under an old total
    as well as a real break out would, and are lines that started.
  \item The total's numbers do not change; a total that restated is an absorb
    or a split with transfer.
  \item If the total runs on in one strand, the candidate is the parts alone.
    If the total's own tag changes at the same boundary, the candidate also
    holds that rename's end and start, as explanation C of Tesla does.
  \item A part whose tag has another life across the boundary is
    admitted only when the candidate's own rename carries that other
    life, as Tesla's reused automotive tag is. Otherwise a fold and a break out of the same tag
    would turn a restatement into two edits and cut the row.
\end{itemize}

\paragraph{Ended and started} are not candidates. They are what an endpoint
costs when nothing covers it.

\subsection{The prices in full}\label{app:prices}

Each price of \cref{sec:cost}, written out. Constants are in
Appendix~\ref{app:constants}. Every price respects six things:
\begin{enumerate}
  \item A number that carries over is free; a changed or new number
    costs its digits, wherever it appears. An edit cannot win by hiding
    numbers.
  \item Names: the closer two labels, and two tags, the cheaper the
    pairing.
  \item Evidence that singles out one pairing among many (the same tag,
    the same label, a number that agrees, a sum that checks) earns a
    credit; evidence shared by many lines earns little; evidence that
    was guaranteed (a sum the filer declared) earns nothing.
  \item One event is charged once: a swap, a carve-out or a change of
    scale is an act, not one charge per number it touched.
  \item No reward for size: every extra line an edit covers pays its
    own way.
  \item A line that ends explains itself, unless its numbers turn up
    elsewhere; a line that starts brings every number it shows.
\end{enumerate}

\paragraph{Ended and started.}
\[
  p_s = \kappa_S + \sum_k c\bigl(v(k)\bigr),
  \qquad
  p_e = \kappa_E + \sum_{k \in \mathit{echo}(e)} c\bigl(v(k)\bigr),
\]
where $\mathit{echo}(e)$ is the set of material numbers of the ended line
that appear on some other line in the next filing but on no other line
in its own.

\paragraph{Rename.} The price of what differs, $D = D_{\mathrm{name}} +
D_{\mathrm{num}}$, less a credit for what matches, $K$:
\[
  w_{e \to s} = D_{\mathrm{name}}(o, n) + D_{\mathrm{num}}(o, n) - K(o, n).
\]

$J$ is the share of label words two labels have in common; $\mathit{sim}$
the share of characters two names have in common.
\begin{align*}
  D_{\mathrm{name}}(o, n) &= \lambda\,\bigl(1 - J(\lambda_o, \lambda_n)\bigr)
    + \sum_{\mathclap{\text{tag parts that differ}}} g(\text{part}), \\
  g(\text{prefix only}) &= \gamma_0, \qquad
  g(\text{concept}) = \gamma_C + \delta_C\,(1 - \mathit{sim}), \\
  g(\text{axis or member}) &= \gamma_A + \delta_A\,(1 - \mathit{sim}), \qquad
  g(\text{breakdown added or removed}) = \gamma_D ; \\[6pt]
  D_{\mathrm{num}}(o, n) &= \sum_{k \notin o} c\bigl(v_n(k)\bigr) + B, \\
  B &= \sum_{\mathclap{k\ \text{disagrees}}} c\bigl(v_n(k) - v_o(k)\bigr),
    \ \text{capped at } \rho \text{ if one ratio moved every shared number;} \\[6pt]
  K(o, n) &= \frac{\kappa(M)}{m_\lambda m'_\lambda}\,[\text{same label}]
    + \frac{\kappa(M)}{m_\nu m'_\nu}\,[\text{same note}]
    + \kappa(M)\,[\text{same tag}] \\
  &\quad + \sum_{\mathclap{k\ \text{agrees}}} \frac{\min\bigl(c(v(k)),\ \kappa(M)\bigr)}{m_k m'_k},
\end{align*}
where $\kappa(M) = \alpha \log_2(1 + M)$ and $m, m'$ count the lines on the
old and new side that share the evidence (a tag is unique on each side).
A pair is refused outright, at any price, when the unit classes conflict,
when one line is per share and the other not, when one line is stated at
instants and the other over periods, or when the two tags carry opposite
words.

\paragraph{Merge and split.}
\[
  w = \max\Bigl(\max\bigl(\mu + X + \delta_N\,(1 - \mathit{sim}),\ w_{\min}\bigr) + U + R - \Phi,\ w_{\min}\Bigr).
\]
\begin{itemize}
  \item $\mu$ is the act; $w_{\min}$ keeps the price positive, before and
    after the credit.
  \item $X$: every part after the first pays its own ended or
    started price. Equal copies repeat the whole, so each extra copy
    pays a flat $\pi$ instead.
  \item $\delta_N (1 - \mathit{sim})$: the best name similarity
    between the whole and any part, or the parts' names run together in any
    order.
  \item $U$: the new side's numbers on periods the check did not
    examine, at $c(v)$ each. A footing split's later parts already paid
    in full through $X$, so only its first part is charged here; an
    equal split's parts all are.
  \item $\Phi$: for each period where the whole is material and
    the sum agrees within a relative tolerance $\varphi$, a credit of
    $\alpha \min\bigl(\bar\sigma, \sigma(v)\bigr)$, where
    $\bar\sigma = \log_{10}(1/\varphi)$ is the number of figures such a check
    can verify. Equal copies earn no such credit.
  \item $R$: when the parts all add into one declared total and the
    whole carries that total's figures (the total under another tag),
    the sum agrees by construction, and the footed numbers are charged
    in full at $c(v)$ as if they were new. Parts that add into the whole
    itself are refused: that is a break out, not a split.
\end{itemize}

\paragraph{Carve-out and absorb}.
\[
  w = w_{e \to s} - B^{*} + \eta + \sum_{j} p_j
      + \sum_{T} \bigl(w_{e_T \to s_T} - B^{*}_T\bigr) .
\]
\begin{itemize}
  \item $B^{*}$ is the part of the parent's restatement on the periods
    the children explain. It is replaced by one act $\eta$; the children
    pay their own ended or started prices.
  \item $T$ ranges over the declared totals above the parent and the
    children that the edit explains (\cref{sec:candidates}): each
    contributes its rename and is credited its restatement $B^{*}_T$ on
    the footing periods.
  \item The candidate exists only when $B^{*} + \sum_T B^{*}_T > \eta$:
    restatements cheaper than the act are left to the renames.
  \item Absorbed into two:
    $w = (w_1 - B^{*}_1) + (w_2 - B^{*}_2) + \eta + p_{e_c}$.
\end{itemize}

\paragraph{Split with transfer.} The split's act, extra parts and name
term, its unchecked numbers, plus the section total's rename less the
restatement the transfer explains:
\[
  w = \max\bigl(\mu + X + \delta_N\,(1 - \mathit{sim}),\ w_{\min}\bigr)
      + U + w_{e_T \to s_T} - B^{*}_T .
\]
$U$, the new side's numbers on periods the check did not examine, is
charged as for a split. There is no $\Phi$, since the parts do not add
up to the old line, and no $R$.

\paragraph{Swap}.
\[
  w = \theta + \theta_1\, m + \sum_{i=1}^{m} \bigl(w_i - G_i - L_i\bigr)
\]
for $m$ legs, where $G_i$ and $L_i$ are leg $i$'s tag and label
charges. The mistagging is one filing error, charged once as
$\theta + \theta_1 m$; each leg still pays for its numbers, so a swap
is chosen only when the numbers carry.

\paragraph{Fold and break out.} The parts pay their own ended or started
prices, less one base charge for each part after the first. A break
out's parts are also not charged for the figure the total implies, one
part per period, to the digits the footing check verifies ($\bar\sigma$,
as for a split's credit):
\[
  w_{\text{break out}} = \sum_j p_{s_j} - (k - 1)\,\kappa_S
    - \alpha \sum_k \min\bigl(\bar\sigma,\ \max_j \sigma(u_j(k))\bigr) ,
  \qquad
  w_{\text{fold}} = \sum_j p_{e_j} - (k - 1)\,\kappa_E ,
\]
plus the total's rename when its tag changes. The footing earns no
evidence credit: it is the filer's declared arithmetic and adds up by
construction.

\subsection{Constants}\label{app:constants}

The prices in \cref{sec:cost} are formulas in the constants below. The
formulas follow from the properties the prices must have; the values do
not, and were tuned. \Cref{tab:constants} gives the values the
implementation uses, set by hand against a bench of known cases and not
fitted to the corpus. They are in units of $c(v)$, and only
their ratios matter.

\begin{table}[htbp]
  \centering
  \includegraphics[width=\linewidth]{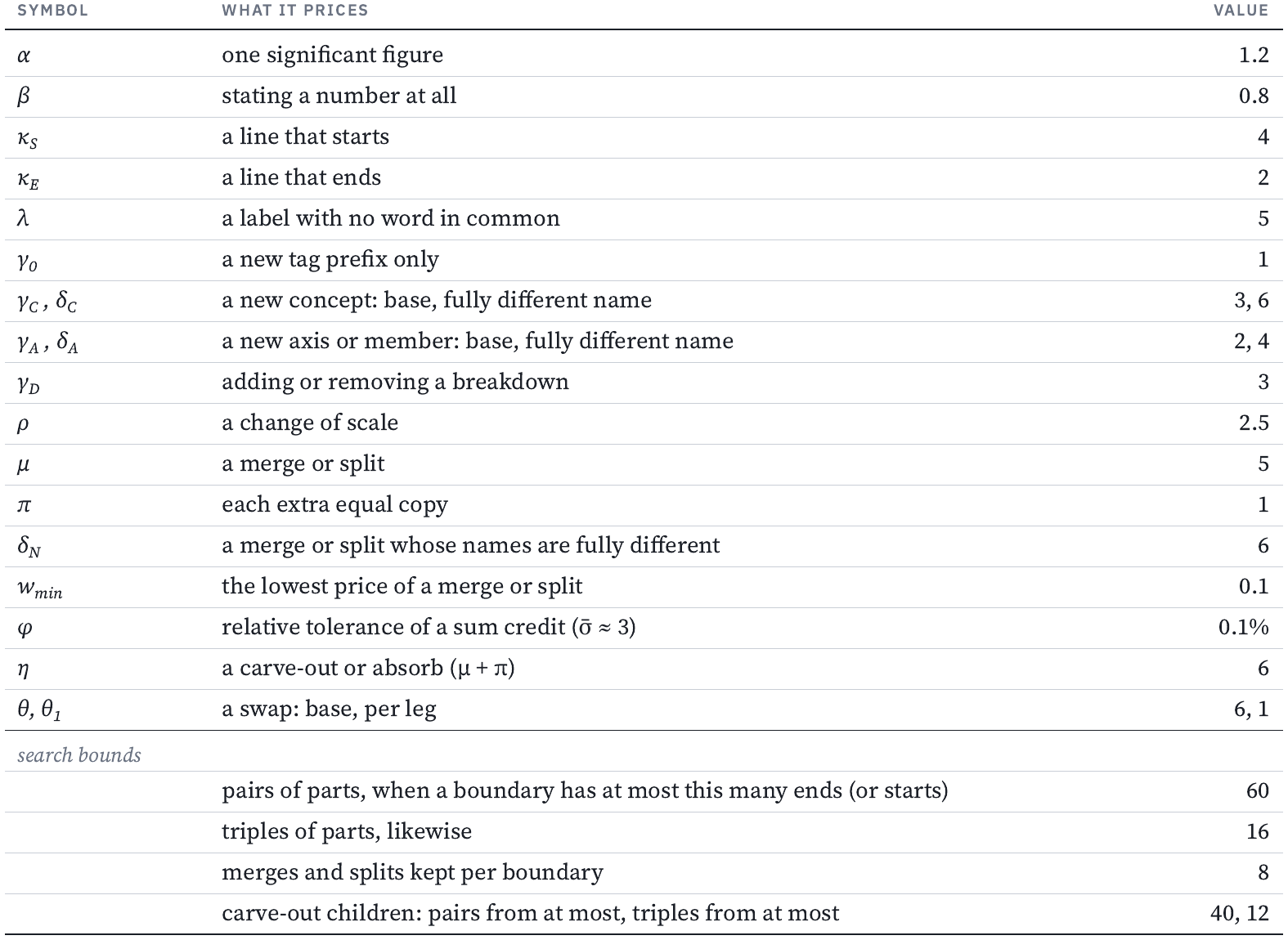}
  \caption{One working set of constants for the prices in \cref{sec:cost}, and
    the bounds on the candidate search in \cref{sec:candidates}.}
  \label{tab:constants}
\end{table}

\subsection{Revisions the random sample prompted}\label{app:revisions}

The random sample of \cref{sec:results-measured} was held out until
the method had been built and first evaluated. Its cases then prompted
three revisions, all in the rules of \cref{sec:candidates}. DXP's
change from net to gross contract positions (Appendix~\ref{app:hard})
prompted two: the sign a child takes from its balance attribute when
it adds into a different total, and the credit a carve-out or absorb
earns for the declared totals its movement explains. The layout check
of \cref{sec:results-measured} prompted the third: it found three
renames of a balance to a flow, Iron Mountain's dividends payable to a
gain on disposal and, at two companies, the change in restricted cash
to the restricted-cash balance in the cash bridge, and a line stated at
instants now never continues as one stated over periods. The same case
showed that the opposite-words rule had refused Iron Mountain's true
rename, from a tag carrying both ``current'' and ``noncurrent'' to one
carrying ``current''; a tag carrying both words of a pair now takes no
side.

%% file: sections/b-vocabulary.tex
\section{The kinds of edit, illustrated}\label{app:vocabulary}

\begin{figure}[tbp]
  \centering
  \includegraphics[width=0.85\linewidth,height=0.42\textheight,keepaspectratio]{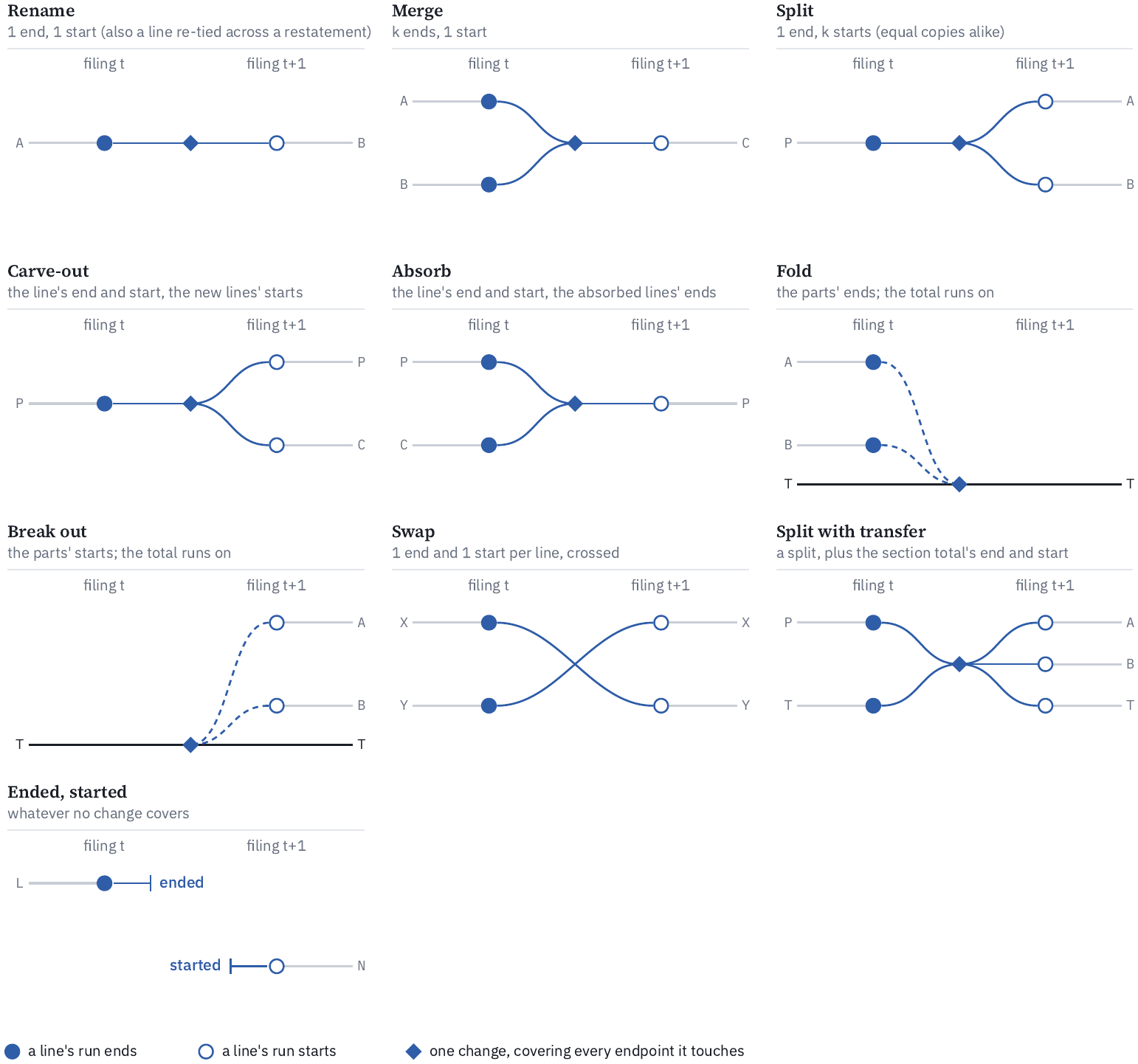}
  \caption{Each kind of edit as the strand endpoints it covers. A filled dot
    is where a line's run ends in filing $t$; an open dot is where a run starts
    in filing $t{+}1$; the diamond is one edit, and it covers every endpoint
    it touches. A fold or break out leaves its total's run untouched (dashed).
    Anything no edit covers has ended or started.}
  \label{fig:edit-shapes}
\end{figure}

Every kind of edit in \cref{tab:edits}, drawn twice: a made-up case
with round numbers, then a real one the method found, with the accession
numbers of both filings and the numbers as filed. What the
company printed is set in serif and the edit is drawn in blue (dashed
where parts of a total are shown or hidden and no number moves); the
shaded column is the year both filings report, where each edit is
checked; a restated number is circled in orange. Tags without a
prefix are \texttt{us-gaap}; some labels are shortened.
\Cref{fig:edit-shapes} first puts every kind side by side as the strand
endpoints it covers.

\paragraph{Rename, ended, started.}

\begin{figure}[tbp]
  \centering
  \includegraphics[width=0.85\linewidth,height=0.42\textheight,keepaspectratio]{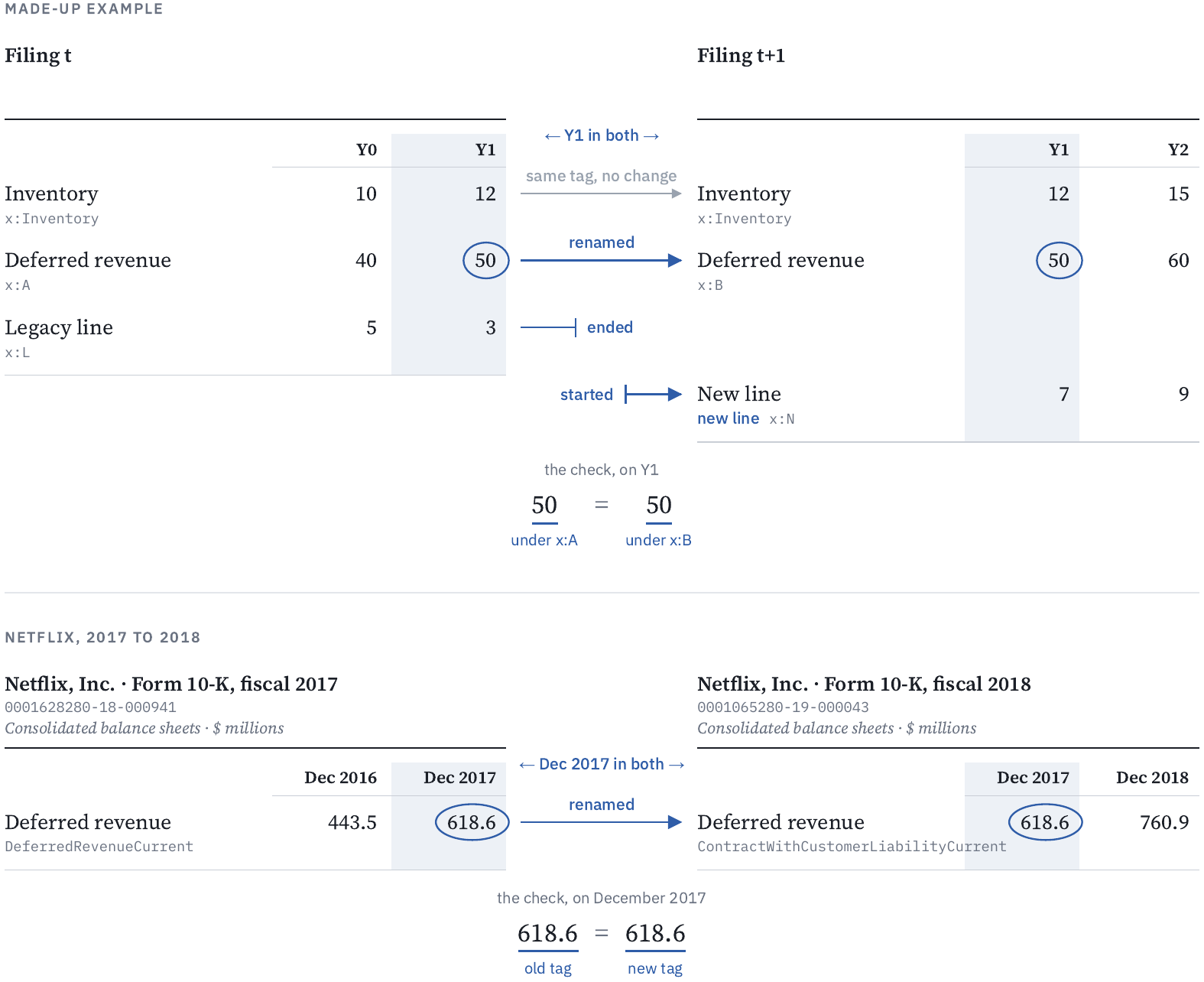}
  \caption{\edit{rename}, and lines that \edit{ended} or \edit{started}.
    Netflix's deferred revenue moves to the tag of the new revenue standard,
    with the same wording and the same December 2017 balance.}
  \label{fig:link}
\end{figure}

A \edit{rename} says that a line that stops appearing under one tag
carries on under another: a different tag, a different product or
segment member, or the same tag after the company redid its prior years.
Nothing has to add up; a rename is priced on what changed: each part of
the tag, each word of the label, and each number that does not carry
over. A line
that \edit{ended} has nothing after it, and one that \edit{started} has
nothing before it. Those two have a fixed price, and every other
explanation of a line has to be cheaper to be chosen (\cref{fig:link}).

\paragraph{Merge and split.}

\begin{figure}[tbp]
  \centering
  \includegraphics[width=0.85\linewidth,height=0.42\textheight,keepaspectratio]{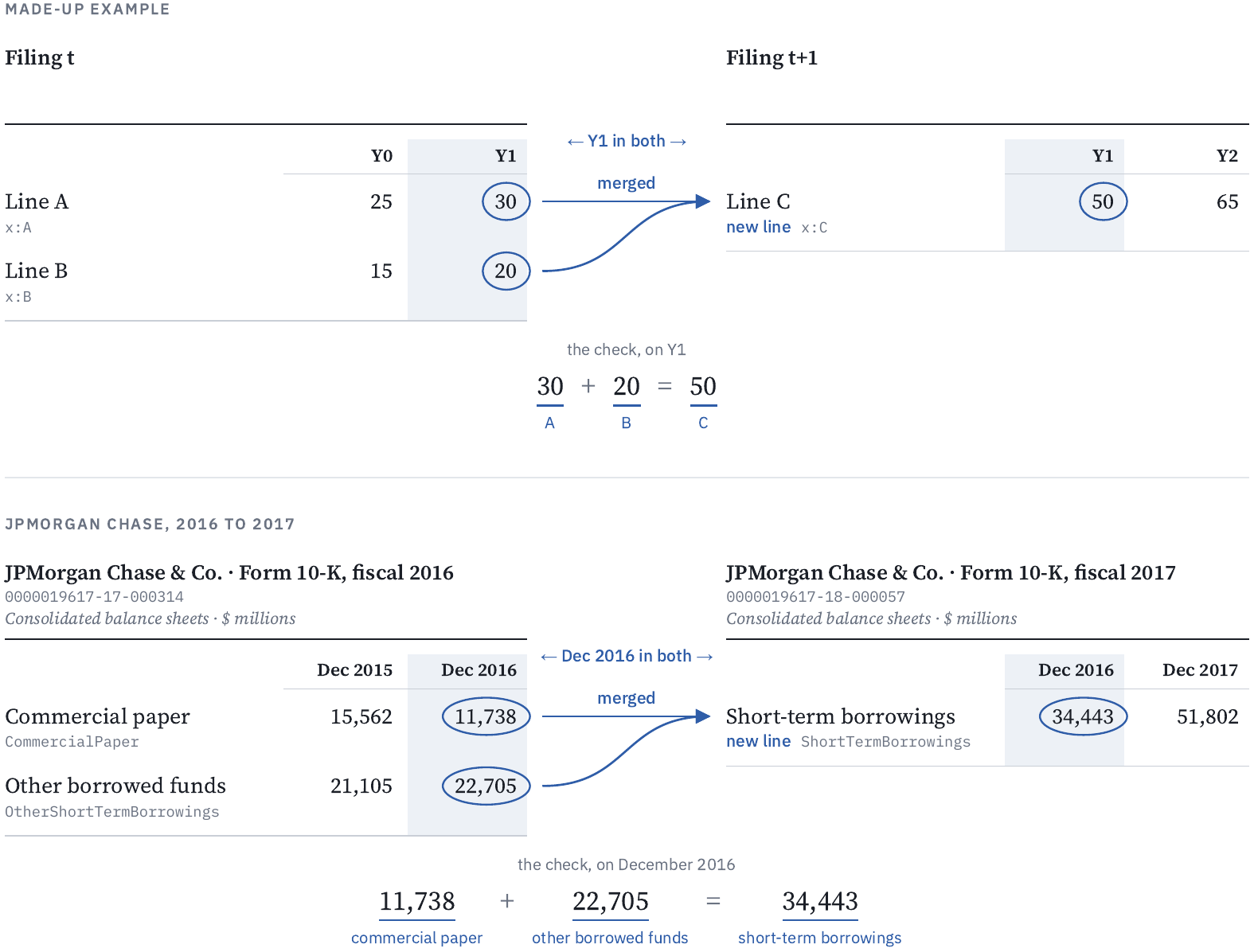}
  \caption{\edit{merge}. JPMorgan Chase shows commercial paper and other
    borrowed funds as one line from its 2017 10-K.}
  \label{fig:merge}
\end{figure}

\begin{figure}[tbp]
  \centering
  \includegraphics[width=0.85\linewidth,height=0.42\textheight,keepaspectratio]{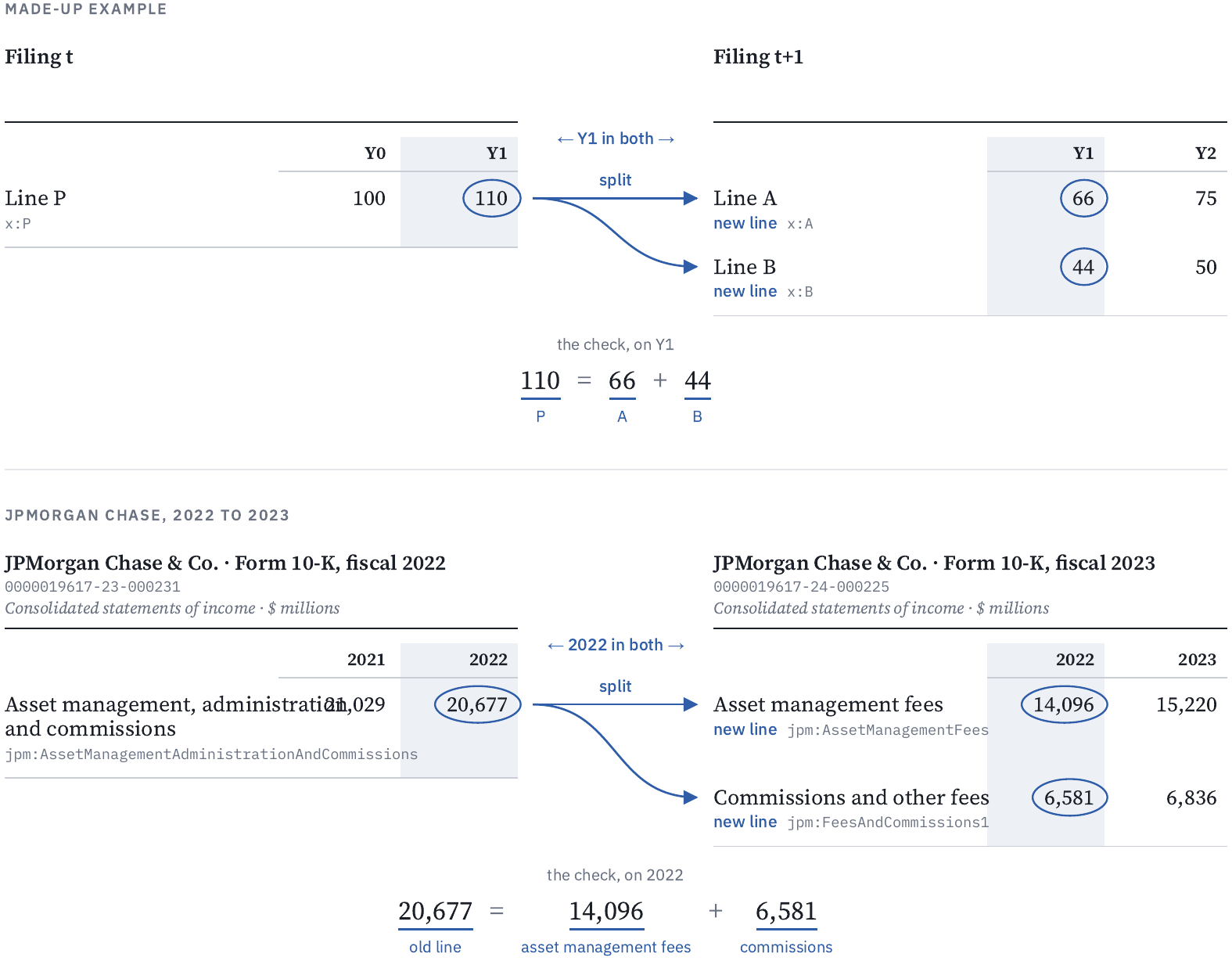}
  \caption{\edit{split}. JPMorgan Chase separates asset-management fees from
    commissions in its 2023 10-K; the sum holds on both prior years the newer
    filing reports.}
  \label{fig:split}
\end{figure}

In a \edit{merge}, several lines end and one new line starts whose
prior years are their total; a \edit{split} is the reverse
(\cref{fig:merge,fig:split}). The sum is checked on every year both
filings report, to the precision the filing states. Every part has to
matter: a merge that would still add up without one of its lines is
refused, which keeps out coincidences.

\paragraph{Equal copies.}

\begin{figure}[tbp]
  \centering
  \includegraphics[width=0.85\linewidth,height=0.42\textheight,keepaspectratio]{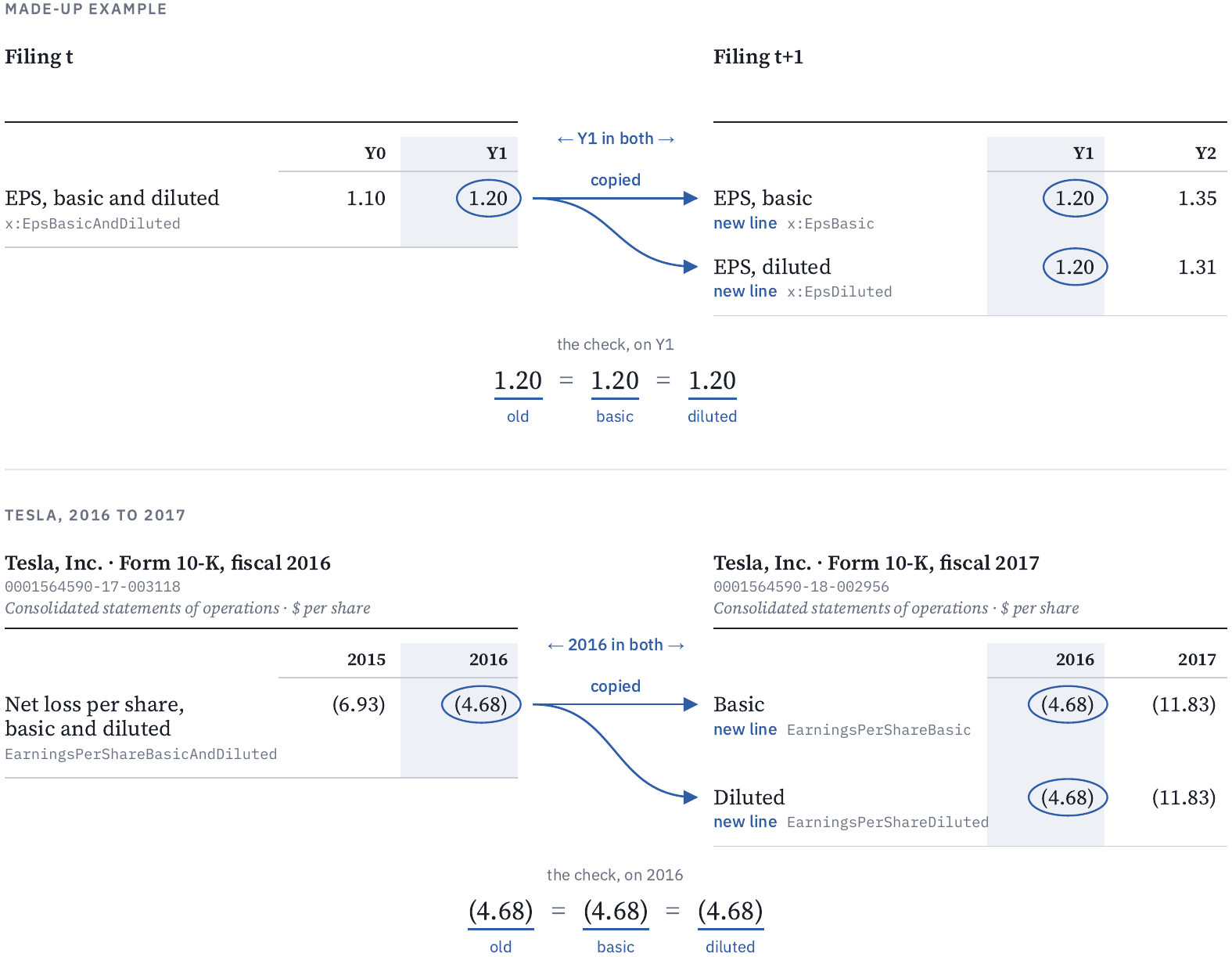}
  \caption{\edit{equal copies}, splitting and merging. Tesla reports basic and
    diluted loss per share separately from its 2017 10-K.}
  \label{fig:eqsplit}
\end{figure}

Some lines repeat rather than add up. When a company with nothing to
dilute its shares reports one earnings-per-share line, ``basic and
diluted'', and later reports two, each new line equals the old one;
treating it as a split would double it (\cref{fig:eqsplit}). The reverse
is rarer, one case on the specimen set, at Kraft Heinz, and 28 on the
random sample, so that direction is drawn with made-up numbers only.

\paragraph{Carve-out and absorb.}

\begin{figure}[tbp]
  \centering
  \includegraphics[width=0.85\linewidth,height=0.42\textheight,keepaspectratio]{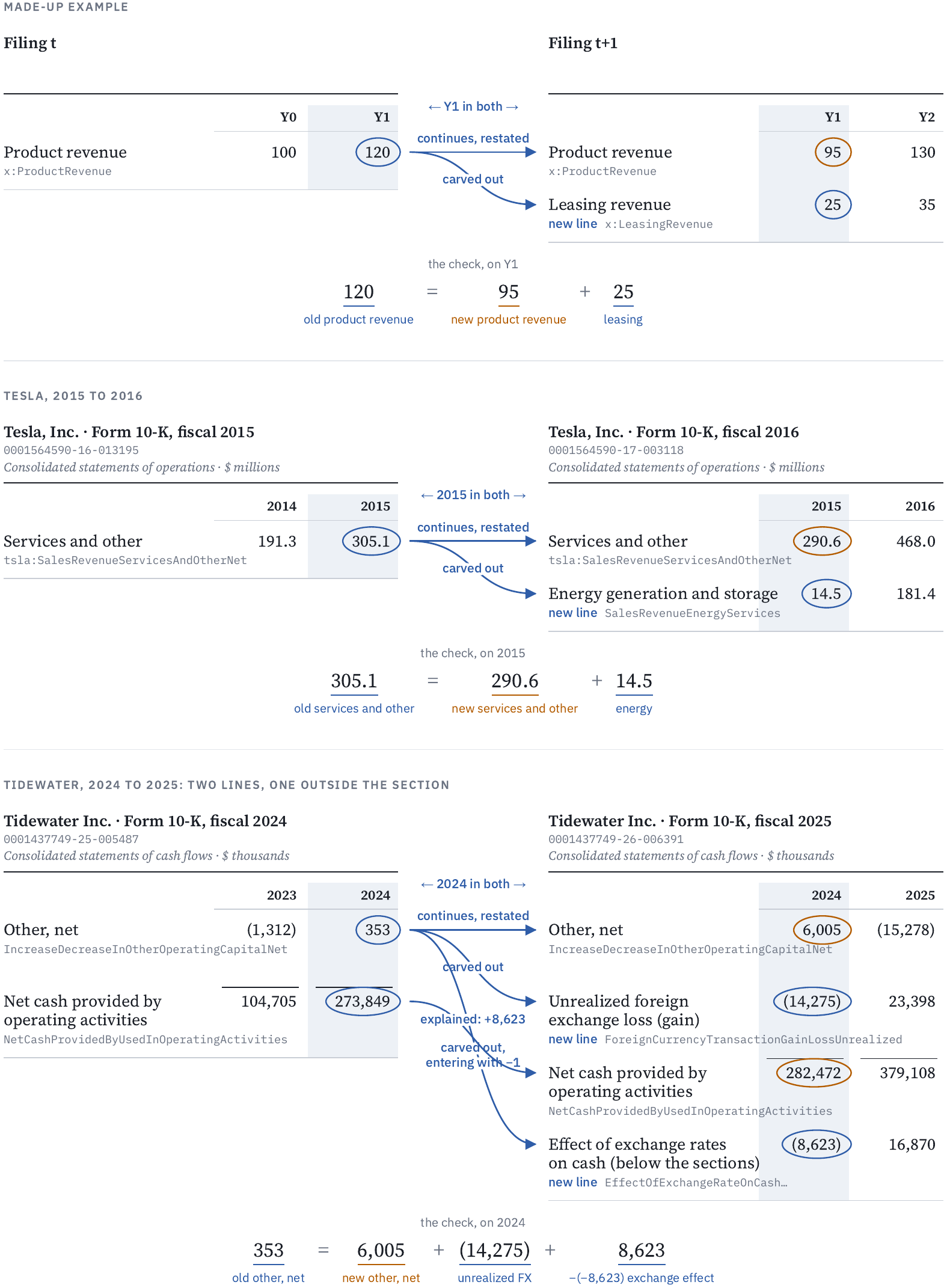}
  \caption{\edit{carve-out}. Tesla takes energy generation and storage out
    of services and other in its 2016 10-K and restates 2015 to match.
    Tidewater carves two lines out of ``Other, net'' in its 2025 10-K,
    one of them below the operating section; operating cash flow moves
    by exactly what left it.}
  \label{fig:carve}
\end{figure}

\begin{figure}[tbp]
  \centering
  \includegraphics[width=0.85\linewidth,height=0.42\textheight,keepaspectratio]{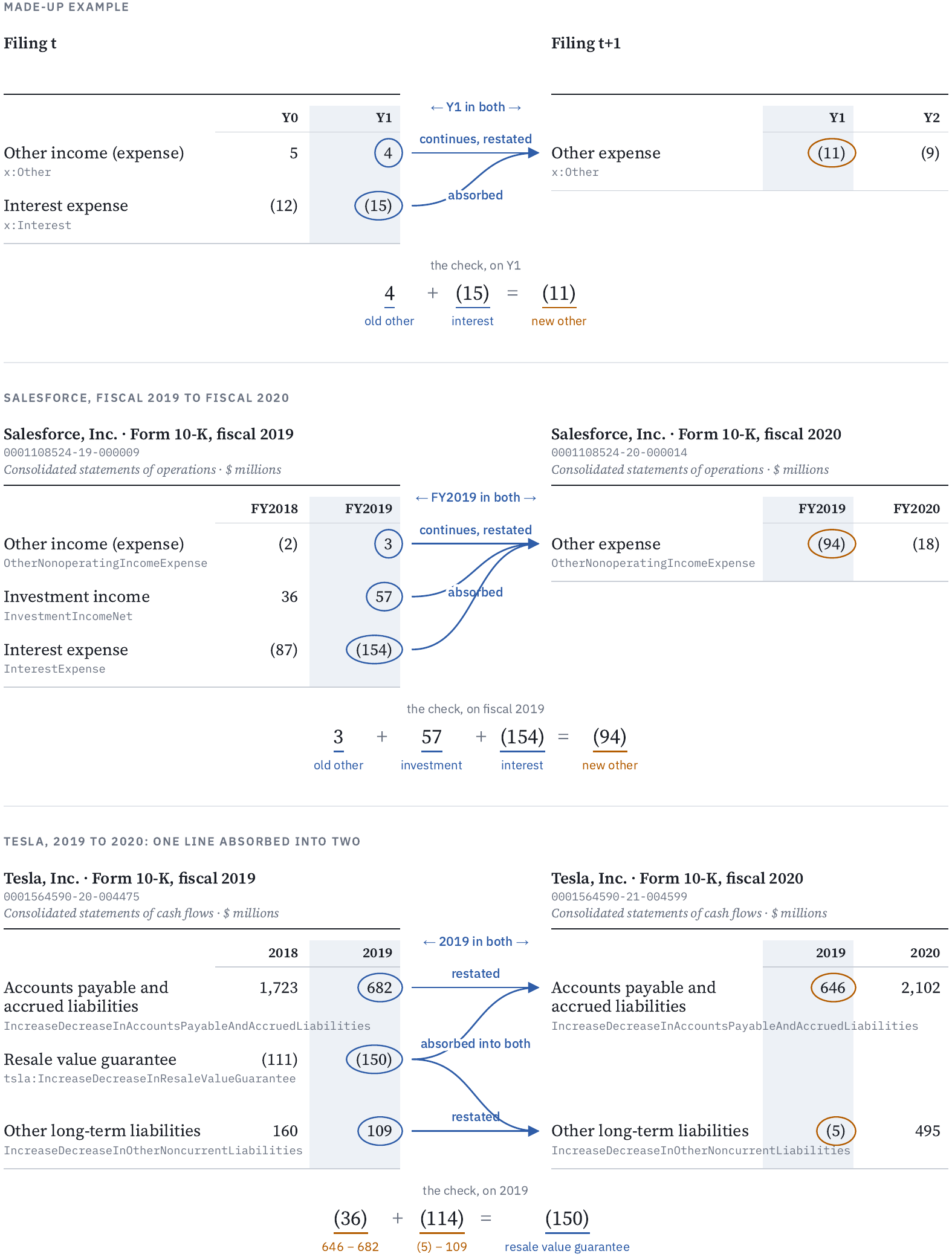}
  \caption{\edit{absorb}. Salesforce moves investment income and interest
    expense into other expense in its fiscal 2020 10-K: one edit, two lines
    absorbed. Tesla's resale value guarantee line is absorbed into two
    continuing lines at once.}
  \label{fig:absorb}
\end{figure}

``Reclassified to conform to the current period presentation'' most
often means one of these two. In a \edit{carve-out}, a line continues, a
new line starts, and the continuing line's prior years are restated down
by exactly the new line's amount (\cref{fig:carve}). An \edit{absorb} is
the reverse: lines end, and a continuing line's prior years are restated
up to include them (\cref{fig:absorb}). Signs follow the statement's own
arithmetic, so interest expense moved into ``other income'' counts
negatively, and a line that leaves for the other side of the balance
sheet, or for below the sections of a cash flow statement, enters with
the sign its balance attribute gives. Several lines absorbed at once are
one edit, since none explains the restatement on its own; one line can
be absorbed into several, when their restatements add up to it. A total
above the parent or a child that moved by exactly what the edit moved
into it is explained by the edit, as Tidewater's operating cash flow is
in \cref{fig:carve}.

\paragraph{Fold and break out.}

\begin{figure}[tbp]
  \centering
  \includegraphics[width=0.85\linewidth,height=0.42\textheight,keepaspectratio]{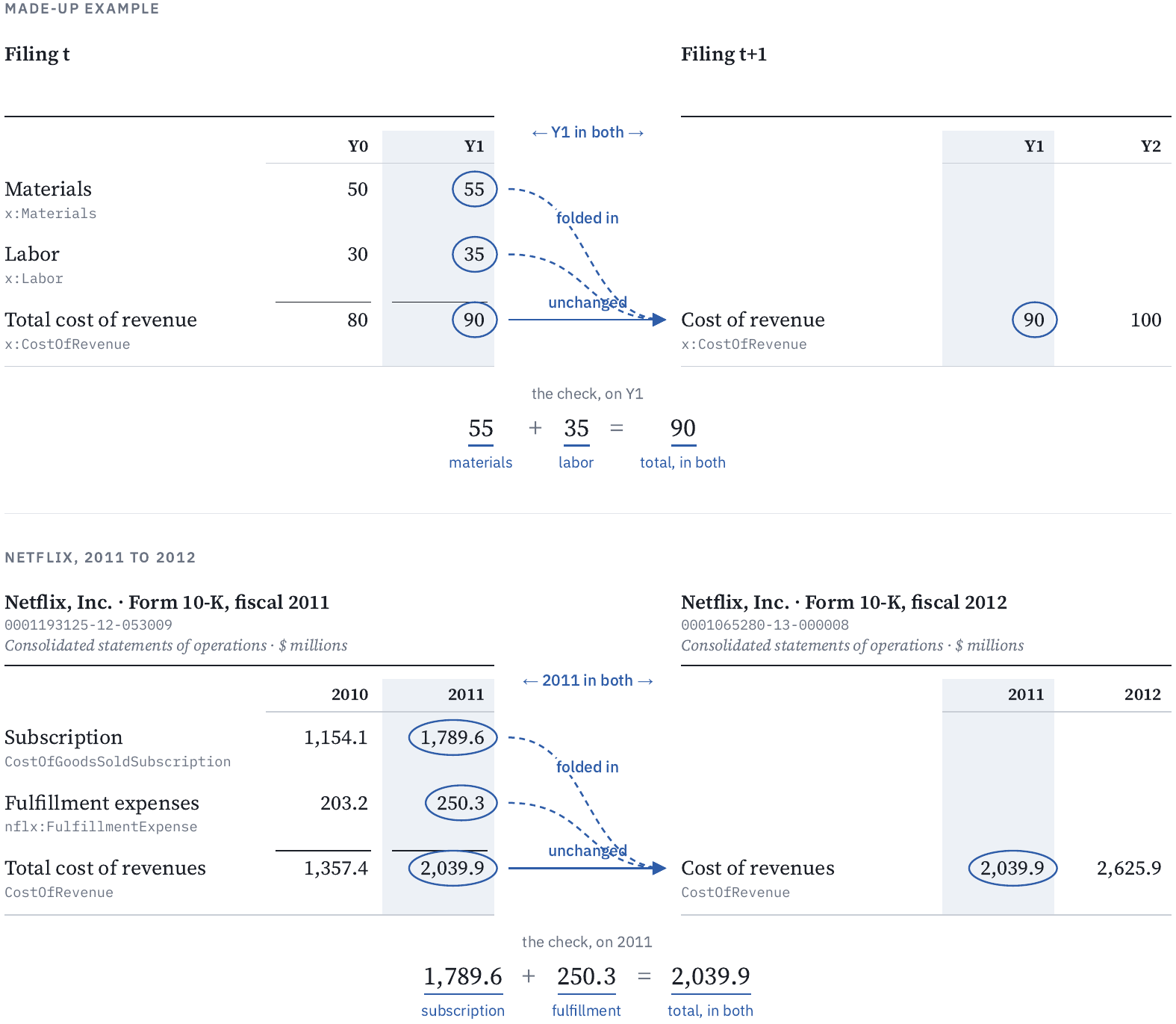}
  \caption{\edit{fold}. Netflix stops splitting cost of revenues into
    subscription and fulfillment costs in its 2012 10-K.}
  \label{fig:fold}
\end{figure}

\begin{figure}[tbp]
  \centering
  \includegraphics[width=0.85\linewidth,height=0.42\textheight,keepaspectratio]{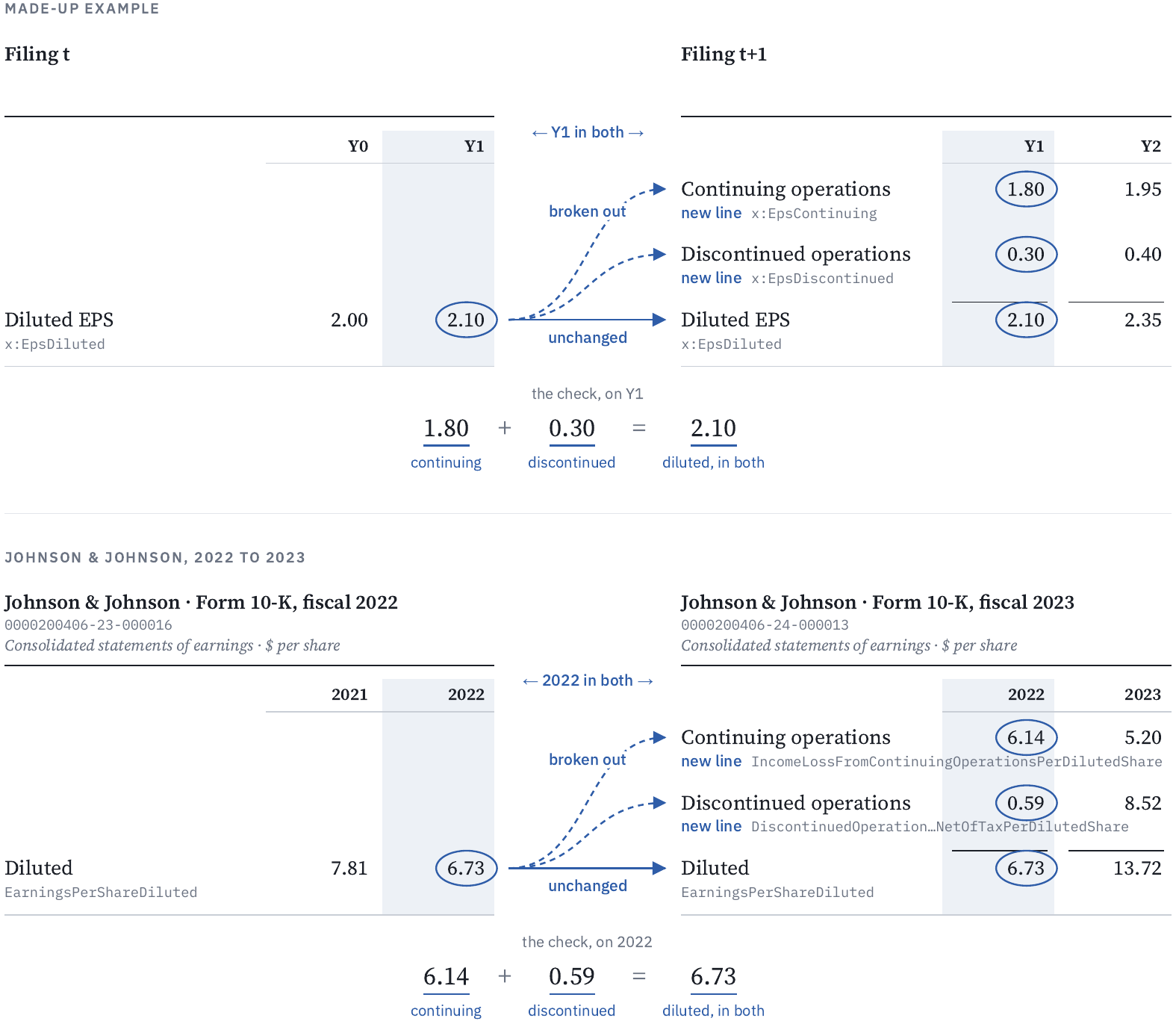}
  \caption{\edit{break out}. Johnson \& Johnson splits diluted earnings per
    share into continuing and discontinued operations after separating
    Kenvue.}
  \label{fig:breakout}
\end{figure}

In a \edit{fold}, the parts of a total stop being shown, but no number
moves: the parts were always inside the total, and the total does not
change (\cref{fig:fold}). In a \edit{break out}, a total that continues
starts showing its parts (\cref{fig:breakout}). Neither restates
anything, so each is checked by the total standing still and by the
filing's own arithmetic adding the parts up to it, in whichever filing
shows the parts.

\paragraph{Swap.}

\begin{figure}[tbp]
  \centering
  \includegraphics[width=0.85\linewidth,height=0.42\textheight,keepaspectratio]{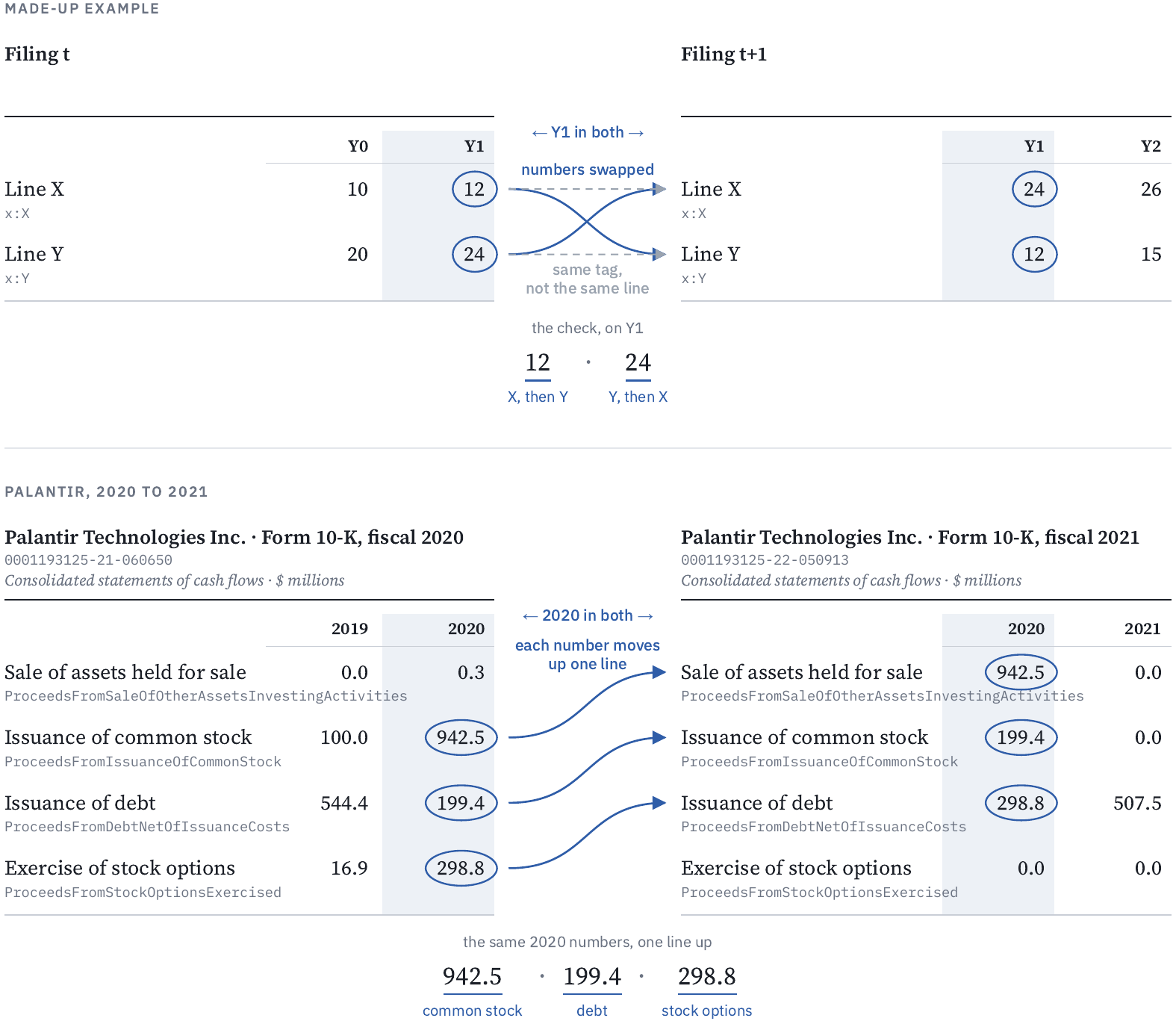}
  \caption{\edit{swap}. In the made-up example, lines X and Y keep their tags,
    but each tag now carries the other line's numbers. In Palantir's 2021 10-K
    each 2020 financing figure appears under the tag of the line above it.}
  \label{fig:swap-appendix}
\end{figure}

A \edit{swap} says the company put a run of filings' numbers on the
wrong tags: the numbers on each tag belong to another tag's line
(\cref{fig:swap-appendix}; Appendix~\ref{sec:palantir} follows
Palantir's case to its correction). It is considered only when the
lines trade places completely, so that every line that lost its numbers
got another line's, and chosen only when naming the mistake makes the
whole history cheaper. A tag that kept its name but not its meaning is
then two lines, one before the swap and one after.

\paragraph{Restatement and transfer.}

\begin{figure}[tbp]
  \centering
  \includegraphics[width=0.85\linewidth,height=0.42\textheight,keepaspectratio]{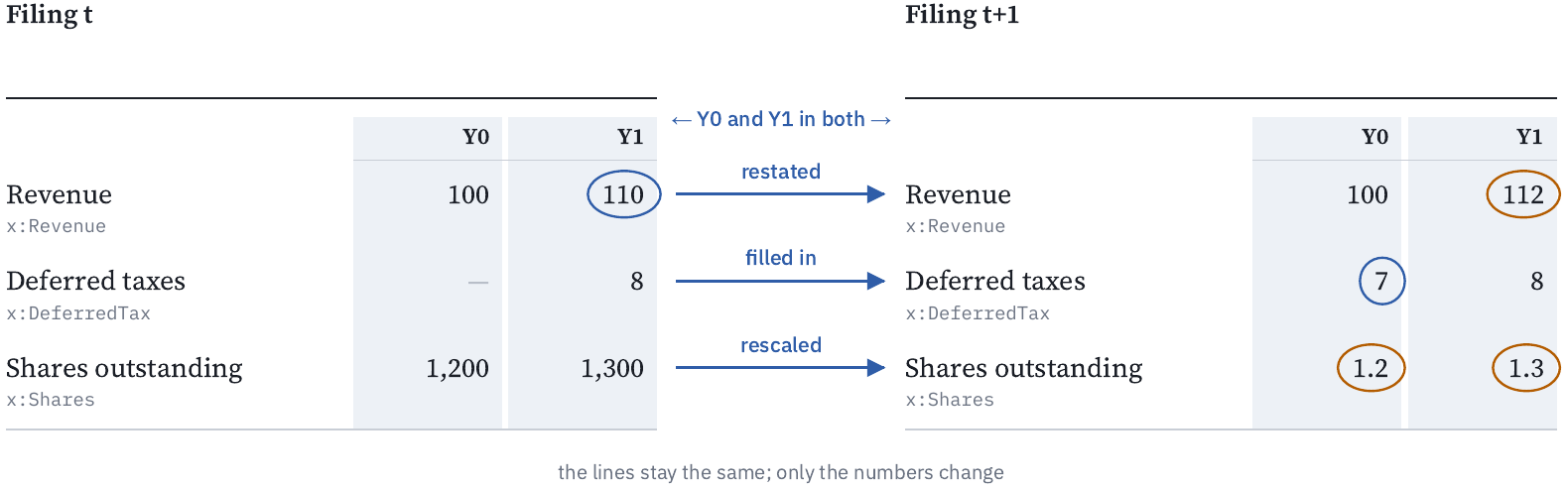}
  \caption{Edits to numbers only, with made-up numbers.}
  \label{fig:value}
\end{figure}

\begin{figure}[tbp]
  \centering
  \includegraphics[width=0.85\linewidth,height=0.42\textheight,keepaspectratio]{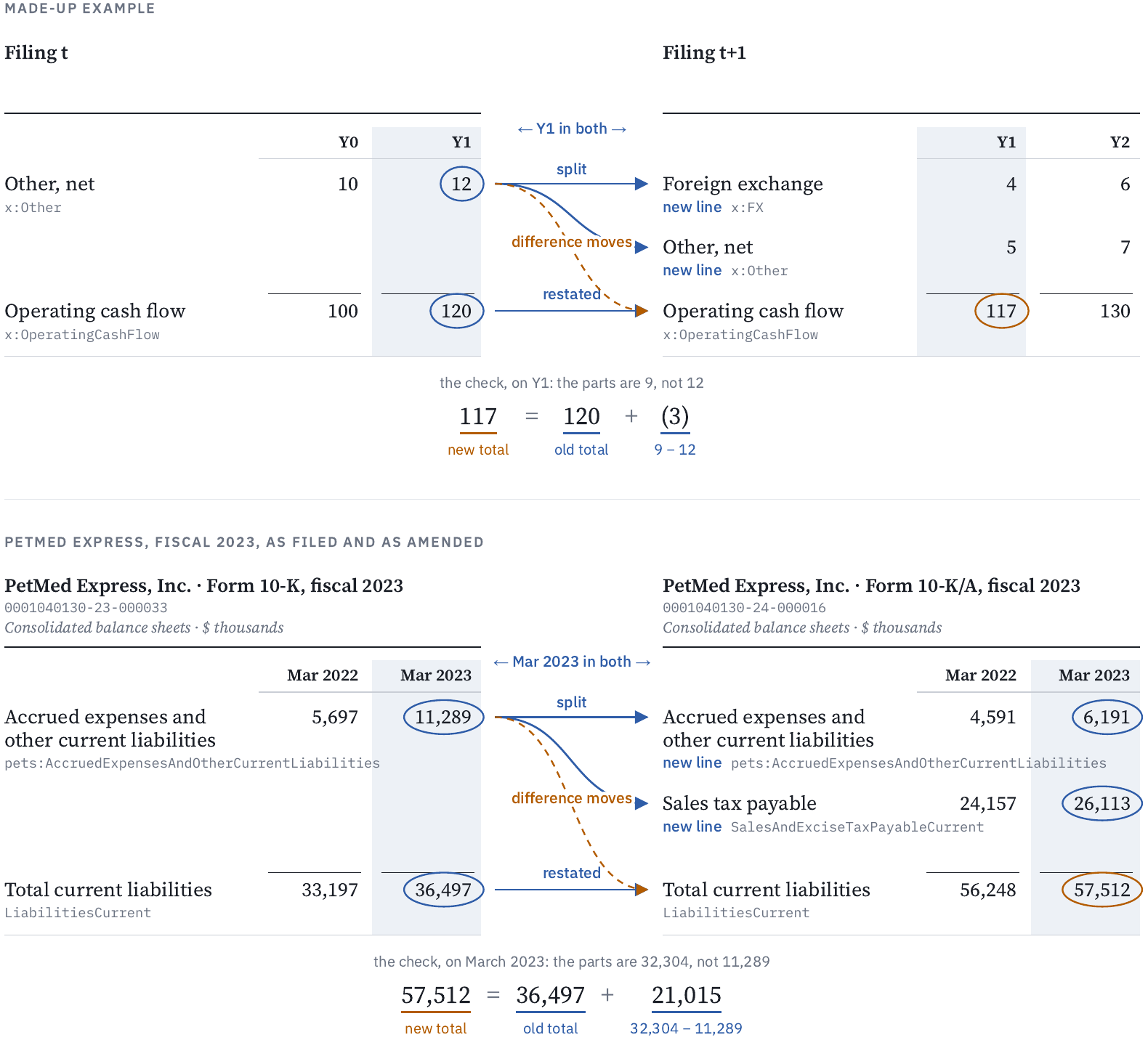}
  \caption{\edit{split} with \edit{transfer}. PetMed Express's amended
    fiscal 2023 10-K shows sales tax payable as its own line. The parts
    exceed the old accrued-liabilities line, and total current liabilities
    for March 2023 is restated by exactly the difference.}
  \label{fig:split-transfer}
\end{figure}

Some edits touch only numbers, on a line whose identity is settled
(\cref{fig:value}). A \edit{restatement} is a prior year's number that
differs from the one first filed. A later filing can also fill in a year
an earlier one left blank, or report a whole series at a different
scale. A \edit{transfer} goes with a split whose parts do not add up to
the old line: the difference moved to another section, whose total was
restated by exactly that amount (\cref{fig:split-transfer}).

%% file: sections/c-examples.tex
\section{Examples}\label{sec:solutions}

A made-up example taken through the four steps of \cref{alg:method}
(\cref{sec:how}), then the four cases of \cref{sec:results-measured} in
full: the facts, the strands, the candidates, the solution, and the rows
that follow. The four are from one run of the implementation, with
nothing set by hand for any company; prices are in the units of
\cref{sec:cost}, numbers in millions.

\subsection{A made-up example}\label{sec:worked}

\begin{figure}[htbp]
  \centering
  \includegraphics[width=\linewidth]{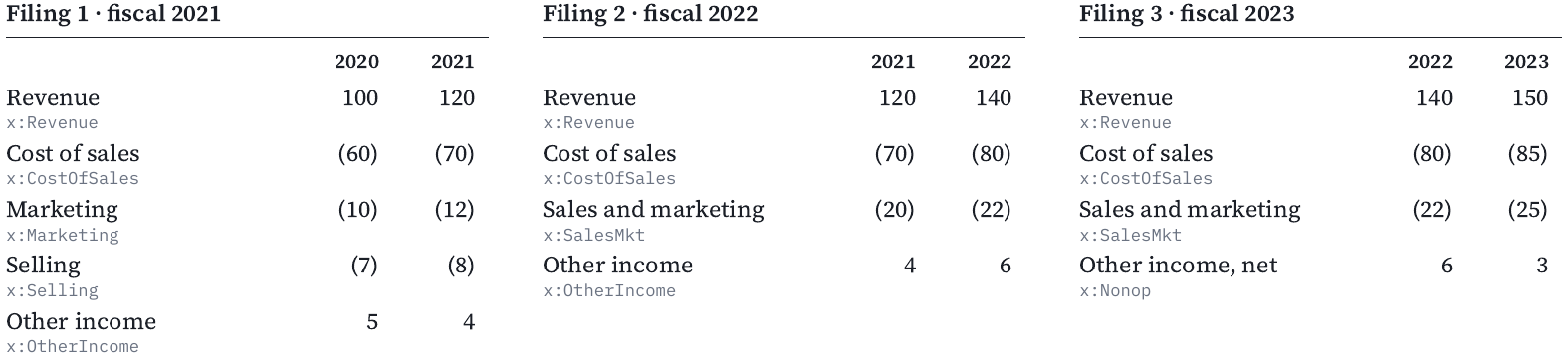}
  \caption{The worked example: three made-up annual reports, each with its lines, their tags, and two years of numbers.}
  \label{fig:worked-filings}
\end{figure}

\paragraph{The filings.} Three annual reports
(\cref{fig:worked-filings}), each showing the current year and the year
before. Marketing and Selling are replaced by Sales and marketing in the second
filing, and Other income is retagged as Other income, net in the third.

\begin{figure}[htbp]
  \centering
  \includegraphics[width=\linewidth]{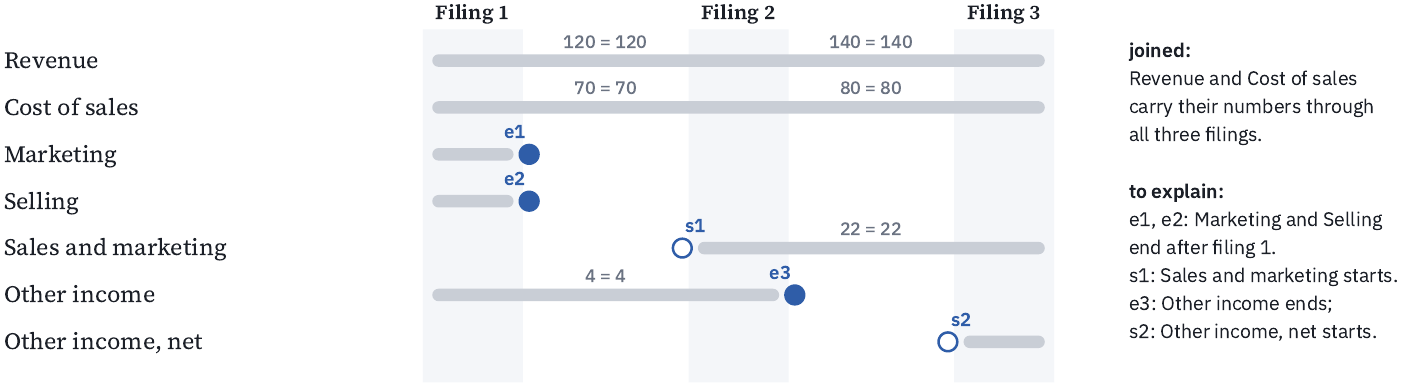}
  \caption{Step 1: strands. A bar is a strand, with the number that joined it across each boundary; a filled dot is where a strand ends ($e$) and an open dot where one starts ($s$).}
  \label{fig:worked-strands}
\end{figure}

\paragraph{Step 1, strands.} Revenue and cost of sales carry their numbers through all three
filings, so each is one strand (\cref{fig:worked-strands}). Everything
else leaves endpoints: Marketing and Selling end after the first filing
($e_1$, $e_2$) and Sales and marketing starts in the second ($s_1$);
Other income ends after the second ($e_3$) and Other income, net starts
in the third ($s_2$).

\begin{figure}[htbp]
  \centering
  \includegraphics[width=\linewidth]{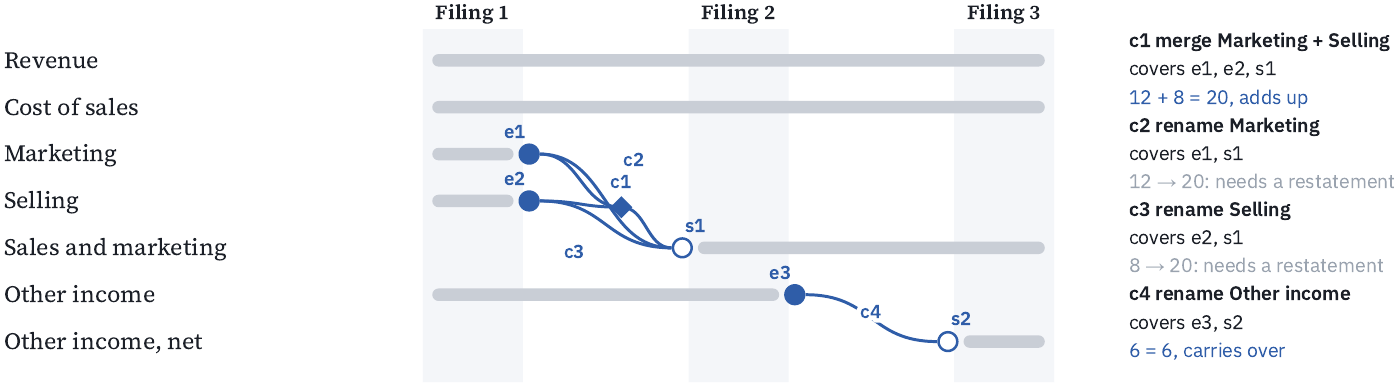}
  \caption{Step 2: candidates. Each candidate joins the endpoints it would explain; a diamond marks an edit that covers more than two.}
  \label{fig:worked-candidates}
\end{figure}

\paragraph{Step 2, candidates.} Four candidates cover those endpoints (\cref{fig:worked-candidates}):
$c_1$, a merge of Marketing and Selling into Sales and marketing
($12 + 8 = 20$, which adds up); $c_2$ and $c_3$, a rename of either line
into Sales and marketing, each needing a restatement ($12 \to 20$,
$8 \to 20$); and $c_4$, a rename of Other income ($6 = 6$, carried
over).

\begin{figure}[htbp]
  \centering
  \includegraphics[width=\linewidth]{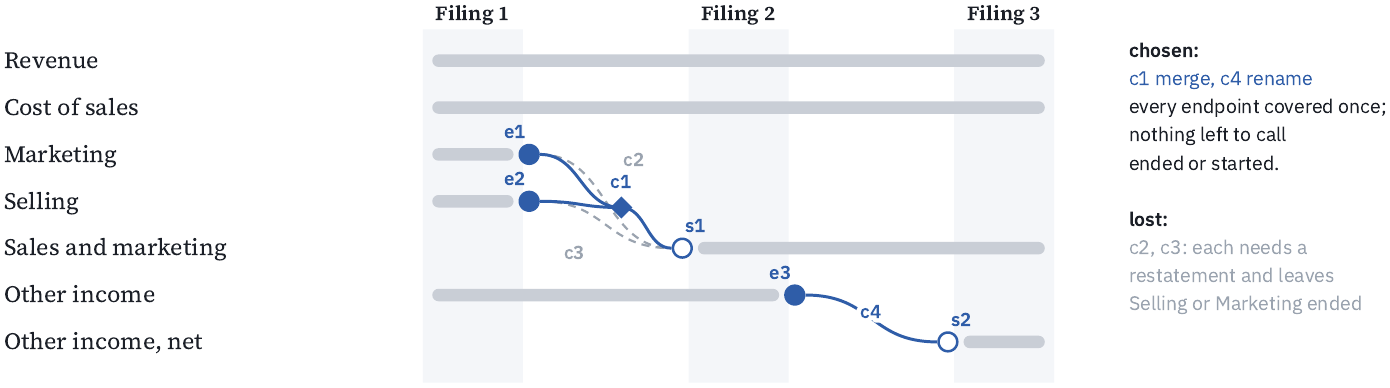}
  \caption{Step 3: the choice. $c_1$ and $c_4$ are chosen (solid); $c_2$ and $c_3$ lose (dashed grey), because each would share $s_1$ with the merge.}
  \label{fig:worked-choice}
\end{figure}

\paragraph{Step 3, choose the edits.} $c_1$, $c_2$ and $c_3$ all need
$s_1$, so at most one can be chosen. The merge adds up, pays nothing for
numbers and leaves no endpoint uncovered; either rename would pay for a
restatement and leave the other line to pay for having ended. $c_4$
shares no endpoint with anything and is decided on its own. The chosen
set, $c_1$ and $c_4$, covers every endpoint once
(\cref{fig:worked-choice}).

\begin{figure}[htbp]
  \centering
  \includegraphics[width=\linewidth]{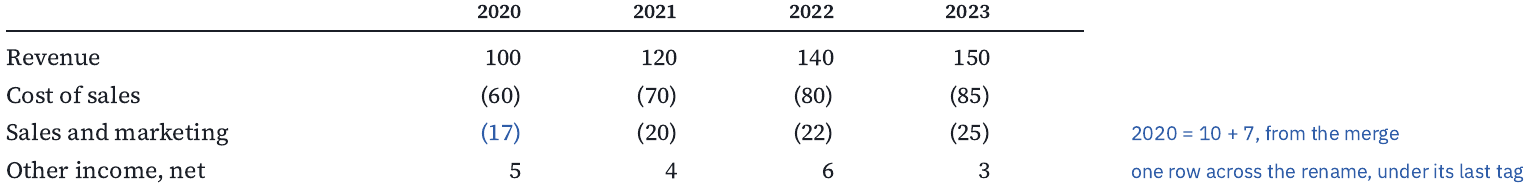}
  \caption{Step 4: the rows that follow from the choice. The figure in blue is not in any filing; the merge supplies it.}
  \label{fig:worked-rows}
\end{figure}

\paragraph{Step 4, rows.} Sales and marketing gets a 2020 figure of $10 + 7 = 17$ from the merge
(\cref{fig:worked-rows}), and Other income is one row across its
rename, under the tag it ends on.

\subsection{Apple}\label{sec:apple}

\paragraph{The facts.} \Cref{tab:apple-facts}: three lines are gone from
the fiscal 2022 report, and the two ``Other'' lines changed by exactly what
they carried ($754 = 880 + (126)$; $(909) = (791) + (210) + 92$). The
report does not say so~\cite{aapl-10k-2021,aapl-10k-2022}.

\begin{table}[htbp]
  \centering
  \includegraphics[width=\linewidth]{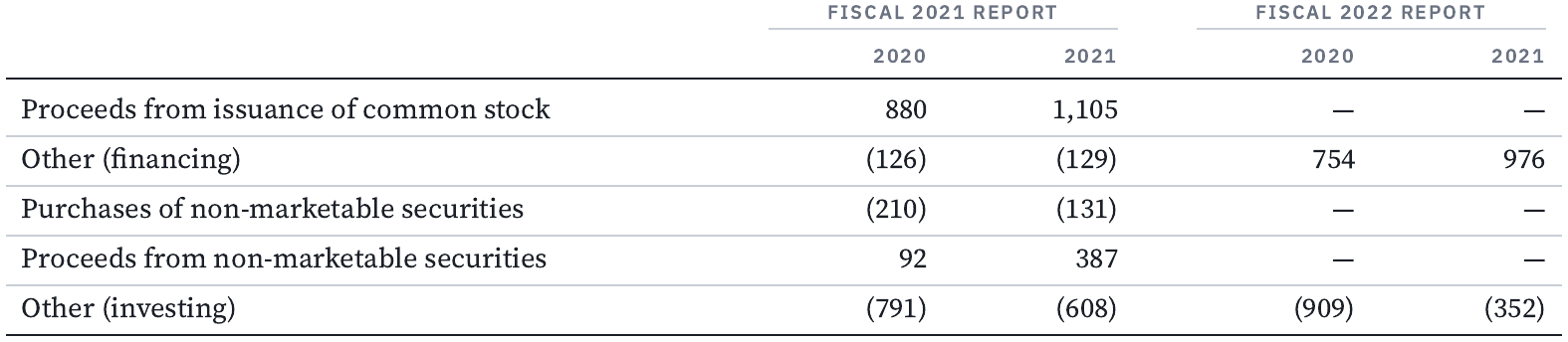}
  \caption{Five lines of Apple's cash flow statement for the same two years,
    as two consecutive reports gave them~\cite{aapl-10k-2021,aapl-10k-2022}.
    Three lines are gone; the two ``Other'' lines changed by exactly their
    amounts.}
  \label{tab:apple-facts}
\end{table}

\begin{figure}[htbp]
  \centering
  \includegraphics[width=\linewidth]{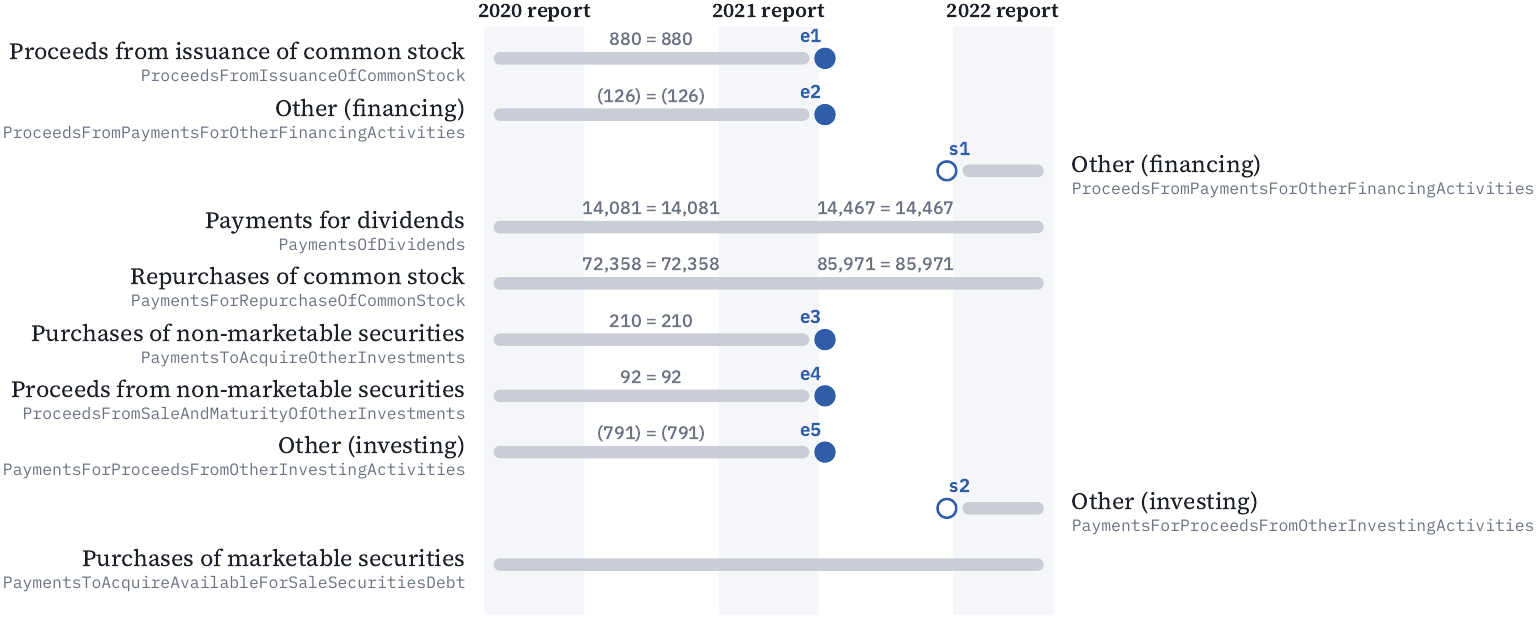}
  \caption{Apple's cash flow lines as strands across the fiscal 2020 to
    2022 reports. A filled dot is a strand that ends, an open dot one that
    starts; the number on a join is the one that carried.}
  \label{fig:apple-strands}
\end{figure}

\paragraph{Strands.} \Cref{fig:apple-strands}. The three gone lines are
strands that end ($e_1$, $e_3$, $e_4$); each ``Other'' is a strand that
ends and one that starts, since no shared number agrees ($e_2$, $s_1$;
$e_5$, $s_2$); every other line runs through. Seven endpoints.

\begin{figure}[htbp]
  \centering
  \includegraphics[width=\linewidth]{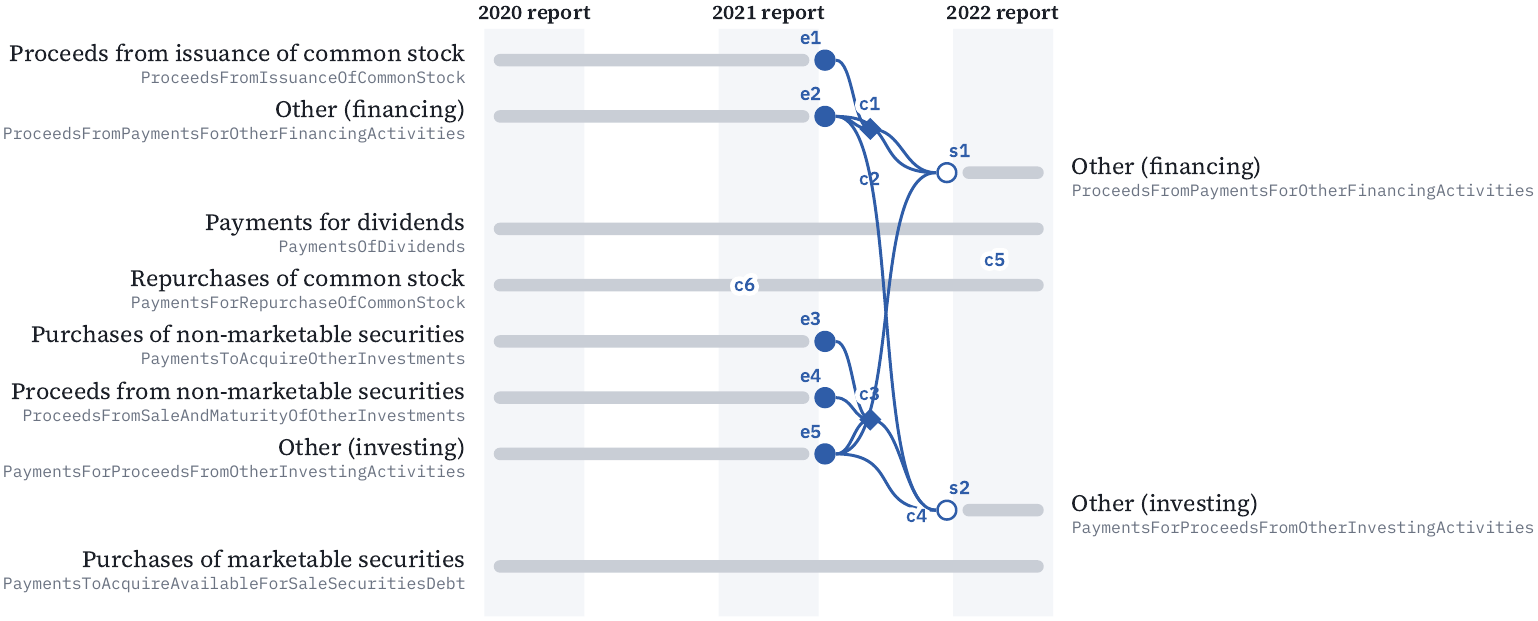}
  \caption{The six candidates at the fiscal 2021 to 2022 boundary, each
    joining the endpoints it would cover.}
  \label{fig:apple-candidates}
\end{figure}

\paragraph{Candidates.} \Cref{fig:apple-candidates}. Six:
\begin{itemize}
  \item $c_1$, $c_3$: each ``Other'' continues restated and \edit{absorbs}
    the lines that vanished ($4.81 + 6.00$; $6.01 + 6.00$).
  \item $c_2$, $c_4$: each ``Other'' continues restated; the gone lines
    ended.
  \item $c_5$, $c_6$: each ``Other'' renamed to the other.
\end{itemize}

\begin{figure}[htbp]
  \centering
  \includegraphics[width=\linewidth]{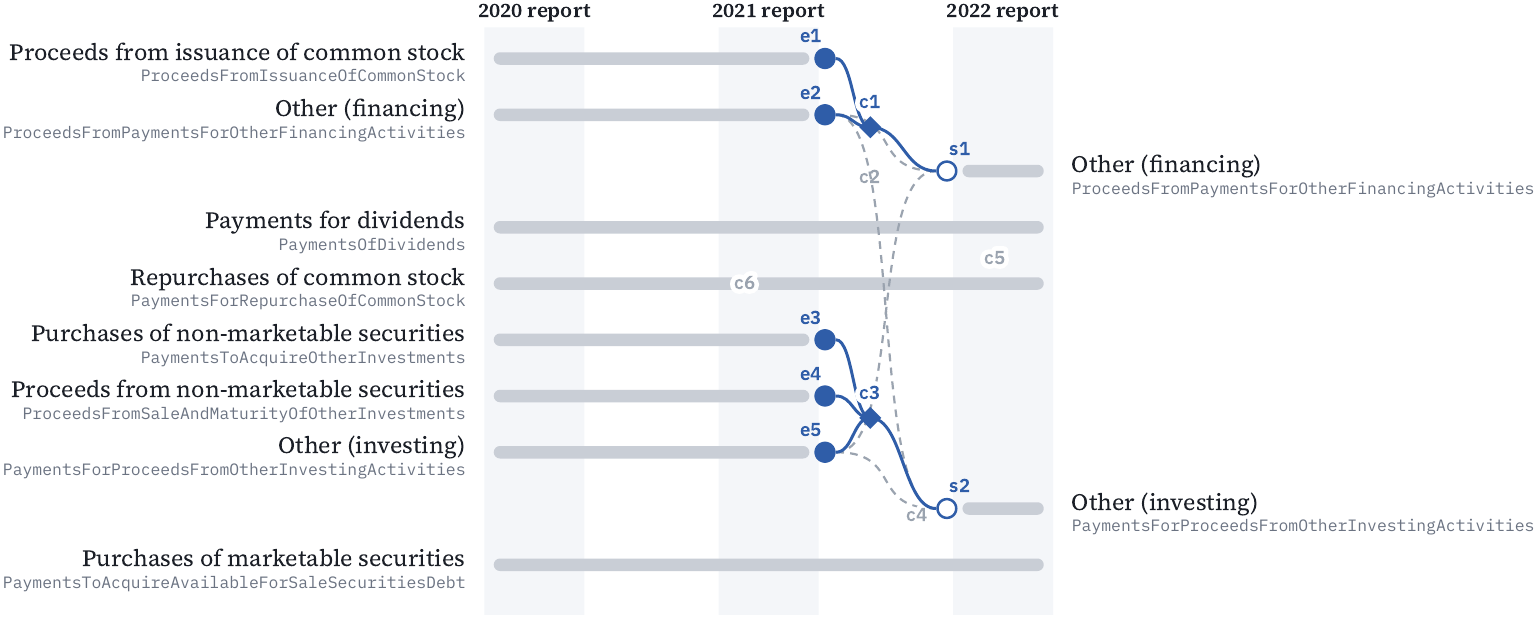}
  \caption{The solution: the two absorbs solid, the four candidates not
    chosen dashed. Every endpoint is covered.}
  \label{fig:apple-choice}
\end{figure}

\paragraph{Solution.} \Cref{fig:apple-choice}: $c_1$ and $c_3$. The
absorb's check explains the restatement, so the restated numbers are not
charged, only the act: each absorb saves $15.99$. The plain continuation
pays for the restated numbers in full and saves $13.19$. The two
renames save $2.49$ and $0.09$, and each needs an endpoint an absorb
holds.

\begin{figure}[htbp]
  \centering
  \includegraphics[width=\linewidth]{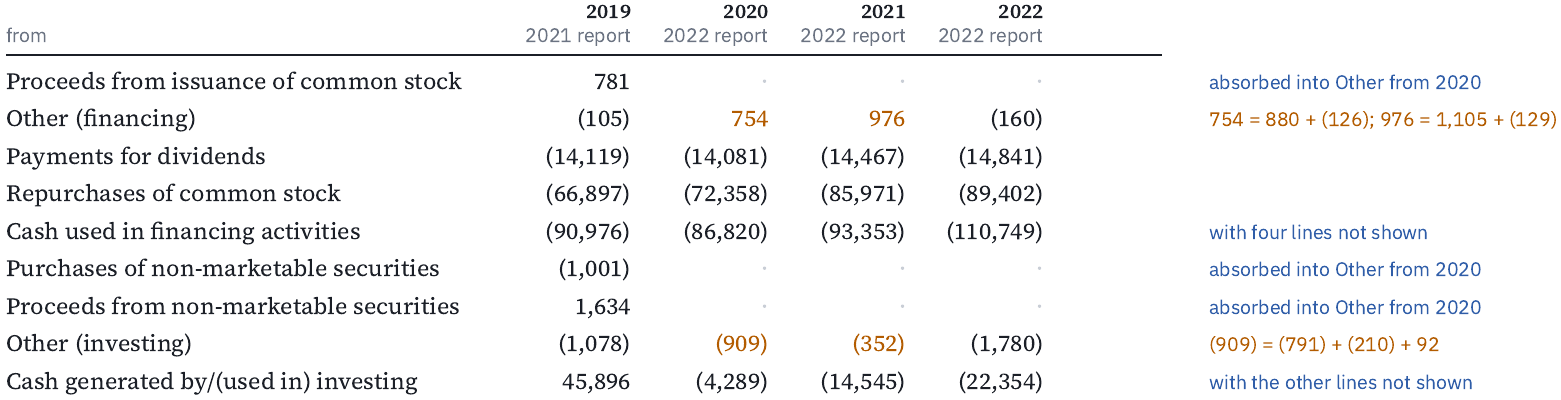}
  \caption{The rows that follow, as last reported (amber: restated across
    an absorb; the totals include lines not shown).}
  \label{fig:apple-rows}
\end{figure}

\paragraph{Rows.} \Cref{fig:apple-rows}, as last reported. 2019 is ruled
by the fiscal 2021 report and keeps the three lines; 2020 onwards by the
fiscal 2022 report, with them inside ``Other''. Every column adds up as its
report did.

\subsection{Palantir}\label{sec:palantir}

\paragraph{The facts.} \Cref{tab:palantir-facts}: in the fiscal 2021
report every column is the 2020 report's column moved down one line. The
preferred-stock line is gone and its numbers sit on the stock-options
tag; the stock-option numbers sit on the debt tag; and so on down to the
sale-of-assets tag, whose own numbers are gone. Labels stayed with their
tags. The 2022 report keeps it so and drops the stock-options tag; the
2023 report puts the stock-option numbers back on their own tag and
drops the other three~\cite{pltr-10k-2020,pltr-10k-2021,pltr-10k-2023}.

\begin{table}[htbp]
  \centering
  \includegraphics[width=\linewidth]{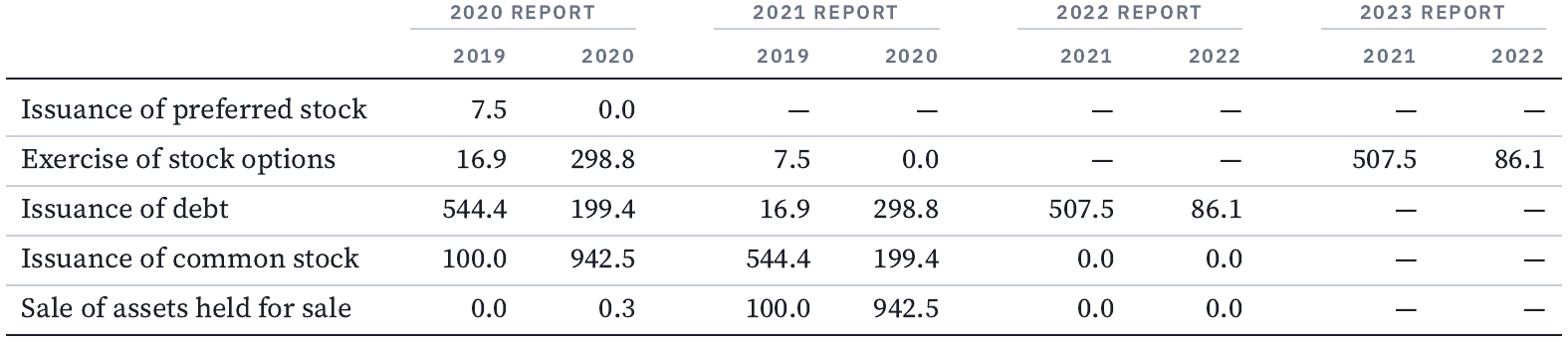}
  \caption{Five financing lines of Palantir's cash flow statement across
    four consecutive reports~\cite{pltr-10k-2020,pltr-10k-2021,pltr-10k-2023}.
    In the 2021 report each column is the 2020 report's moved down one
    line; the 2022 report keeps it so; the 2023 report moves the
    stock-option numbers back.}
  \label{tab:palantir-facts}
\end{table}

\begin{figure}[htbp]
  \centering
  \includegraphics[width=\linewidth]{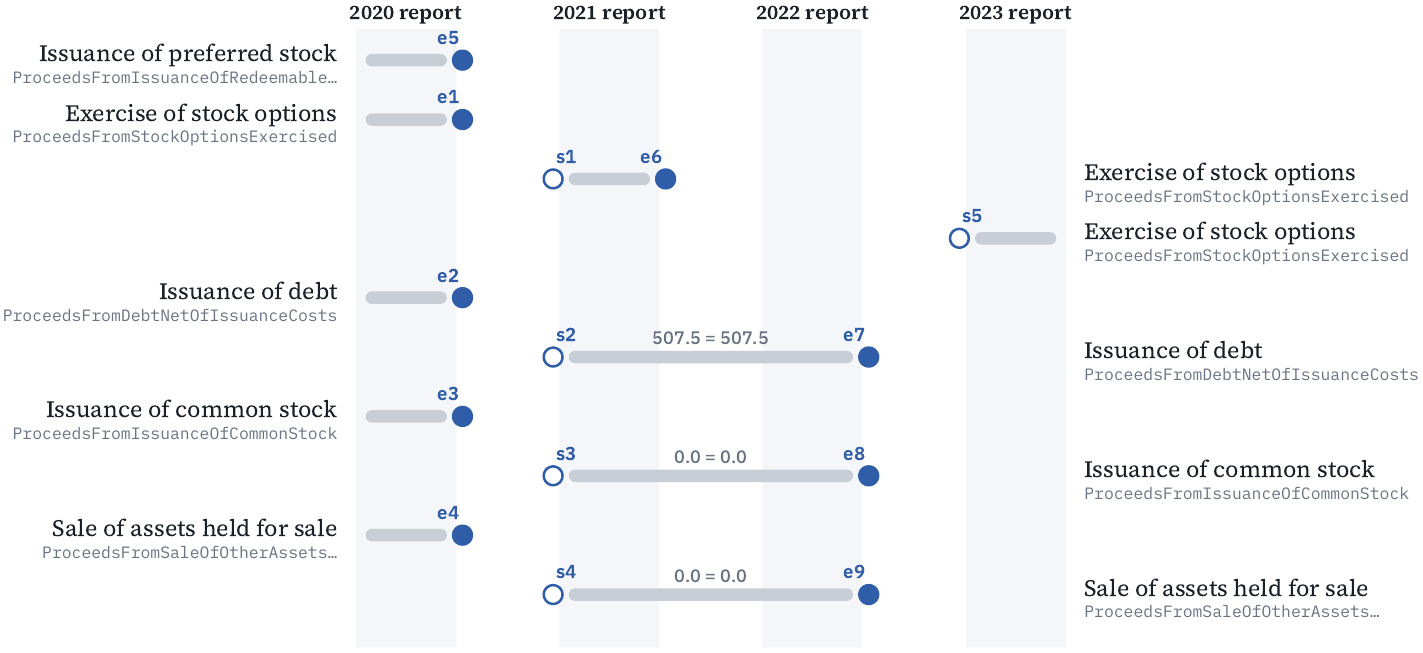}
  \caption{Palantir\textquotesingle s financing lines as strands across the fiscal 2020 to 2023 reports.}
  \label{fig:palantir-strands}
\end{figure}

\paragraph{Strands.} \Cref{fig:palantir-strands}. At the first boundary no
tag's numbers agree, so every tag is a strand that ends and one that starts
($e_1$ to $e_4$, $s_1$ to $s_4$), and the preferred-stock line ends
($e_5$). The new strands run on into the 2022 report, where the
stock-options one ends ($e_6$); at the last boundary the other three end
($e_7$, $e_8$, $e_9$) and the stock-options tag starts again ($s_5$).
Fourteen endpoints.

\begin{figure}[htbp]
  \centering
  \includegraphics[width=\linewidth]{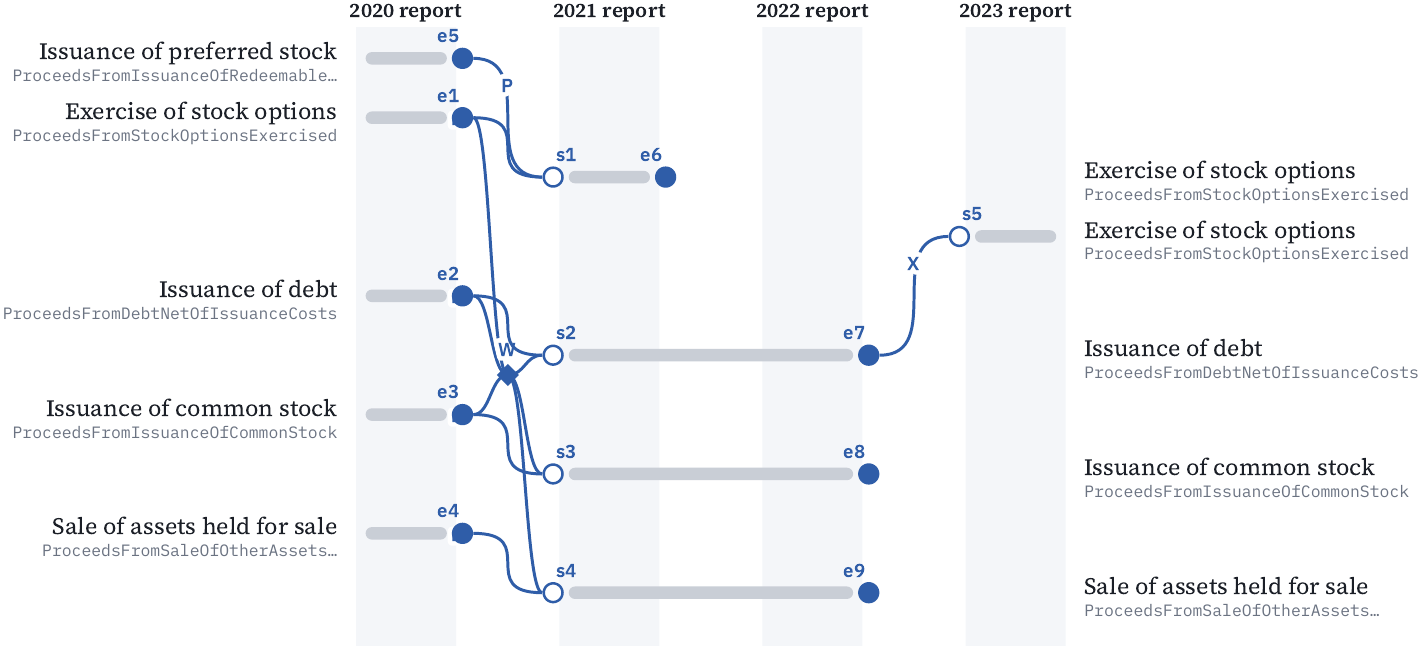}
  \caption{The candidates at the first and last boundaries.}
  \label{fig:palantir-candidates}
\end{figure}

\paragraph{Candidates.} \Cref{fig:palantir-candidates}. At the first
boundary:
\begin{itemize}
  \item W: a \edit{swap} of three legs, stock options to debt, debt to
    common stock, common stock to sale of assets. One act of $9.00$; each
    leg at the lowest price, $w_{\min}$, since along it every number
    carries verbatim.
  \item P: the preferred-stock line renamed to the stock-options tag,
    $8.17$: the fourth leg, priced as a rename because the old tag ends.
  \item $r_1$ to $r_4$: each tag continuing under its own name, with
    every number restated.
\end{itemize}
At the last boundary, X: the debt tag's strand renamed to the returning
stock-options tag, at $w_{\min}$, since $507.5$ and $86.1$ carry.

\begin{figure}[htbp]
  \centering
  \includegraphics[width=\linewidth]{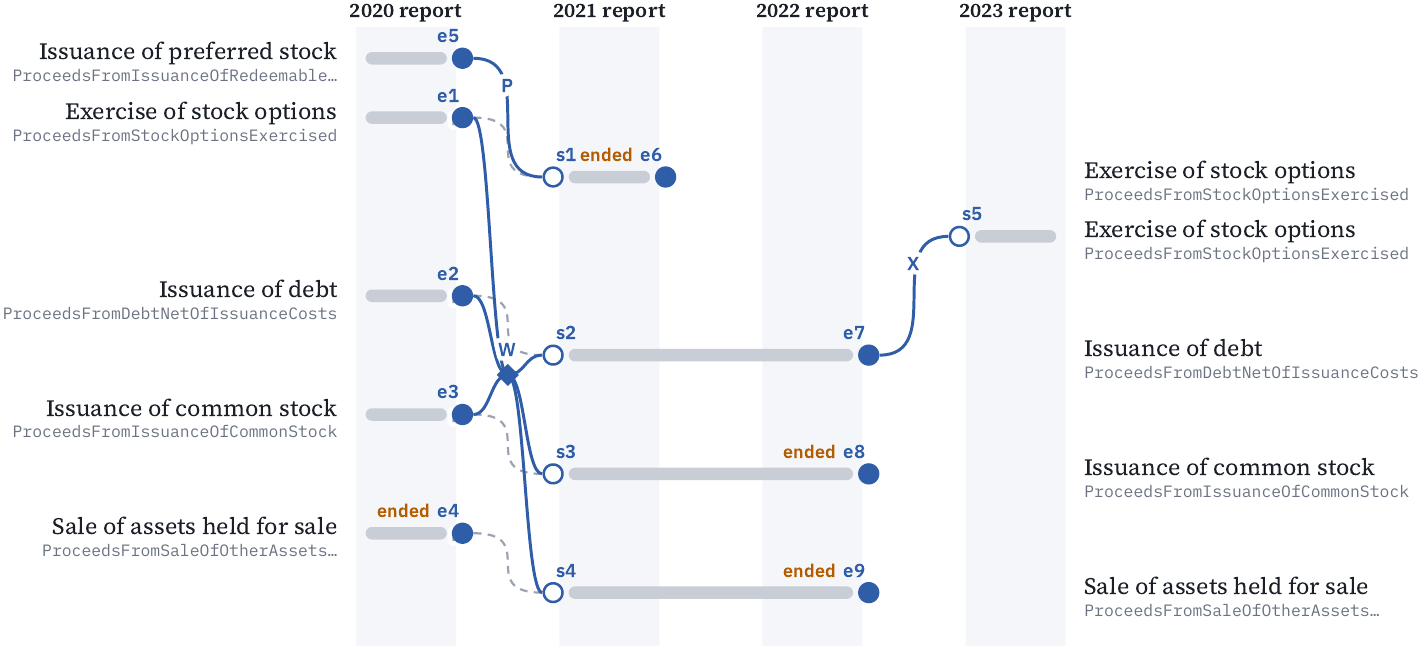}
  \caption{The solution: the swap, the preferred-stock rename and the correction solid; the same-name continuations dashed; four strands ended.}
  \label{fig:palantir-choice}
\end{figure}

\paragraph{Solution.} \Cref{fig:palantir-choice}: W, P and X. The swap
covers six endpoints at once and saves $126.46$; the same-name
continuations save $24.75$, $36.75$, $30.75$ and $19.55$, and each needs an
endpoint the swap holds. X is the correction, and nothing
competes with it. Four strands end uncovered: the old sale-of-assets
series at the first boundary, the stock-options tag's run at the
second, and the common-stock and sale tags at the last. The swap
separates each tag's run before the boundary from its run after.

\begin{figure}[htbp]
  \centering
  \includegraphics[width=\linewidth]{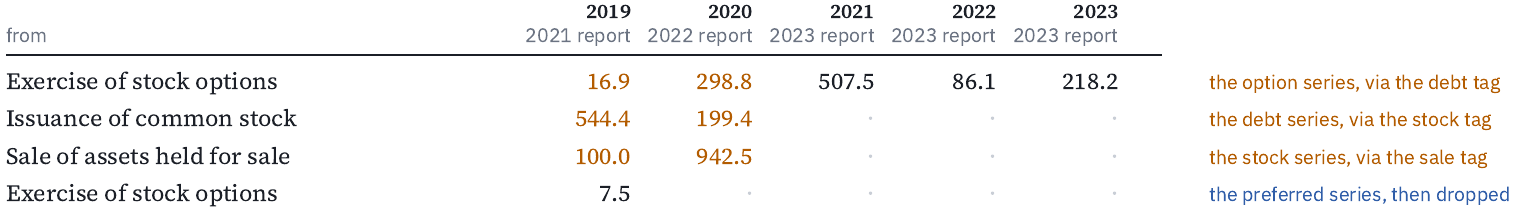}
  \caption{The rows that follow, as last reported: one series per row,
    each named by the tag it ends on (amber: years on a borrowed tag,
    flagged).}
  \label{fig:palantir-rows}
\end{figure}

\paragraph{Rows.} \Cref{fig:palantir-rows}. Each row is a chain named
by the tag it ends on, and each is one series: the stock-option proceeds
run $16.9$, $298.8$, $507.5$, $86.1$, $218.2$ across three tags and end
under their own tag; the debt and common-stock series end where their
borrowed tags were dropped. Years on a borrowed tag are flagged.
As first reported, the 2021 column shows the numbers on the tags
Palantir put them on.

\subsection{Tesla}\label{sec:tesla}

\paragraph{The facts.} \Cref{tab:tesla-facts}: the 2016 report restates
Automotive ($3{,}741.0$ to $3{,}431.6$), adds Automotive leasing ($309.4$)
and a total ($3{,}741.0$) under a new tag, and does the same to Services
and other, restated with Energy generation and storage carved out
($305.1 = 290.6 + 14.5$)~\cite{tsla-10k-2015,tsla-10k-2016}.
\Cref{sec:example} set out three explanations
of the automotive change.

\begin{table}[htbp]
  \centering
  \includegraphics[width=\linewidth]{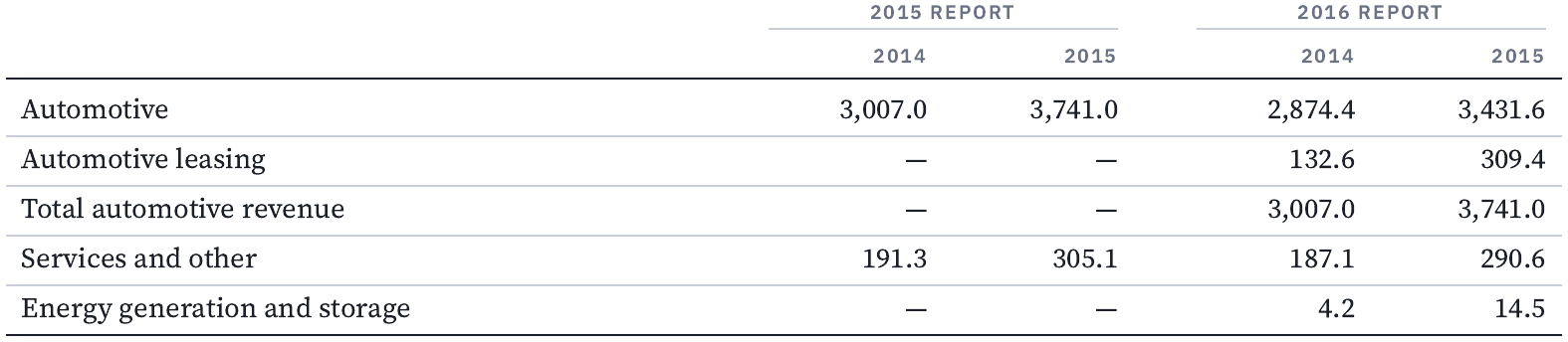}
  \caption{Five revenue lines of Tesla's income statement for the same two
    years, as two consecutive reports gave
    them~\cite{tsla-10k-2015,tsla-10k-2016}.}
  \label{tab:tesla-facts}
\end{table}

\begin{figure}[htbp]
  \centering
  \includegraphics[width=\linewidth]{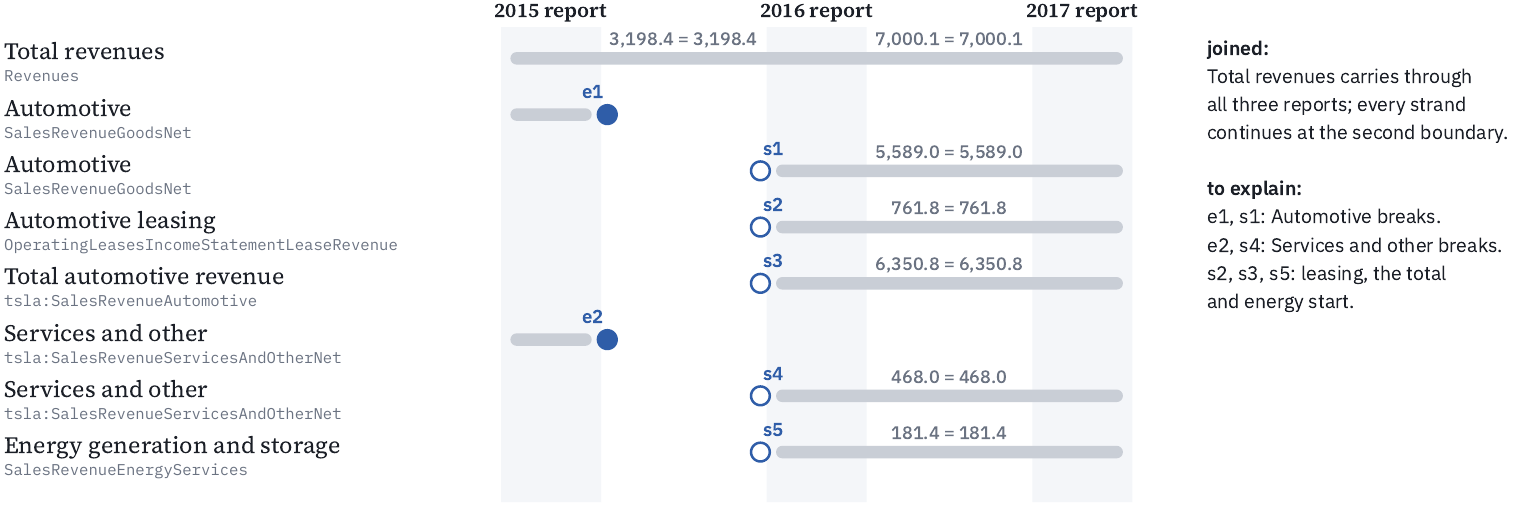}
  \caption{Tesla\textquotesingle s revenue lines as strands across the 2015 to 2017 reports.}
  \label{fig:tesla-rows-strands}
\end{figure}

\paragraph{Strands.} \Cref{fig:tesla-rows-strands}, across the 2015 to
2017 reports~\cite{tsla-10k-2017}. Automotive and Services and other each end at the 2015
report and start again at the 2016 ($e_1$, $s_1$; $e_2$, $s_4$); leasing,
the total and energy start there ($s_2$, $s_3$, $s_5$); Total revenues
runs through, and at the 2016 to 2017 boundary every strand continues.
Seven endpoints.

\begin{figure}[htbp]
  \centering
  \includegraphics[width=\linewidth]{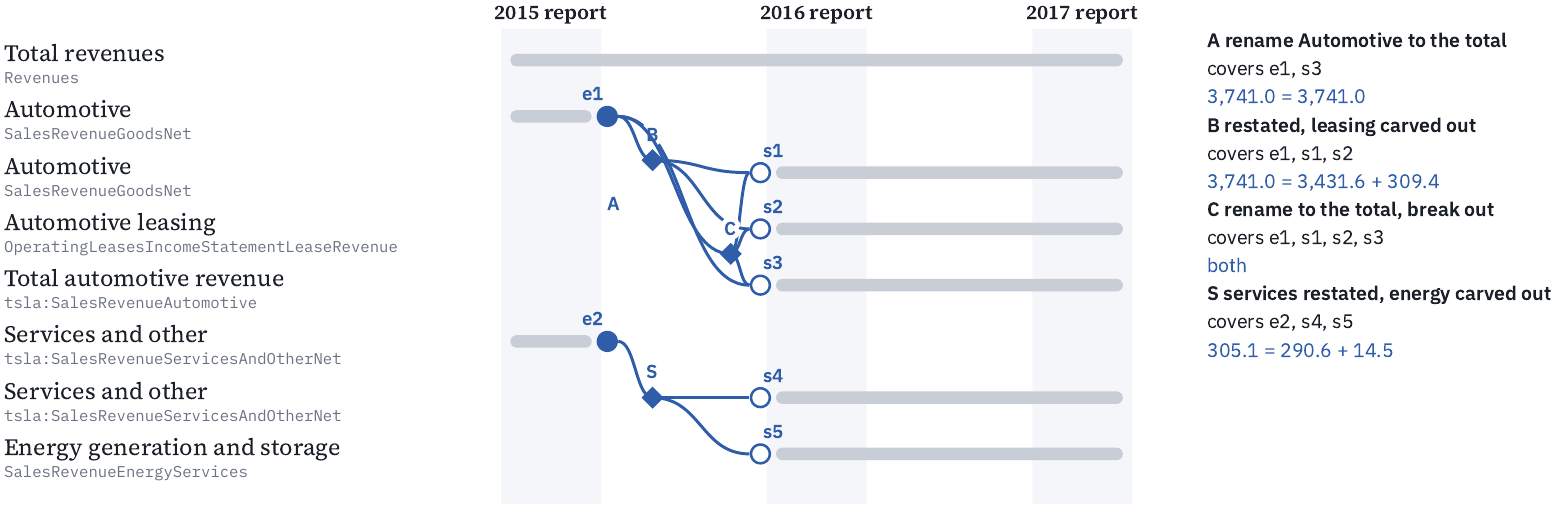}
  \caption{Every candidate at the 2015 to 2016 boundary.}
  \label{fig:tesla-rows-candidates}
\end{figure}

\paragraph{Candidates.}\label{app:tesla-choice}
\Cref{fig:tesla-rows-candidates}: explanations A, B and C of
\cref{sec:example} and the plain restated rename for the automotive
endpoints $e$, $s_1$, $s_2$, $s_3$, and S for the services endpoints.
\Cref{tab:tesla-choice} lists the automotive ones as sets with the check
each passed; all four hold $e$, so they are one group and at most one
can be chosen. \Cref{tab:tesla-prices} prices them: a candidate's $w_c$
plus the $p_n$ of what it leaves uncovered.

\begin{table}[htbp]
  \centering
  \includegraphics[width=\linewidth]{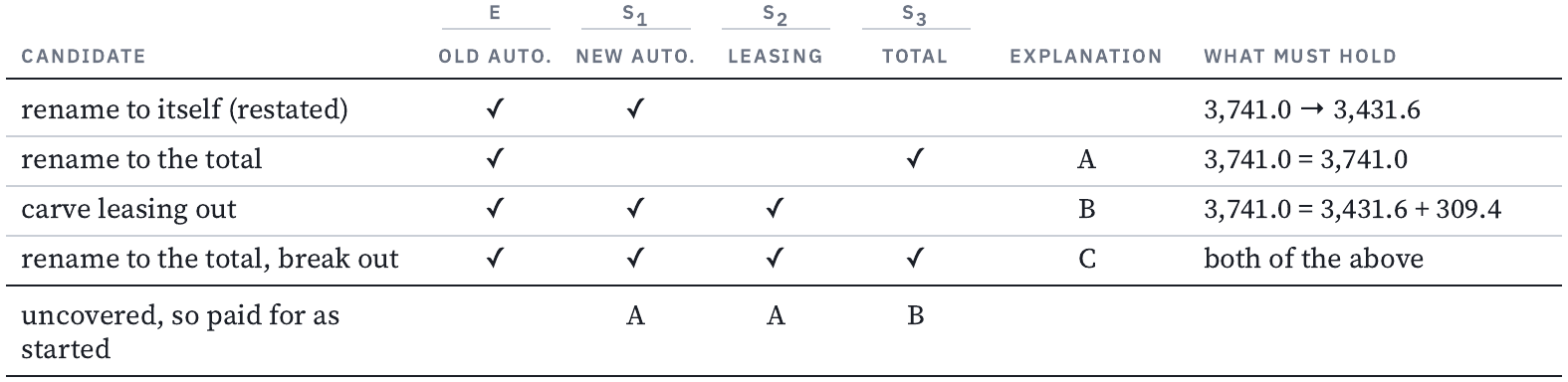}
  \caption{Tesla's 2015 to 2016 boundary as a choice, with some of its
    candidates. A tick means the candidate covers that endpoint; no two chosen
    candidates may share one. The last column is the check each candidate
    passed; the last row shows which explanation leaves each endpoint uncovered.}
  \label{tab:tesla-choice}
\end{table}
\begin{table}[htbp]
  \centering
  \includegraphics[width=\linewidth]{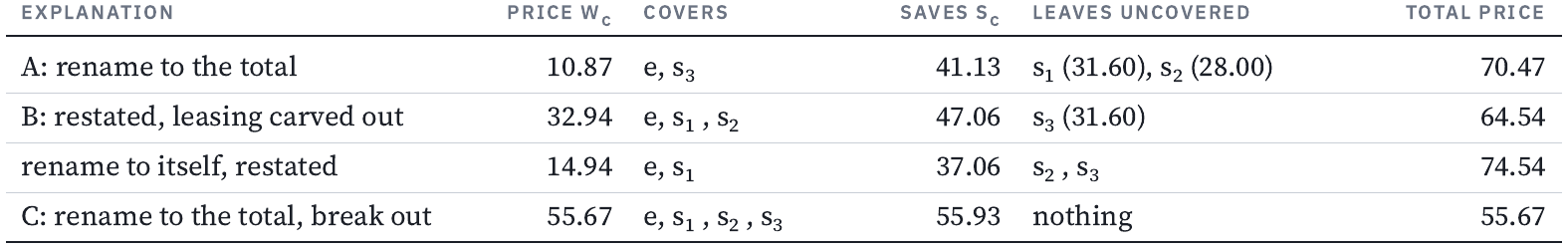}
  \caption{Tesla's 2015 to 2016 boundary, priced. The endpoints cost
    $p_e = 20.40$, $p_{s_1} = 31.60$, $p_{s_2} = 28.00$, $p_{s_3} = 31.60$
    uncovered. An explanation's total is its candidate's price plus the price of
    what it leaves uncovered; C wins by $8.87$.}
  \label{tab:tesla-prices}
\end{table}

\begin{figure}[htbp]
  \centering
  \includegraphics[width=\linewidth]{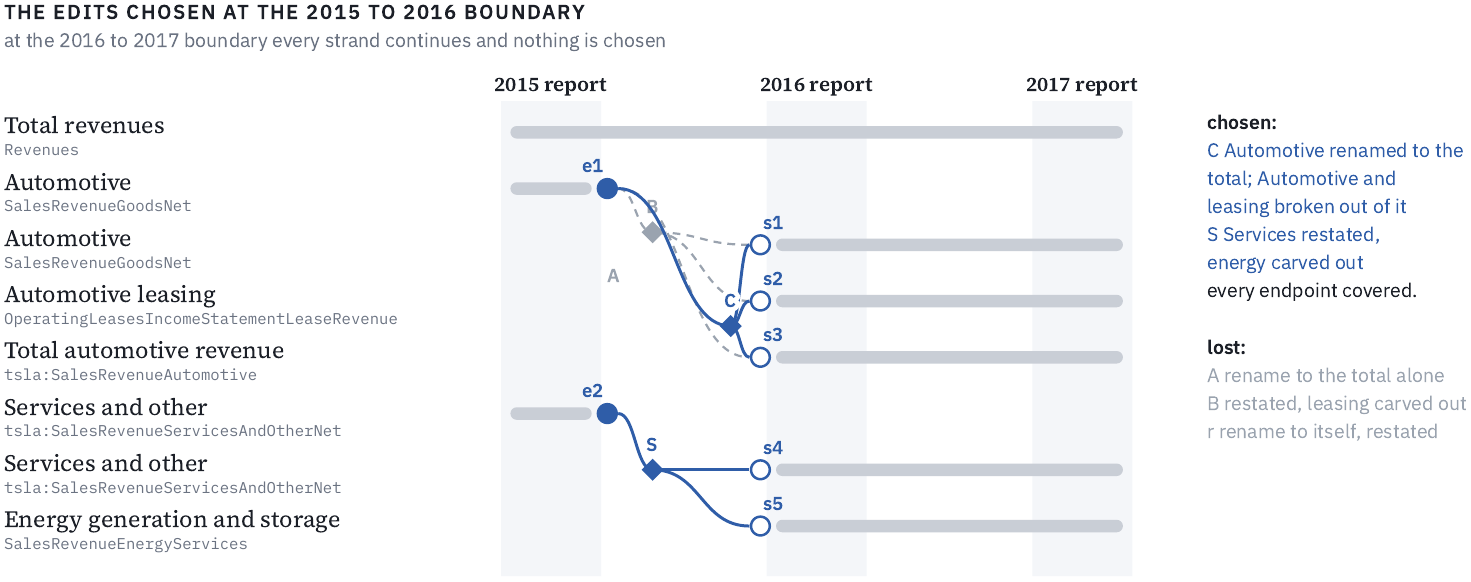}
  \caption{Tesla's revenue lines across three annual reports as strands,
    with the edits chosen at the 2015 to 2016 boundary. Chosen candidates
    are solid, the others dashed; every endpoint is covered. At the next
    boundary every strand continues.}
  \label{fig:tesla-rows-edits}
\end{figure}

\paragraph{Solution.} C and S (\cref{fig:tesla-rows-edits}). A is the
cheapest single candidate but leaves the most uncovered; C pays the most
and leaves nothing. C's price is A's rename plus the two broken-out
lines priced as lines that start, less one flat starting charge and
less a credit for the sum the total verifies over three years. A is
cheaper line by line; C is cheaper for the boundary as a whole, and
under it the total's history is like for like.

\paragraph{Rows.} \Cref{fig:tesla-rows-edits} shows the chosen edits
over all three reports; at the 2016 to 2017 boundary every
strand continues and nothing is chosen.

\begin{figure}[htbp]
  \centering
  \includegraphics[width=\linewidth]{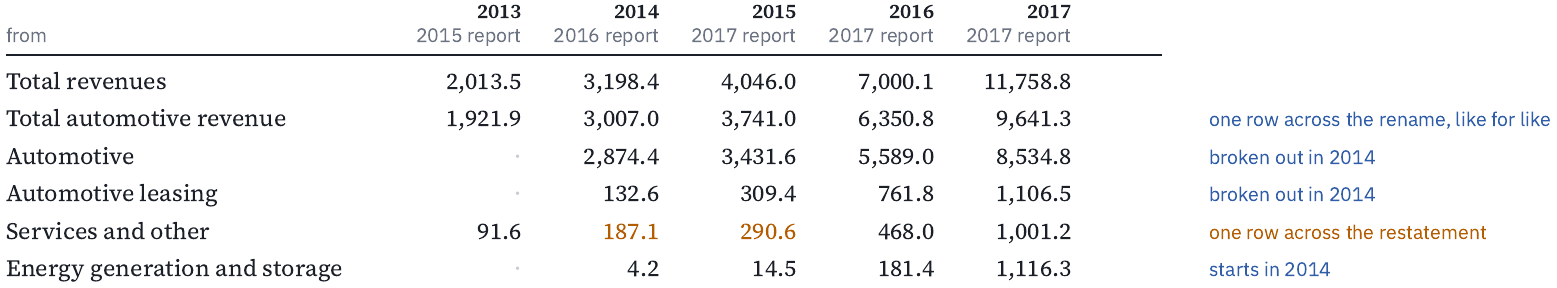}
  \caption{The rows that follow, as last reported, in millions. Each
    column's numbers are its ruling report's, the latest that shows the
    year. A rename joins two strands into one row; a number in amber was
    restated across the join; a dot is a year before a line was broken out
    or started.}
  \label{fig:tesla-rows-table}
\end{figure}

\Cref{fig:tesla-rows-table} is the table that follows, as last
reported. Each rename joins two strands into one row: the old Automotive
line and the new total are one row, like for like across all five years,
and Services and other is one row across its restatement. The three
broken-out or carved-out lines begin in 2014. Each column is ruled by the
latest report that shows it: 2013 by the 2015 report, 2014 by the 2016
report, 2015 onwards by the 2017 report; where a joined row crosses a
restatement, the ruling filing's numbers are shown and marked.

\begin{figure}[htbp]
  \centering
  \includegraphics[width=\linewidth]{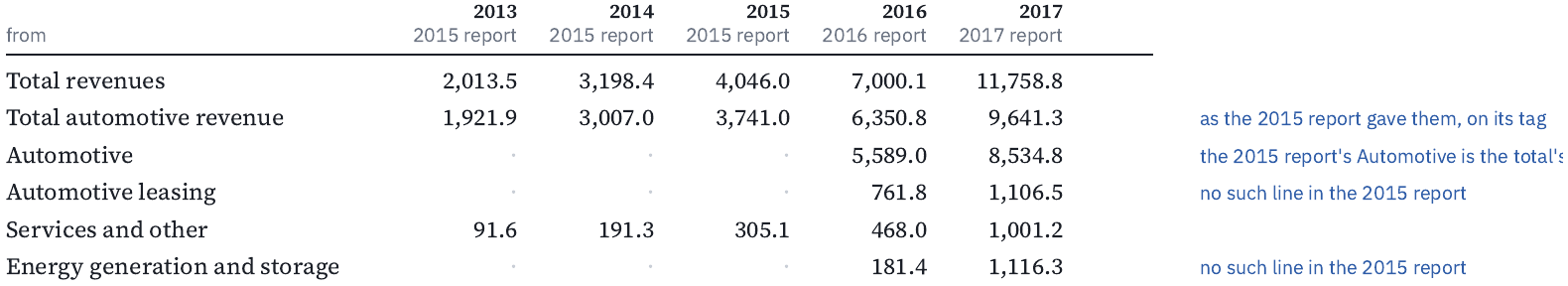}
  \caption{The same rows as first reported: each column ruled by the report
    whose year it is. A line the ruling report did not show is blank; a
    number is never taken from a later report.}
  \label{fig:tesla-rows-table-first}
\end{figure}

\Cref{fig:tesla-rows-table-first} is the same rows as first reported,
each year ruled by its own report. The total's 2014 and 2015 are the
2015 report's 3,007.0 and 3,741.0 under the tag it had then; lines the
2015 report did not show are blank there rather than filled from later,
including the new Automotive line, whose tag the 2015 report used for
what is now the total. Every column of either table adds up as its
ruling report did.

\subsection{JPMorgan}\label{sec:jpm}

\paragraph{The facts.} \Cref{tab:jpm-facts}: JPMorgan's 2023 report
replaces one line of non-interest revenue with two, and restates both
comparative years into the parts: $14{,}405 + 6{,}624 = 21{,}029$ and
$14{,}096 + 6{,}581 = 20{,}677$~\cite{jpm-10k-2022,jpm-10k-2023}.

\begin{table}[htbp]
  \centering
  \includegraphics[width=\linewidth]{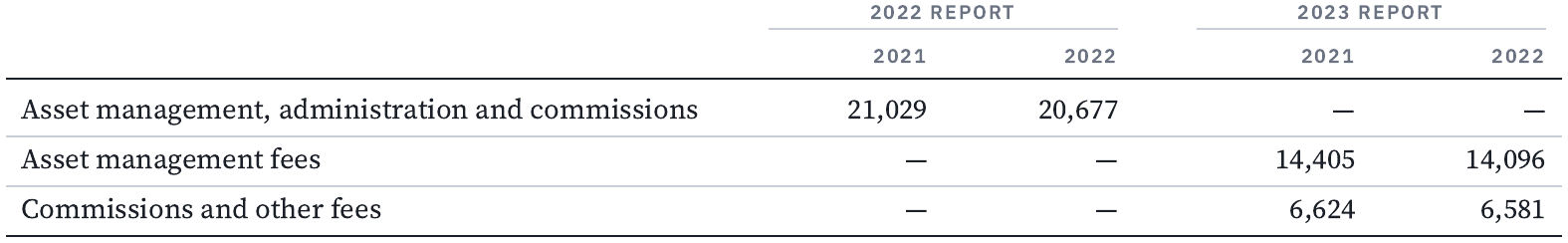}
  \caption{One line of JPMorgan's income statement becoming two, for the
    same two years in two consecutive
    reports~\cite{jpm-10k-2022,jpm-10k-2023}.}
  \label{tab:jpm-facts}
\end{table}

\begin{figure}[htbp]
  \centering
  \includegraphics[width=\linewidth]{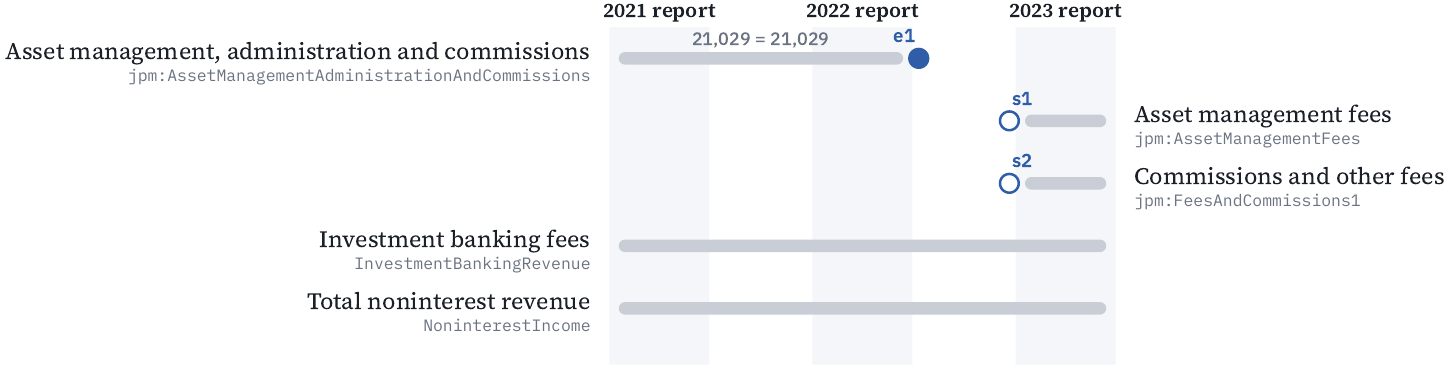}
  \caption{JPMorgan\textquotesingle s asset management line as strands across the 2021 to 2023 reports.}
  \label{fig:jpm-strands}
\end{figure}

\paragraph{Strands.} \Cref{fig:jpm-strands}. The old line joins its strand
from the 2021 to the 2022 report and ends there ($e_1$); the two parts
start in the 2023 report ($s_1$, $s_2$). Three endpoints.

\begin{figure}[htbp]
  \centering
  \includegraphics[width=\linewidth]{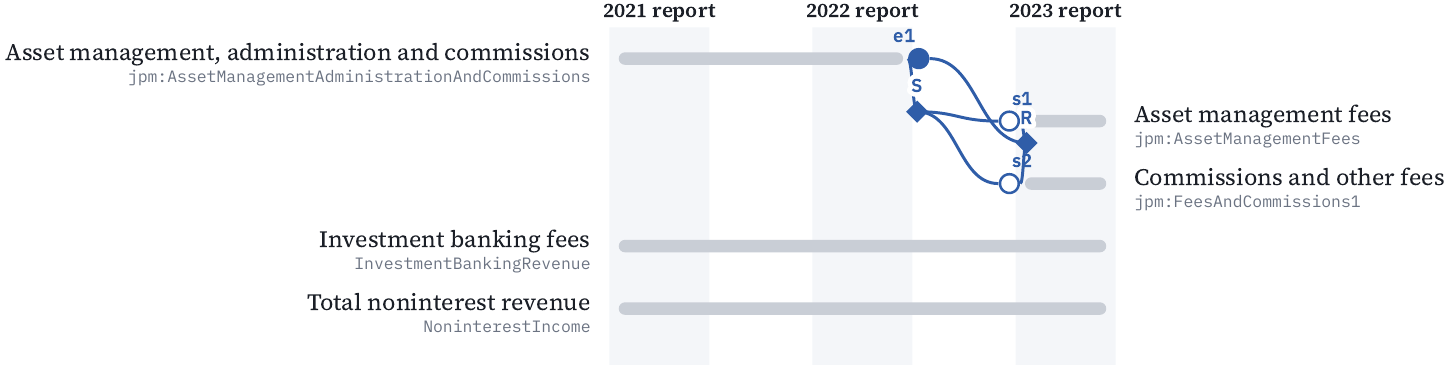}
  \caption{The two candidates at the 2022 to 2023 boundary, both covering the same three endpoints.}
  \label{fig:jpm-candidates}
\end{figure}

\paragraph{Candidates.} \Cref{fig:jpm-candidates}. Two, for the same three
endpoints:
\begin{itemize}
  \item S: a \edit{split} of the old line into the two parts, $26.15$.
  \item R: the old line renamed to asset management fees, with commissions
    carved out of it.
\end{itemize}

\begin{figure}[htbp]
  \centering
  \includegraphics[width=\linewidth]{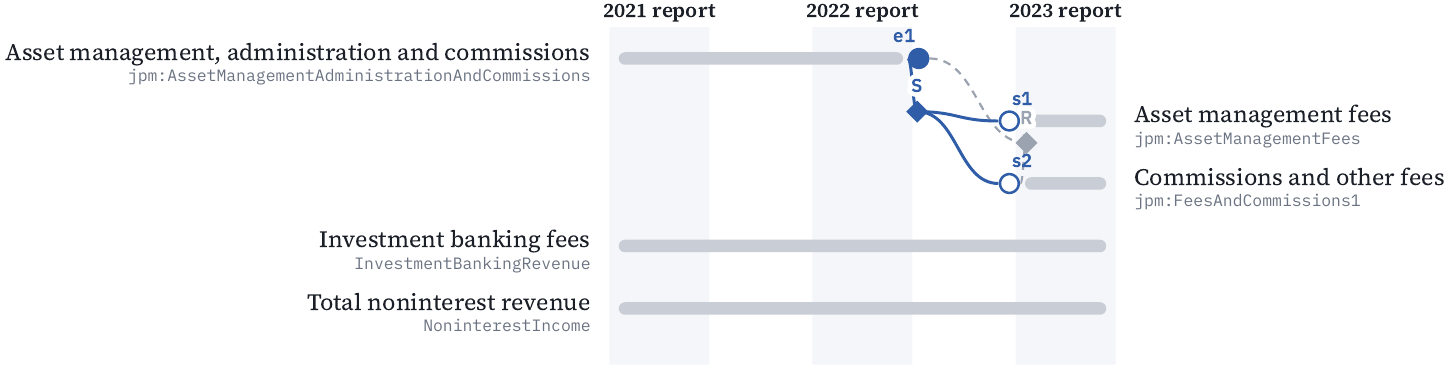}
  \caption{The solution: the split solid, the rename with a carve-out dashed.}
  \label{fig:jpm-choice}
\end{figure}

\paragraph{Solution.} \Cref{fig:jpm-choice}: S, which saves $19.85$.
Both pass the same sums, but R pays the name term of a rename to a line
with a different label, while a split's name term is judged against
both parts together; R saves $3.83$.

\begin{figure}[htbp]
  \centering
  \includegraphics[width=\linewidth]{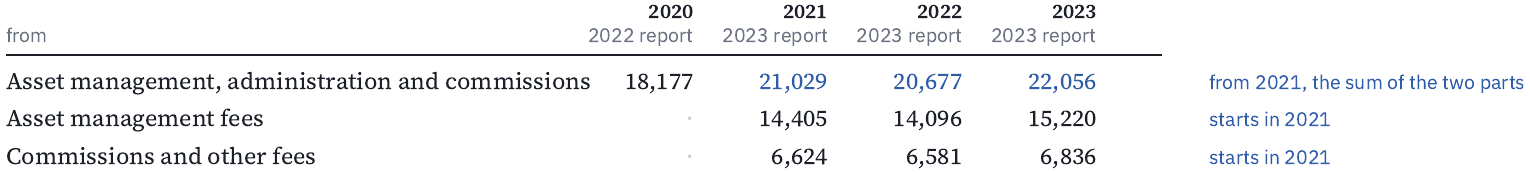}
  \caption{The rows that follow (blue: derived as the sum of the parts).}
  \label{fig:jpm-rows}
\end{figure}

\paragraph{Rows.} \Cref{fig:jpm-rows}. The old line continues as the sum of
its parts, marked as derived and cited to both; $22{,}056$ for 2023 is in
no filing. The parts are blank before 2021.

%% file: sections/d-hard-cases.tex
\section{One hard case per kind}\label{app:hard}

The four cases of Appendix~\ref{sec:solutions} were chosen because they
are instructive. The cases here were chosen from the edits the
optimization made on the two sets so that each can be checked: one edit
at one boundary, where the names give little or no help, small enough
to verify against the two filings named by accession number, and with
the boundary inside the fifteen-year window so the rows that follow can
be shown. Each case has two figures. The first shows the two
reports side by side with the edit drawn over them, in the conventions
of Appendix~\ref{app:vocabulary}: what the company printed is in serif,
the edit in blue, the shared year shaded, a restated number circled in
orange. The second is the excerpt of the finished statement, as last
reported: a restated number in orange, a derived number in blue, a dot
where a line has no number in that year. Numbers are as filed, at the
filer's own scale; a number in parentheses is negative. \edit{Ended}
and \edit{started} need no example.

\subsection{Rename, to a company's own tag: Tupperware Brands}

\begin{figure}[htbp]
  \centering
  \includegraphics[width=\linewidth]{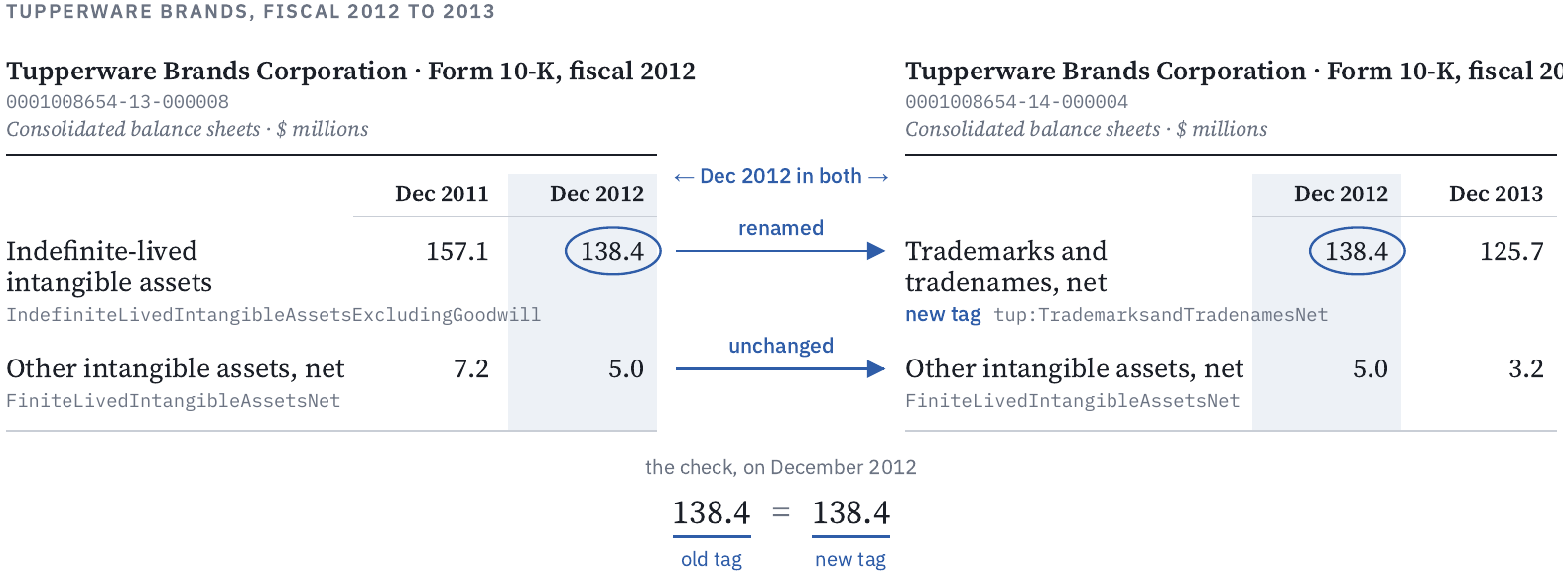}
  \caption{Tupperware Brands, balance sheet, fiscal 2012 to 2013:
    indefinite-lived intangible assets renamed to trademarks and
    tradenames, net, under a company tag. The December 2012 balance
    carries.}
  \label{fig:hard-rename-tupperware}
\end{figure}

\begin{figure}[htbp]
  \centering
  \includegraphics[width=0.9\linewidth]{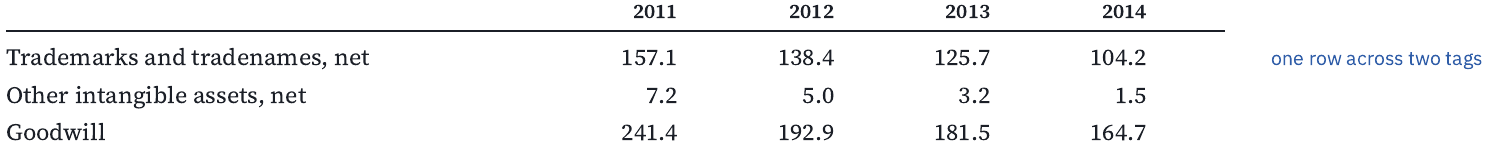}
  \caption{The rows that follow, \$ millions: one line across the
    rename.}
  \label{fig:hard-rename-tupperware-rows}
\end{figure}

Tupperware's fiscal 2013 report drops \emph{Indefinite-lived intangible
assets} and shows \emph{Trademarks and tradenames, net} under a company
tag instead (\cref{fig:hard-rename-tupperware}). Tag and label both
change and share no word; the only evidence is that the December 2012
balance, 138.4, carries, while the neighbouring lines do not move. The
rows (\cref{fig:hard-rename-tupperware-rows}) are one line across the
rename, and across a second rename in a later report.


\subsection{Rename, to a standard's tag: NOV}

\begin{figure}[htbp]
  \centering
  \includegraphics[width=\linewidth]{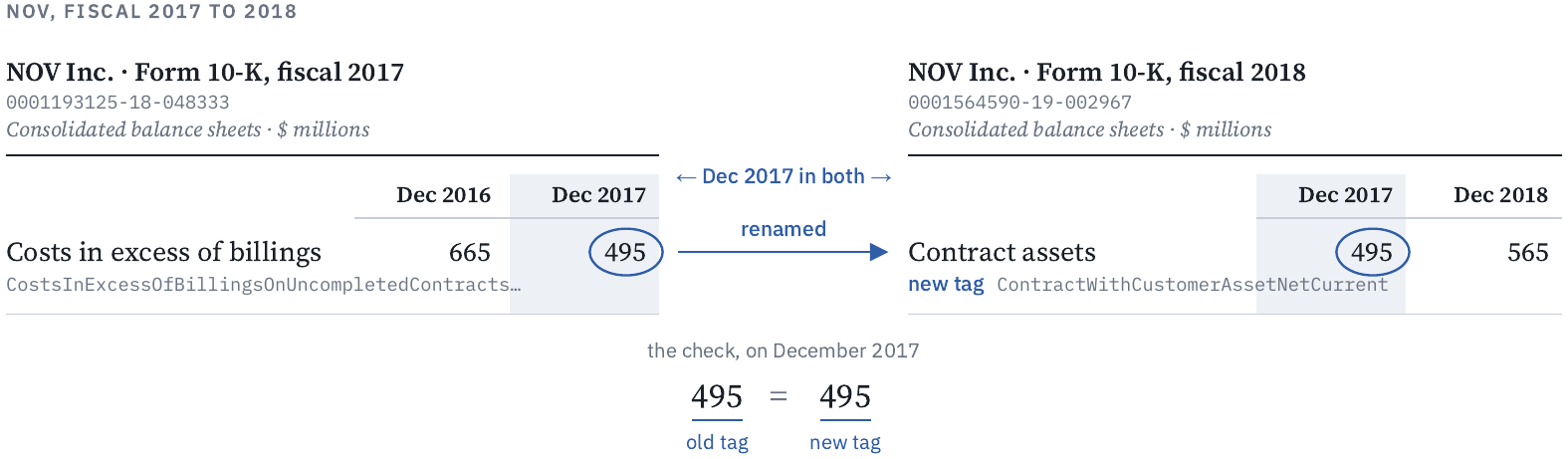}
  \caption{NOV, balance sheet, fiscal 2017 to 2018: costs in excess of
    billings renamed to contract assets under the new standard's tag.}
  \label{fig:hard-rename-nov}
\end{figure}

\begin{figure}[htbp]
  \centering
  \includegraphics[width=0.9\linewidth]{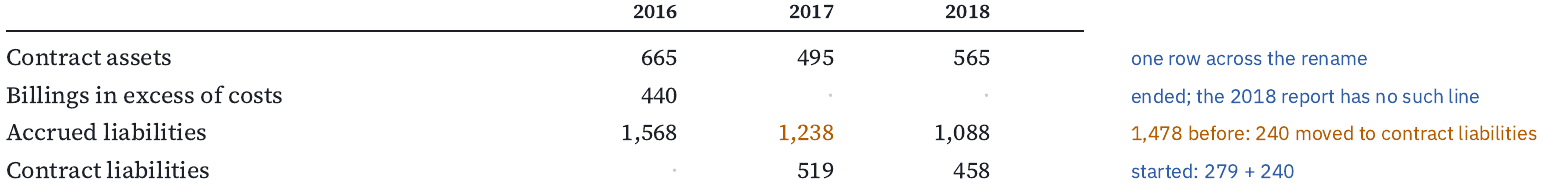}
  \caption{The rows that follow, \$ millions. The asset line is one row.
    On the liability side the method leaves a line ended, a line started
    and a bare restatement, where the filing describes one move of 240.}
  \label{fig:hard-rename-nov-rows}
\end{figure}

The new revenue standard gave contract positions new tags and new
names. NOV's \emph{Costs in excess of billings} becomes \emph{Contract
assets} in the fiscal 2018 report, with the December 2017 balance of
495 carried (\cref{fig:hard-rename-nov}). The liability side is a
limit of the method (\cref{fig:hard-rename-nov-rows}). \emph{Billings
in excess of costs} was 279 for December 2017 and \emph{Contract
liabilities} is 519, and the report says why: 240 of advance payments
and deferred revenue moved out of \emph{Accrued liabilities}, which is
restated from 1,478 to 1,238. That is one move between two continuing
lines, and no kind of edit in \cref{sec:candidates} describes it; a
carve-out's children are lines that start, not lines that continue. So
the solver leaves a line ended, a line started and a bare restatement
of accrued liabilities, three edits where the filing describes one.
Every cell holds the right number, but the record of changes is
incomplete.


\subsection{Absorb: Palantir}

\begin{figure}[htbp]
  \centering
  \includegraphics[width=\linewidth]{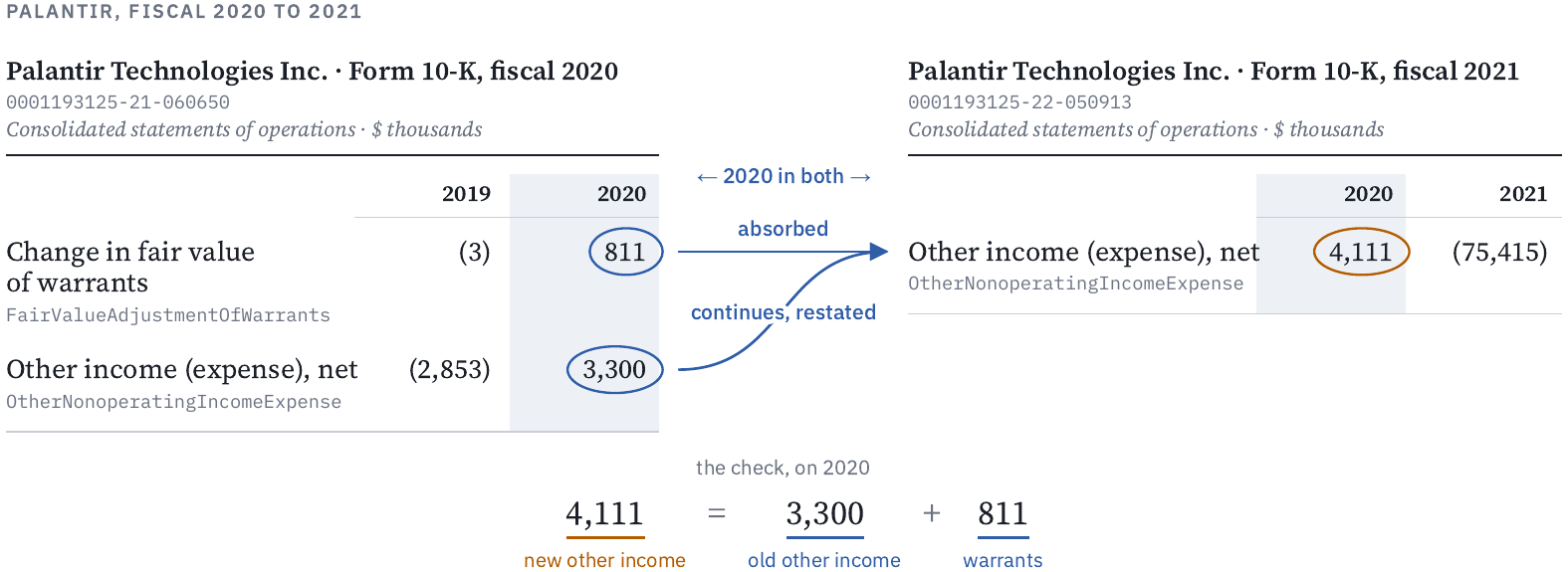}
  \caption{Palantir, income statement, fiscal 2020 to 2021: the warrant
    line absorbed into other income (expense), net. Amounts as they
    enter income.}
  \label{fig:hard-absorb-palantir}
\end{figure}

\begin{figure}[htbp]
  \centering
  \includegraphics[width=0.9\linewidth]{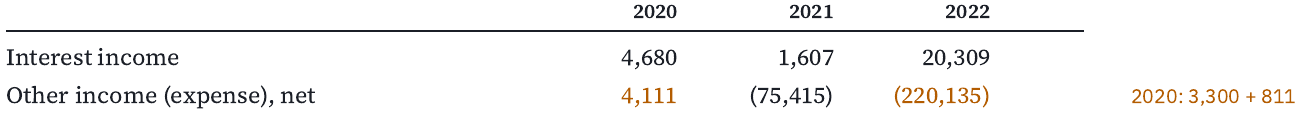}
  \caption{The rows that follow, \$ thousands, as last reported.}
  \label{fig:hard-absorb-palantir-rows}
\end{figure}

Palantir's fiscal 2021 report has no \emph{Change in fair value of
warrants} line, and \emph{Other income (expense), net} for 2020 reads
4,111 where the year before it read 3,300
(\cref{fig:hard-absorb-palantir}). The difference is the warrant line,
811, to the dollar, and 2019 moves by the same line's 3. Nothing in
either name relates the two, and 0.8 million is small in an income
statement measured in billions. Every other line is unchanged. The rows
(\cref{fig:hard-absorb-palantir-rows}) carry the restated 4,111, marked;
the warrant line has no row.


\subsection{Carve-out: Enterprise Financial Services}

\begin{figure}[htbp]
  \centering
  \includegraphics[width=\linewidth]{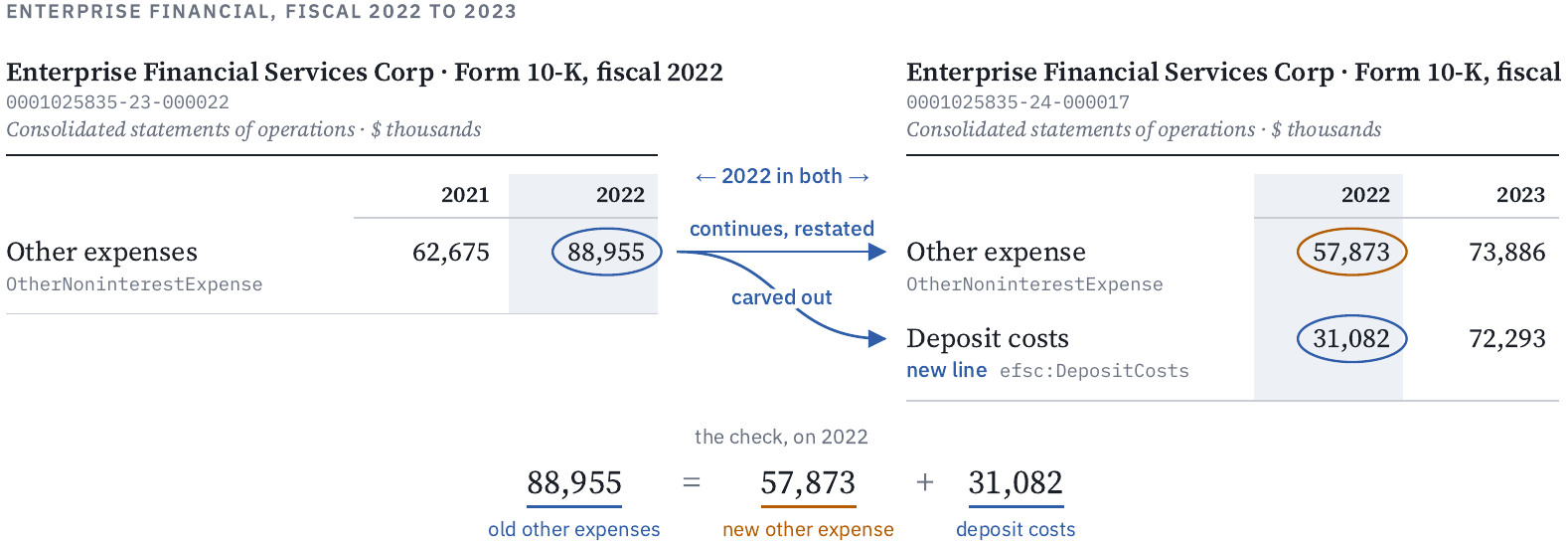}
  \caption{Enterprise Financial Services, income statement, fiscal 2022
    to 2023: deposit costs carved out of other expenses.}
  \label{fig:hard-carve-efsc}
\end{figure}

\begin{figure}[htbp]
  \centering
  \includegraphics[width=0.9\linewidth]{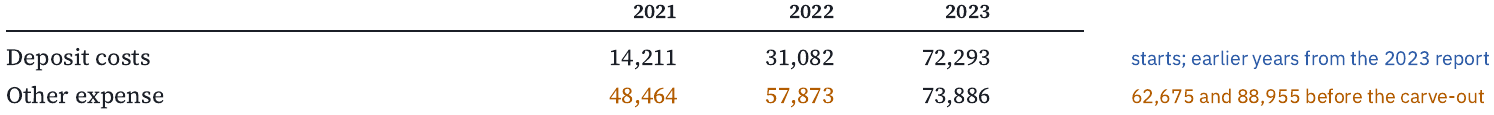}
  \caption{The rows that follow, \$ thousands, as last reported.}
  \label{fig:hard-carve-efsc-rows}
\end{figure}

Enterprise Financial's fiscal 2023 report adds a line, \emph{Deposit
costs}, and restates \emph{Other expenses} down by exactly it in both
comparative years: 88,955 becomes 57,873 plus 31,082 for 2022, and
62,675 becomes 48,464 plus 14,211 for 2021
(\cref{fig:hard-carve-efsc}). The carved-out line shares no word with
the line it came from, and the line it came from is called
\emph{Other}. The rows (\cref{fig:hard-carve-efsc-rows}) show other
expense restated for the two years the carve-out reached and deposit
costs starting with the earliest year the 2023 report shows.


\subsection{Merge: Costco}

\begin{figure}[htbp]
  \centering
  \includegraphics[width=\linewidth]{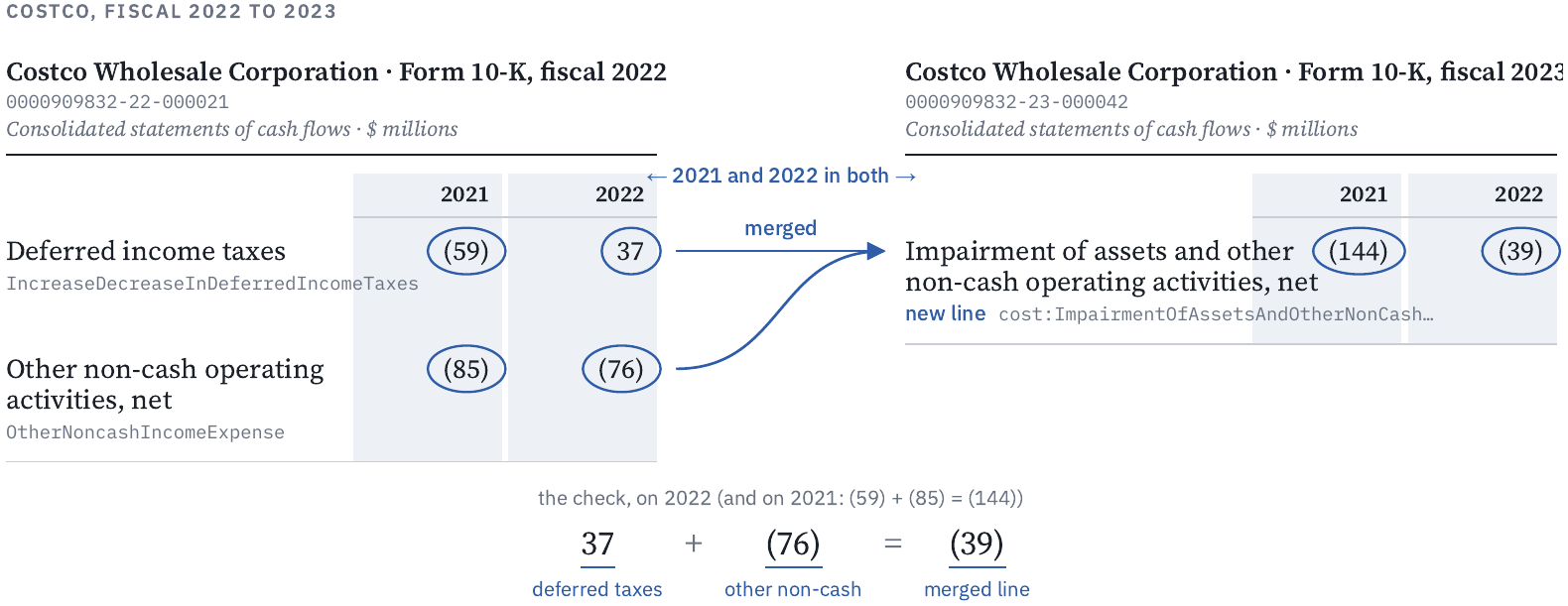}
  \caption{Costco, cash flow statement, fiscal 2022 to 2023: two
    operating adjustments merged into one line, checked on both shared
    years.}
  \label{fig:hard-merge-costco}
\end{figure}

\begin{figure}[htbp]
  \centering
  \includegraphics[width=0.9\linewidth]{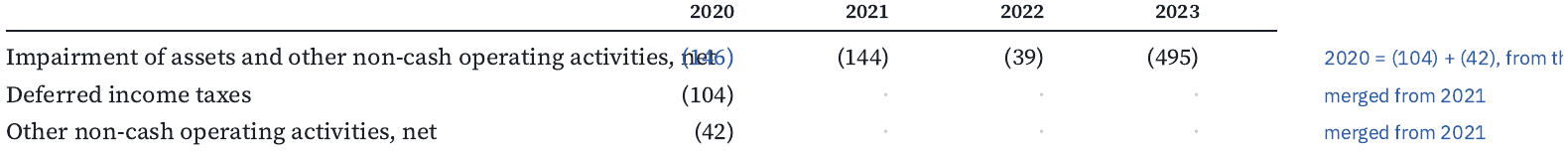}
  \caption{The rows that follow, \$ millions.}
  \label{fig:hard-merge-costco-rows}
\end{figure}

Costco's fiscal 2023 report shows one line, \emph{Impairment of assets
and other non-cash operating activities, net}, where the year before
showed \emph{Deferred income taxes} and \emph{Other non-cash operating
activities, net}. The sum holds on both years the two reports share:
$(59) + (85) = (144)$ for 2021 and $37 + (76) = (39)$ for 2022
(\cref{fig:hard-merge-costco}). The new tag is the company's own and
shares one word with one of the old ones. In the rows
(\cref{fig:hard-merge-costco-rows}) the merged line runs back through
2020 as the cited sum of its parts, and the parts are blank after the
merge. A merge on a balance sheet can be checked on one instant only,
since consecutive reports share one; here two years foot, which is what
the derived 2020 figure rests on.


\subsection{Split: Macatawa Bank}

\begin{figure}[htbp]
  \centering
  \includegraphics[width=\linewidth]{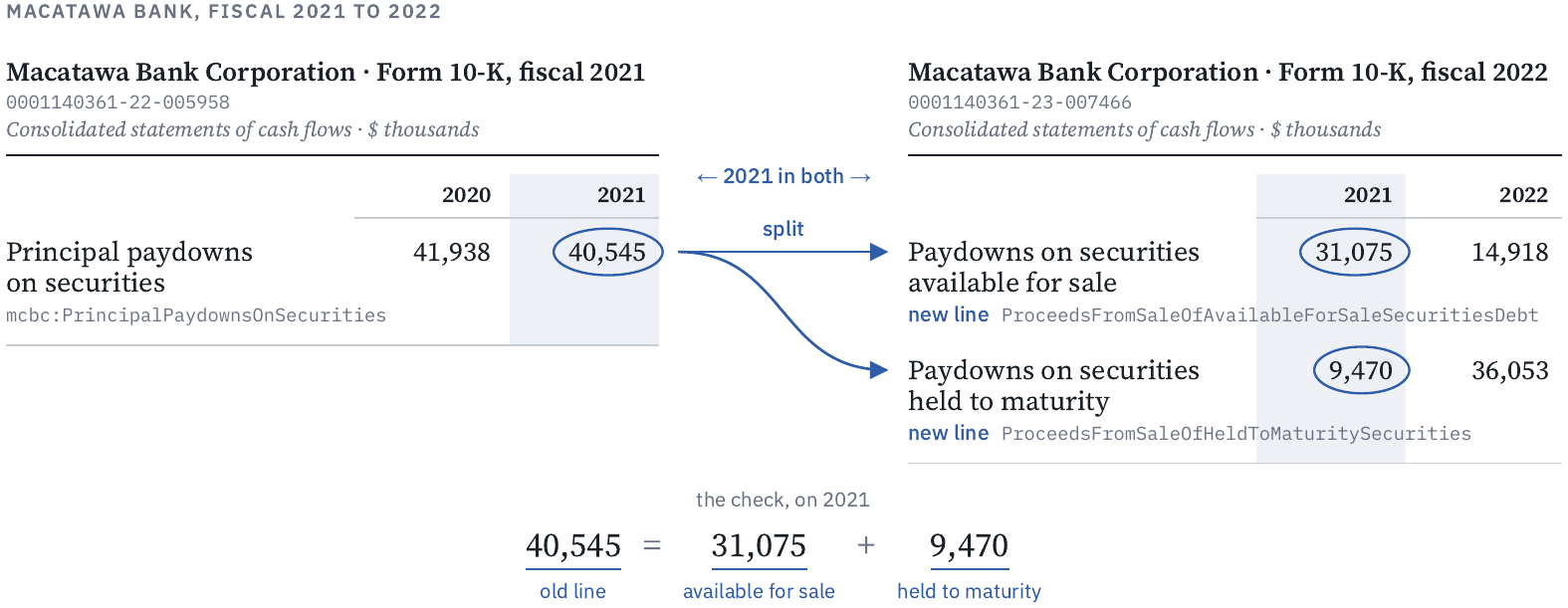}
  \caption{Macatawa Bank, cash flow statement, fiscal 2021 to 2022: one
    paydowns line split into two.}
  \label{fig:hard-split-macatawa}
\end{figure}

\begin{figure}[htbp]
  \centering
  \includegraphics[width=0.9\linewidth]{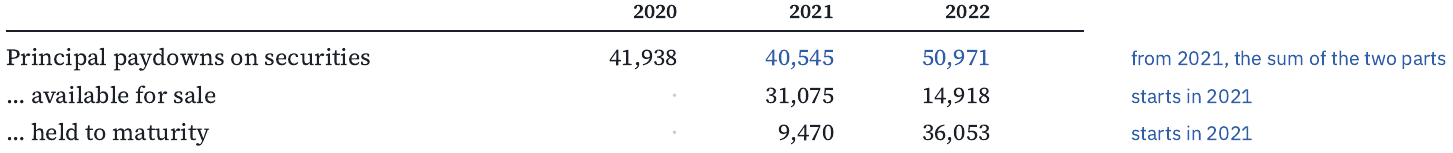}
  \caption{The rows that follow, \$ thousands.}
  \label{fig:hard-split-macatawa-rows}
\end{figure}

Macatawa's fiscal 2022 report replaces \emph{Principal paydowns on
securities} with two lines, one for securities available for sale and
one for securities held to maturity; for 2021 the two are 31,075 and
9,470, and the old line was 40,545 (\cref{fig:hard-split-macatawa}). The
old tag is the company's own and the new ones are standard tags for
proceeds from sales, which is not what the lines say. The rows
(\cref{fig:hard-split-macatawa-rows}) keep the old line as the cited
sum of its parts from 2021 on.


\subsection{Fold: Berkshire Hathaway}

\begin{figure}[htbp]
  \centering
  \includegraphics[width=\linewidth]{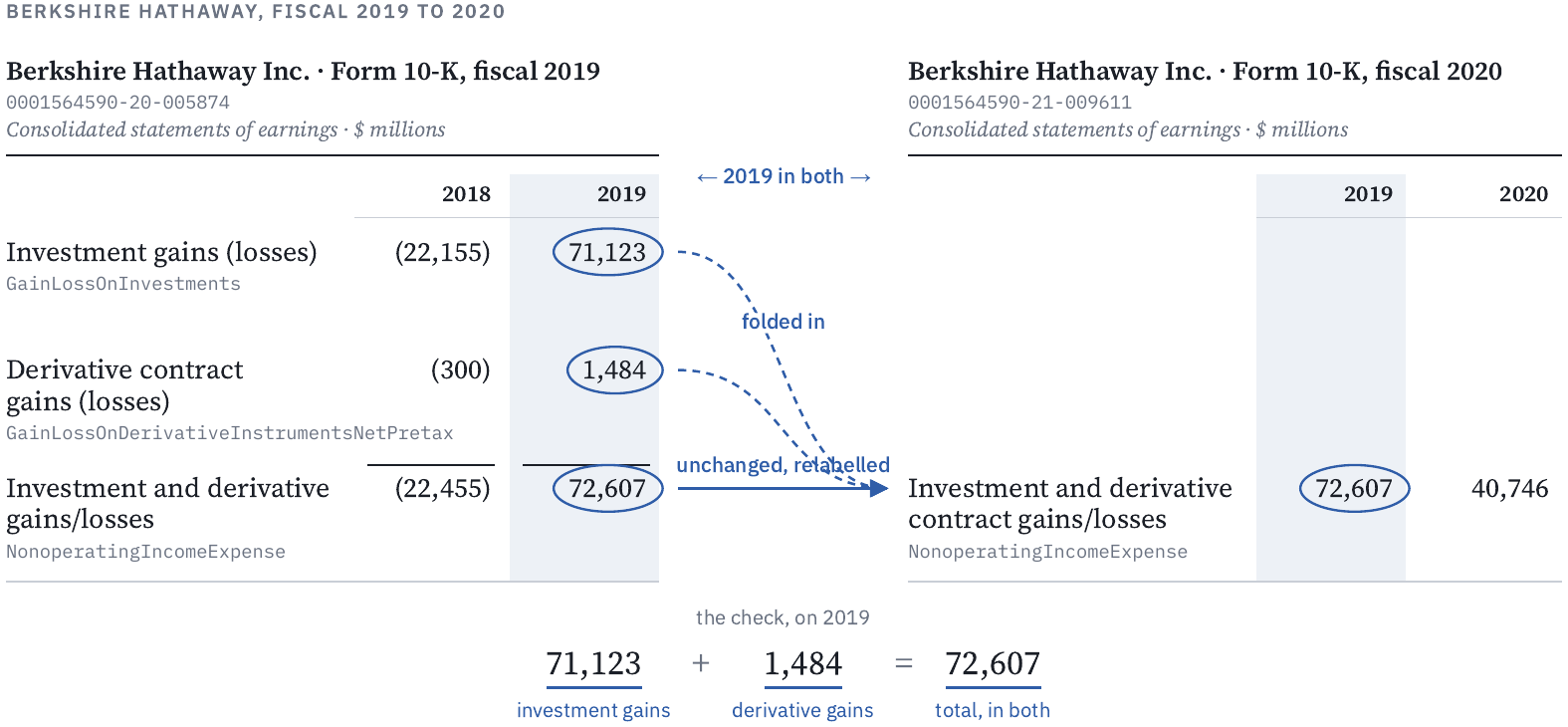}
  \caption{Berkshire Hathaway, income statement, fiscal 2019 to 2020:
    two gain lines folded into their total.}
  \label{fig:hard-fold-berkshire}
\end{figure}

\begin{figure}[htbp]
  \centering
  \includegraphics[width=0.9\linewidth]{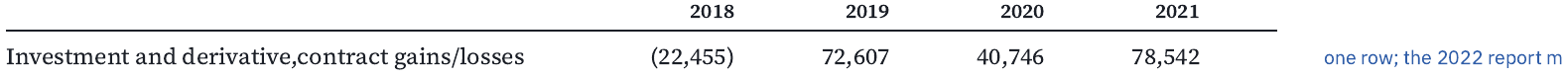}
  \caption{The rows that follow, \$ millions, as last reported.}
  \label{fig:hard-fold-berkshire-rows}
\end{figure}

Berkshire's fiscal 2019 report shows \emph{Investment gains (losses)}
and \emph{Derivative contract gains (losses)} as the two parts of
\emph{Investment and derivative gains/losses}. The fiscal 2020 report
shows only the total, relabelled, with 2019 still 72,607
(\cref{fig:hard-fold-berkshire}). Two lines vanish and the surviving
line's label changes in the same report, which reads as a rename and two
lines ended until the declared arithmetic is checked. As last reported
(\cref{fig:hard-fold-berkshire-rows}) the parts have no row and the
total is one row, through the fold and through a change of tag in the
2022 report.


\subsection{Break out: JPMorgan Chase}

\begin{figure}[htbp]
  \centering
  \includegraphics[width=\linewidth]{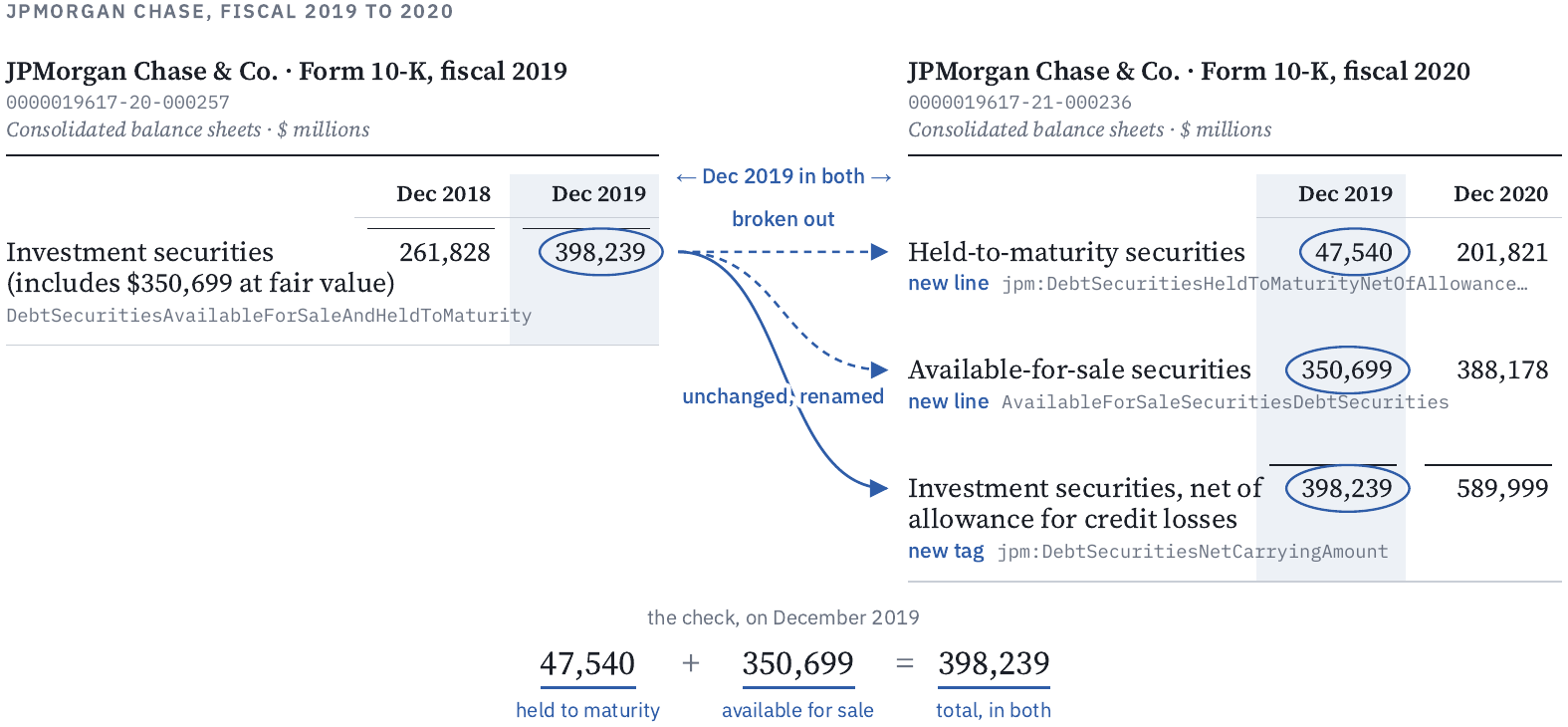}
  \caption{JPMorgan Chase, balance sheet, fiscal 2019 to 2020:
    investment securities broken out into held-to-maturity and
    available-for-sale, under a total that changes tag and does not
    move.}
  \label{fig:hard-breakout-jpm}
\end{figure}

\begin{figure}[htbp]
  \centering
  \includegraphics[width=0.9\linewidth]{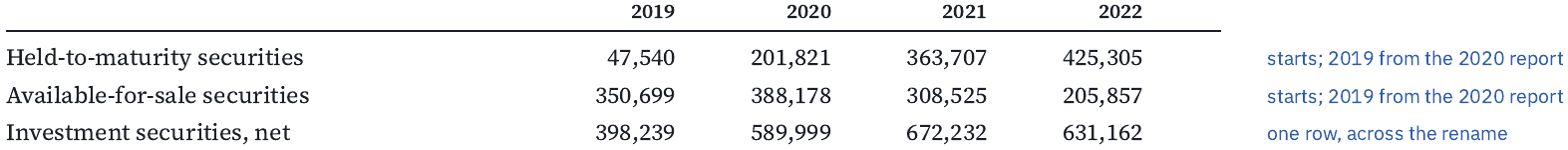}
  \caption{The rows that follow, \$ millions.}
  \label{fig:hard-breakout-jpm-rows}
\end{figure}

JPMorgan's fiscal 2019 report shows one line of investment securities,
398,239 at December 2019, with the available-for-sale portion in a
parenthetical. The fiscal 2020 report shows the same total under a new
tag and, beneath it, held-to-maturity and available-for-sale securities
as separate lines, 47,540 and 350,699 for December 2019
(\cref{fig:hard-breakout-jpm}). The parts are material in the shared
year and add up to a total that did not move; that is what a break out
is, and what tells it apart from new lines that are zero in the shared
years, which are lines that started (\cref{sec:candidates}). The rows
(\cref{fig:hard-breakout-jpm-rows}) carry the parts from 2019 and the
total as one row across its change of tag.


\subsection{Equal copies: KalVista Pharmaceuticals}

\begin{figure}[htbp]
  \centering
  \includegraphics[width=\linewidth]{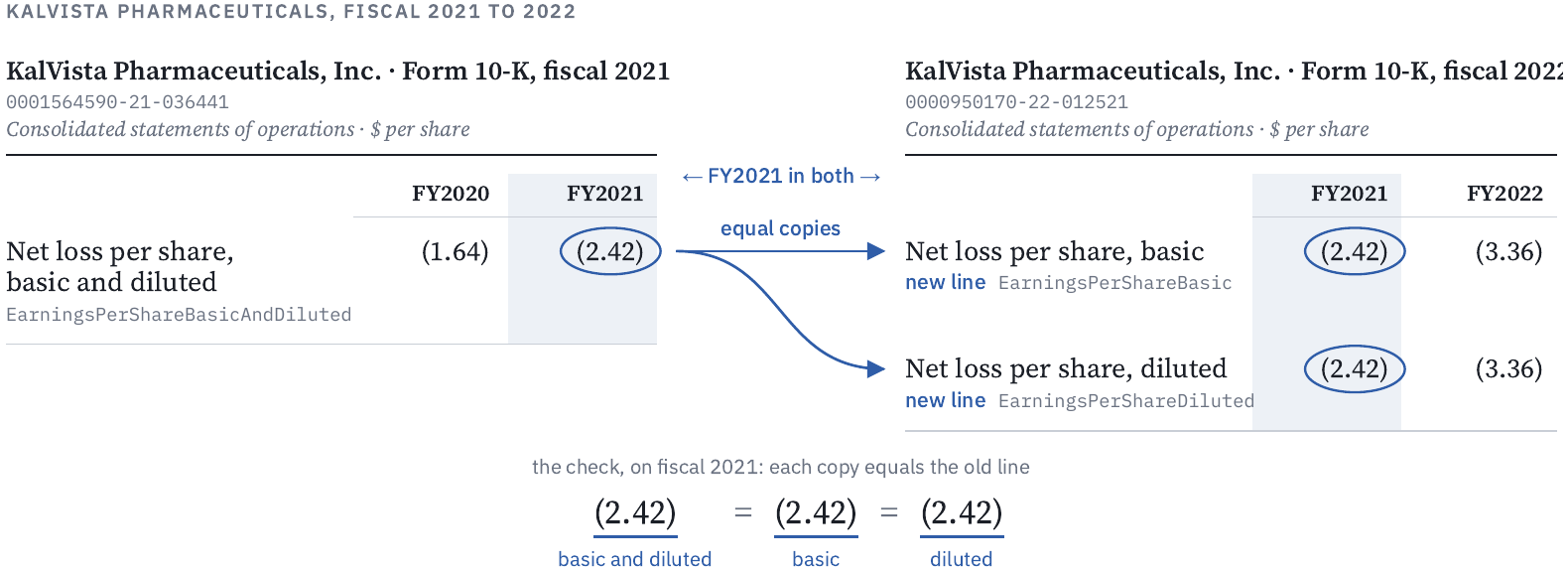}
  \caption{KalVista Pharmaceuticals, income statement, fiscal 2021 to
    2022: one per-share line becomes two equal copies.}
  \label{fig:hard-equal-kalvista}
\end{figure}

\begin{figure}[htbp]
  \centering
  \includegraphics[width=0.9\linewidth]{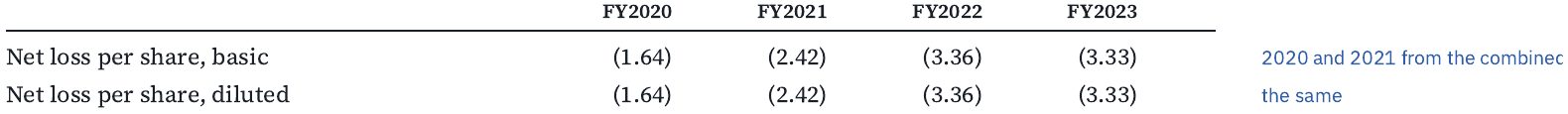}
  \caption{The rows that follow, \$ per share.}
  \label{fig:hard-equal-kalvista-rows}
\end{figure}

KalVista reported one line, \emph{Net loss per share, basic and
diluted}, until its fiscal 2022 report, which reports basic and diluted
on two lines; for fiscal 2021 each reads (2.42), as the one line did
(\cref{fig:hard-equal-kalvista}). A split would test whether the two
add up to the old line, and they do not; each equals it. The rows
(\cref{fig:hard-equal-kalvista-rows}) are two lines, each running back
through the combined line's years.


\subsection{Swap: FormFactor}

\begin{figure}[htbp]
  \centering
  \includegraphics[width=\linewidth]{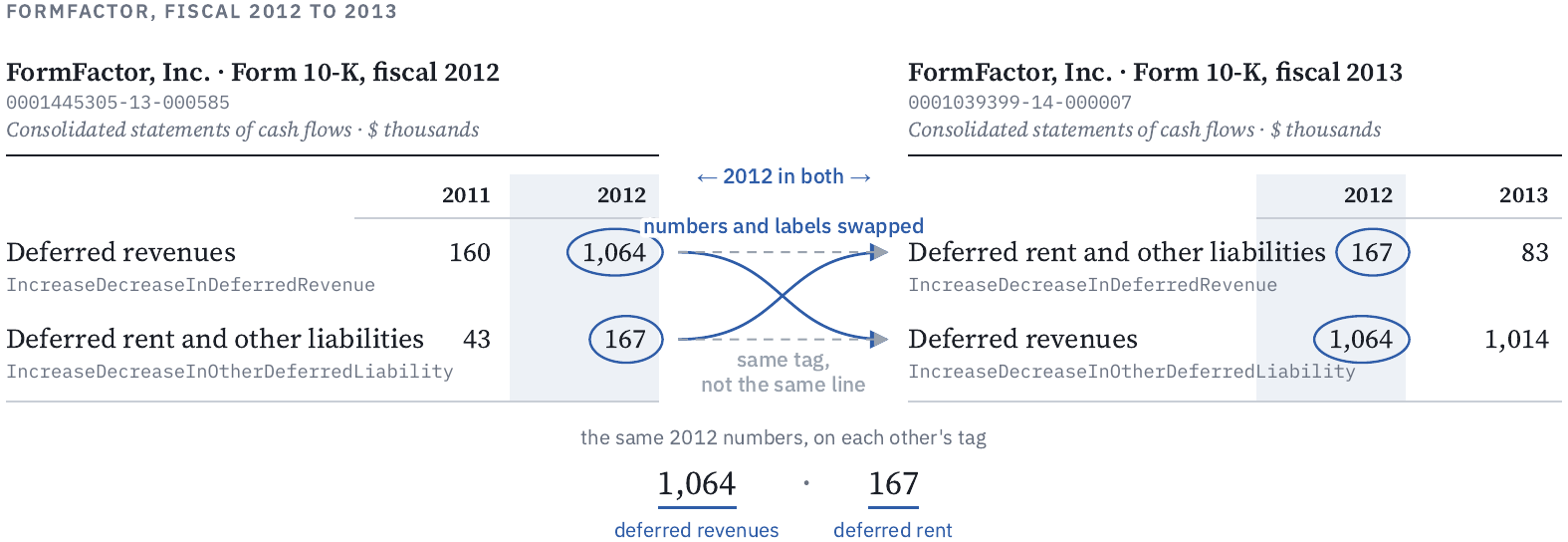}
  \caption{FormFactor, cash flow statement, fiscal 2012 to 2013: two
    tags exchange their labels and their numbers.}
  \label{fig:hard-swap-formfactor}
\end{figure}

\begin{figure}[htbp]
  \centering
  \includegraphics[width=0.9\linewidth]{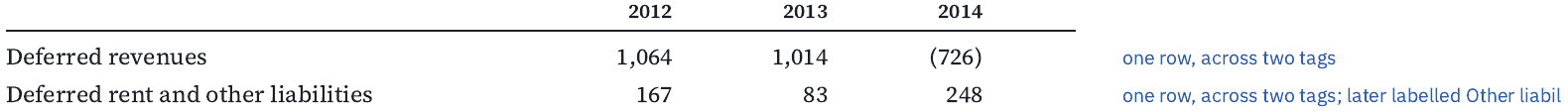}
  \caption{The rows that follow, \$ thousands: each line follows its
    numbers, not its tag.}
  \label{fig:hard-swap-formfactor-rows}
\end{figure}

In FormFactor's fiscal 2012 report the deferred-revenue tag carries
\emph{Deferred revenues}, 1,064 for 2012, and the other-deferred-
liability tag carries \emph{Deferred rent and other liabilities}, 167.
In the fiscal 2013 report the two tags have exchanged both labels and
numbers (\cref{fig:hard-swap-formfactor}). Following the tag gives each
line the other's history; following the label is right, but only a
swap says why every number under each tag disagrees with the year
before. The rows (\cref{fig:hard-swap-formfactor-rows}) follow the
numbers.


\subsection{Carve-out across the balance sheet: DXP Enterprises}

\begin{figure}[htbp]
  \centering
  \includegraphics[width=\linewidth]{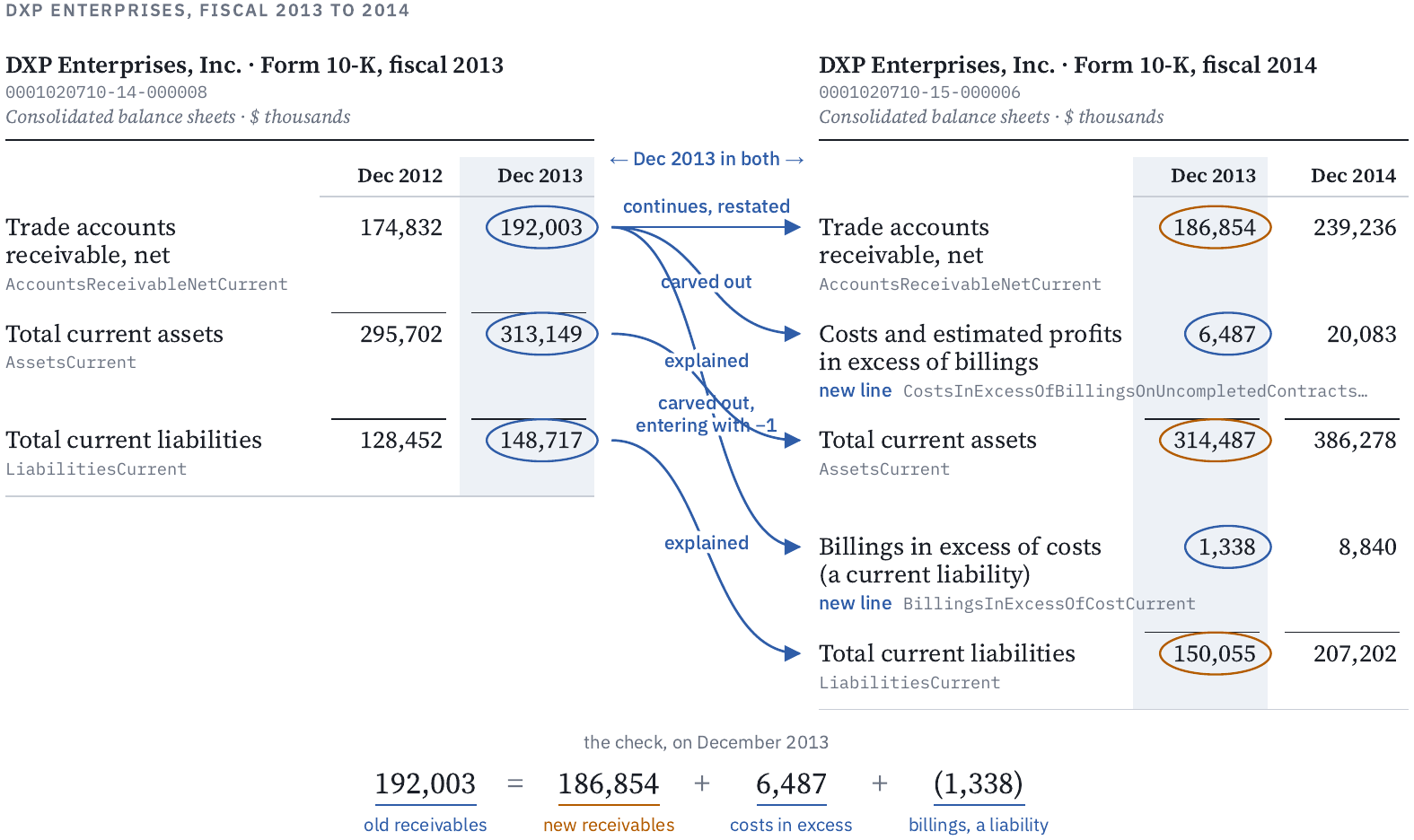}
  \caption{DXP Enterprises, balance sheet, fiscal 2013 to 2014: two lines
    carved out of receivables, one on each side of the balance sheet, with
    both section totals explained.}
  \label{fig:hard-transfer-dxp}
\end{figure}

\begin{figure}[htbp]
  \centering
  \includegraphics[width=0.9\linewidth]{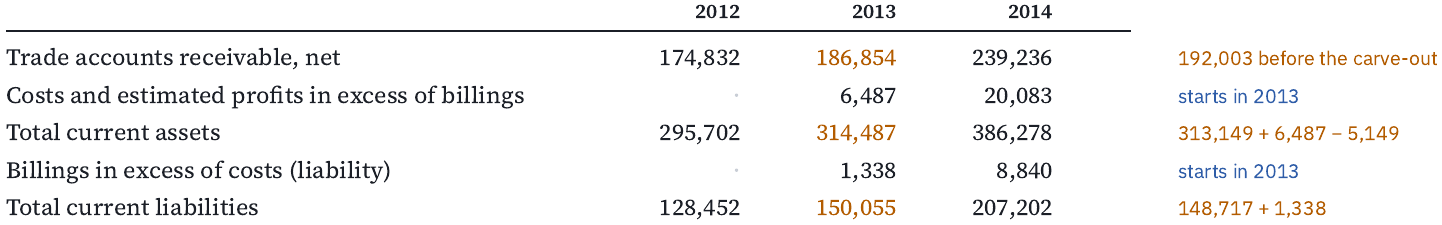}
  \caption{The rows that follow, \$ thousands, as last reported.}
  \label{fig:hard-transfer-dxp-rows}
\end{figure}

DXP's fiscal 2014 report changes its contract positions from a net to a
gross presentation. \emph{Trade accounts receivable, net} for December
2013 was 192,003; the new report shows receivables of 186,854, a new
asset line, \emph{Costs and estimated profits in excess of billings},
6,487, and a new liability, \emph{Billings in excess of costs}, 1,338,
and it restates both section totals up by 1,338
(\cref{fig:hard-transfer-dxp}). The edit is a carve-out of both new
lines from receivables: 192,003 is 186,854 plus 6,487 less 1,338, the
liability entering with $-1$ because its balance attribute is the
opposite of the receivable's (\cref{sec:candidates}), and the three
totals above them, current assets, current liabilities and total assets,
are explained by what the edit moved into each. No name relates any of
the three lines to the others. The rows (\cref{fig:hard-transfer-dxp-rows})
keep receivables as one row, restated where the carve-out reached, and
show both totals' restatements marked.


%% file: sections/e-measurements.tex
\section{Measurements in detail}\label{app:measurements}

\subsection{How often the declared arithmetic holds}\label{app:footing}

The arithmetic says which lines add up to which totals, and with what
sign, stated once per tag in a separate file of the filing
(\cref{fig:xbrl-tag}). A filing must declare it wherever a
statement shows lines and their total~\cite{edgar-xbrl-guide}, but a
total that does not add up has no effect on whether EDGAR accepts the
filing~\cite{edgar-xbrl-guide,efm-vol2}, so how often the equations
hold is a measurement. We made it over the 200 companies drawn at
random for \cref{sec:results} and their 2,826 annual reports. Over 99\% of
declared totals state a number. A check is one declared total
in one period, asking whether the parts the same filing states add up
to it; there are 104,408. Read as declared, with each breakdown row
counted as a part, 87.7\% hold exactly; folding the breakdown rows into
their line and allowing the filer's rounding brings that to 89.8\%.
Most of the rest have a declared part with no number in that period;
among totals whose every part carries a number, 98.7\% add up
(\cref{tab:footing}). Breakdowns are the exception: a
declared sum relates one tag to another in the same context, so a
line's breakdown is never declared to add up to it~\cite{xbrl21}, and a
breakdown can confirm a number but cannot prove one wrong.

\begin{figure}[htbp]
  \small
  \textsf{\textbf{(a)} The label the cash flow statement shows for Apple's tag
  in \cref{fig:xbrl}(a), from the label file \texttt{aapl-20210925\_lab.xml}}
\begin{Verbatim}
<link:label id="lab_us-gaap_ProceedsFromIssuanceOfCommonStock_77a3…_terseLabel_en-US"
  xlink:label="lab_us-gaap_ProceedsFromIssuanceOfCommonStock"
  xlink:role="http://www.xbrl.org/2003/role/terseLabel"
  xlink:type="resource"
  xmlns:xml="http://www.w3.org/XML/1998/namespace"
  xml:lang="en-US">Proceeds from issuance of common stock</link:label>
\end{Verbatim}
  \medskip
  \textsf{\textbf{(b)} The declared sum it takes part in, from the calculation
  file \texttt{aapl-20210925\_cal.xml}: it adds into net cash from financing
  with weight $+1$}
\begin{Verbatim}
<link:calculationArc order="1"
  weight="1.0"
  xlink:arcrole="http://www.xbrl.org/2003/arcrole/summation-item"
  xlink:from="loc_us-gaap_NetCashProvidedByUsedInFinancingActivities_5903…"
  xlink:to="loc_us-gaap_ProceedsFromIssuanceOfCommonStock_63d0…"
  xlink:type="arc"/>
\end{Verbatim}
  \caption{What a filing says about a tag, as filed in Apple's 10-K for fiscal
    2021. The label and the arithmetic are stated once per tag, in separate
    files of the same filing, not on each fact~\cite{aapl-10k-2021}.}
  \label{fig:xbrl-tag}
\end{figure}

\begin{table}[htbp]
  \centering
  \includegraphics[width=\linewidth]{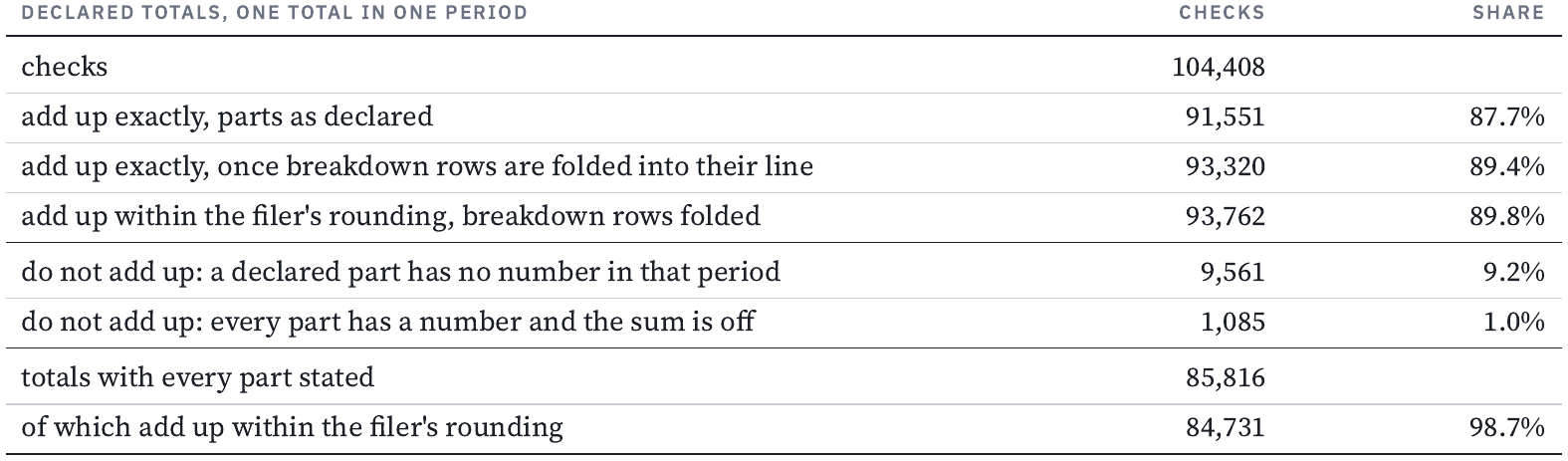}
  \caption{How often a declared total adds up inside its own filing, and
    why the rest do not, over 200 companies drawn at random and their
    2,826 annual reports. The last share is of the totals with
    every part stated.}
  \label{tab:footing}
\end{table}

\subsection{The misses against the vendor, sorted}\label{app:misses}

The vendor of \cref{sec:results-measured} has a different goal from
ours. It standardizes so that companies can be compared with
one another~\cite{fmp-income-statement}, the purpose the SEC gave for
standard tags~\cite{sec-33-9002}; we standardize one company across
time and keep its own lines. The vendor maps each filer's lines onto its
own fields, nets some, gathers what is left into residual buckets,
adjusts per-share figures for later splits, and applies its own
definitions of net income and of payables. Each of those puts a number
in its series that no filing states for the line we anchored, and
\cref{tab:misses} sorts the misses by them first: 258 of the 334 on the
specimen set and 2,683 of the 5,234 on the random sample. Of the rest,
the vendor's figure is elsewhere on our page, in another row of the
same statement, for 30 and 1,380 cells; in 2 and 129 more the anchored
row itself shows the vendor's number in the rendered table, a cell the
table has that the lineage test did not reach. The two rows' labels sort
those: in about half, 17 and 802, they are one line under two tags, a
history the optimization left in two rows, which is our error. It
happens where a line's definition changed and the filer kept both lines
for a while (cash with and without restricted cash, from 2018: the old
tag continues beside the new total, so by our own rule they are two
lines, which the vendor joins), and where the tag changed
across a gap in our filings, so the two lives share no year and the
rename has only names to go on. In the other half the row
that has the number is a different line that happens to carry it, or
one the vendor prefers (total equity with and without noncontrolling
interests, cash with and without restricted cash). What is neither the
vendor's standardization nor explained is 44 cells on the specimen set
and 1,042 on the random sample; with the split lines, the identity
errors the comparison finds or cannot rule out come to 61 cells, 0.6\%,
and 1,844, 3.1\%. These are agreement rates, not accuracy: two series
can agree on a wrong number, and a cell that agrees by coincidence is
not found by examining disagreements. What the comparison does
establish is the coverage gain, and that every cell where the vendor's
standardization departs from the filings counts against us. The method
and the tag alone are scored
against the same series with the same anchors, so the rows of
\cref{tab:fmp} can be compared; the absolute level is set by the
vendor's standardization, not by us.

\begin{table}[htbp]
  \centering
  \includegraphics[width=\linewidth]{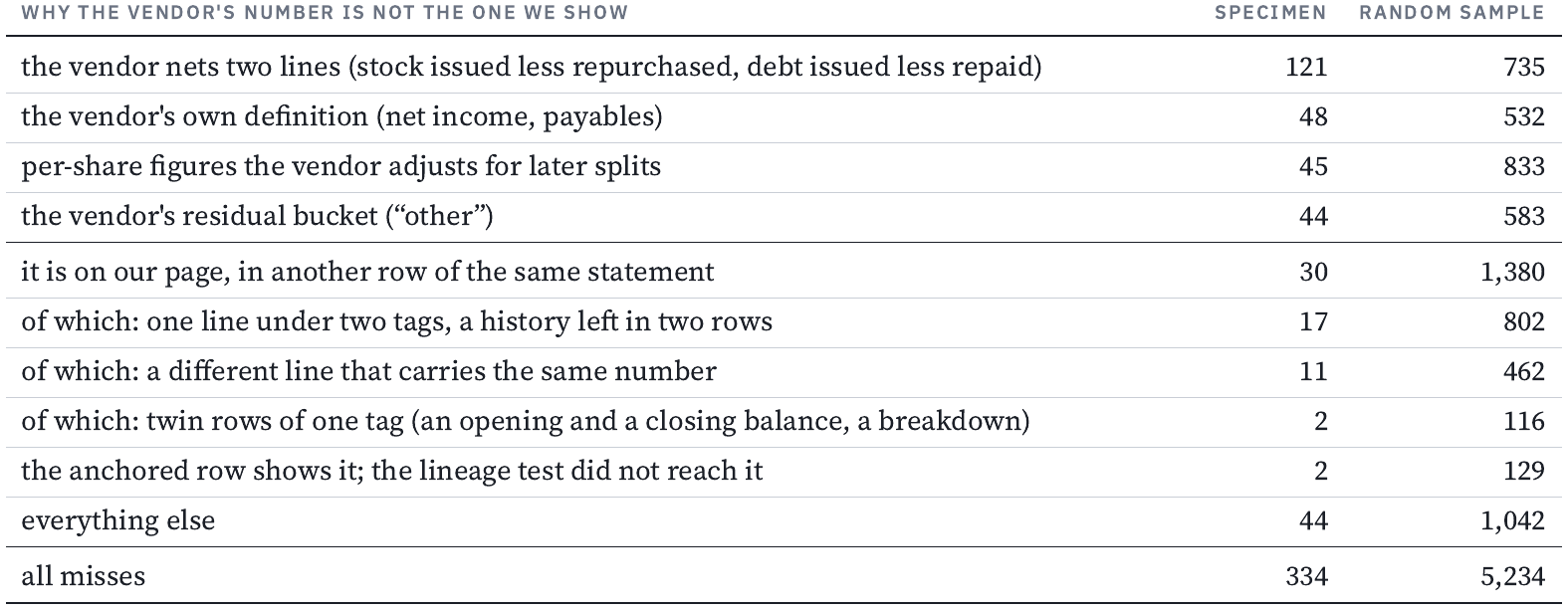}
  \caption{The method's misses against the vendor, sorted: first by what
    the vendor did to the field, then by whether the vendor's number is
    elsewhere on our page, and what the two rows are.}
  \label{tab:misses}
\end{table}

\subsection{The layout differences, and the checks not run}\label{app:layout}

Of the two specimen statements that differ between the newest
filing's layout and the oldest presenter's over fifteen years
(\cref{sec:results-measured}), GE's cash flow statement differs by one
row: the older layout also lists under a segment a figure the newest
shows once (share repurchases under ``GE''), with the same numbers in
both rows, and the checks that use that row move with it; JPMorgan's
balance sheet differs only in which checks run, as do 20 of the random
sample's 33. The other 13 were rendered under both layouts and every
differing row read. In 11 the same numbers sit under a differently
named row, one of them a total the newest layout shows twice; in 2,
both at one company, a merged line that the newest layout derives as
the sum of its parts appears in the older layout as the parts alone.
The script that reads the renders marks four rows as changed values:
the total shown twice, a dividend the newest layout also lists under
the operating partnership, and the two derived rows; in each, wherever
both layouts show a number for the row, it is the same number. An
earlier run of the check found six more, each a row whose lineage
joined a line stated at instants to one stated over periods; under one
layout such a row showed the balances and under the other the flows.
Three were wrong renames and three were one line a filer had tagged
both ways (Appendix~\ref{app:revisions}).

The checks that cannot be run record why. On the random sample,
14,808 of 19,668 are years for which no filing declares arithmetic for
the line, spread across the window; 2,496 have a part with no number
in that column; 369 are the first column of a history.